\documentclass{article}

\usepackage[T1]{fontenc}
\usepackage[letterpaper,margin=1in]{geometry}
\usepackage{amsmath, amsthm, amssymb, mathtools, mleftright}
\usepackage[shortlabels]{enumitem}
\usepackage[colorlinks=true]{hyperref}
\mleftright
\usepackage{tikz-cd}

\newtheorem{thm}{Theorem}[section]
\newtheorem{lem}[thm]{Lemma}
\newtheorem{prop}[thm]{Proposition}
\newtheorem{cor}[thm]{Corollary}
\theoremstyle{definition}
\newtheorem{nota}[thm]{Notation}
\newtheorem{defn}[thm]{Definition}
\newtheorem{qn}[thm]{Question}
\newtheorem{examp}[thm]{Example}
\theoremstyle{remark}
\newtheorem{rem}[thm]{Remark}

\DeclareMathOperator{\tr}{tr}
\DeclareMathOperator{\id}{id}
\DeclareMathOperator{\gr}{gr}
\DeclareMathOperator{\ch}{char}
\DeclareMathOperator{\coef}{coef}
\DeclareMathOperator{\disc}{disc}
\DeclareMathOperator{\inv}{inv}
\DeclareMathOperator{\res}{res}
\DeclareMathOperator{\supp}{supp}
\let\div\undefined \DeclareMathOperator{\div}{div}
\let\th\undefined \DeclareMathOperator{\th}{th}

\newcommand{\alg}{\mathrm{alg}}
\newcommand{\ur}{\mathrm{ur}}
\newcommand{\ram}{\mathrm{ram}}
\newcommand{\cyc}{\mathrm{cyc}}
\newcommand{\nonmax}{\mathrm{nonmax}}
\newcommand{\sep}{\mathrm{sep}}
\newcommand{\zoneI}{\mathrm{I}}
\newcommand{\zoneII}{\mathrm{II}}
\newcommand{\zoneIII}{\mathrm{III}}
\newcommand{\SR}{\mathrm{SR}}
\DeclareMathOperator{\Res}{res}
\DeclareMathOperator{\Hom}{Hom}
\DeclareMathOperator{\Disc}{Disc}
\DeclareMathOperator{\Stab}{Stab}
\DeclareMathOperator{\Hol}{Hol}
\DeclareMathOperator{\Cl}{Cl}
\DeclareMathOperator{\Sel}{Sel}
\DeclareMathOperator{\Sym}{Sym}
\DeclareMathOperator{\Aut}{Aut}
\DeclareMathOperator{\Gal}{Gal}
\DeclareMathOperator{\End}{End}
\newcommand{\GL}{\mathrm{GL}}
\newcommand{\SL}{\mathrm{SL}}
\newcommand{\PGL}{\mathrm{PGL}}
\newcommand{\PSL}{\mathrm{PSL}}
\newcommand{\SO}{\mathrm{SO}}
\renewcommand{\AA}{\mathbb{A}}
\newcommand{\CC}{\mathbb{C}}
\newcommand{\FF}{\mathbb{F}}
\newcommand{\GG}{\mathbb{G}}
\newcommand{\NN}{\mathbb{N}}
\newcommand{\QQ}{\mathbb{Q}}
\newcommand{\PP}{\mathbb{P}}
\newcommand{\RR}{\mathbb{R}}
\newcommand{\ZZ}{\mathbb{Z}}
\renewcommand{\aa}{\mathfrak{a}}
\newcommand{\cc}{\mathfrak{c}}
\newcommand{\mm}{\mathfrak{m}}
\newcommand{\nn}{\mathfrak{n}}
\newcommand{\pp}{\mathfrak{p}}
\renewcommand{\tt}{\mathfrak{t}}
\newcommand{\A}{\mathcal{A}}
\newcommand{\B}{\mathcal{B}}
\newcommand{\C}{\mathcal{C}}
\newcommand{\F}{\mathcal{F}}
\newcommand{\G}{\mathcal{G}}
\renewcommand{\L}{\mathcal{L}}
\newcommand{\OO}{\mathcal{O}}
\newcommand{\V}{\mathcal{V}}
\newcommand{\X}{\mathfrak{X}}
\newcommand{\1}{\mathbf{1}}
\newcommand{\bs}{\backslash}
\newcommand{\cross}{\times}
\newcommand{\tensor}{\otimes}
\renewcommand{\to}{\mathop{\rightarrow}\limits}
\newcommand{\longto}{\mathop{\longrightarrow}\limits}
\newcommand{\textand}{\quad \text{and} \quad}
\newcommand{\size}[1]{\lvert #1 \rvert}
\newcommand{\Size}[1]{\left\lvert #1 \right\rvert}
\newcommand{\floor}[1]{\left\lfloor #1 \right\rfloor}
\newcommand{\ceil}[1]{\left\lceil #1 \right\rceil}
\newcommand{\intsec}{\cap}
\newcommand{\union}{\cup}
\newcommand{\isom}{\cong}
\newcommand{\<}{\left\langle}
\renewcommand{\>}{\right\rangle}
\renewcommand{\(}{\left(}
\renewcommand{\)}{\right)}
\newcommand{\ignore}[1]{}

\renewcommand{\epsilon}{\varepsilon}

\newcommand{\bbq}[8]{
  \begin{minipage}{0.18\linewidth}
    \centering
    \begin{tikzpicture}[
      x=1.9em,
      y=1.7em,
      baseline=(current bounding box.center),
      every node/.style={inner sep=1pt}
      ]
      \node (a1) at (0,0) {$#1$};
      \node (a2) at (2,0) {$#2$};
      \node (a3) at (0,-2) {$#3$};
      \node (a4) at (2,-2) {$#4$};
      \node (a5) at (1,1) {$#5$};
      \node (a6) at (3,1) {$#6$};
      \node (a7) at (1,-1) {$#7$};
      \node (a8) at (3,-1) {$#8$};
      
      \draw (a5)--(a6) (a5)--(a7) (a6)--(a8)
      (a1)--(a5) (a1)--(a2) (a1)--(a3)
      (a2)--(a6) (a2)--(a4) (a7)--(a8)
      (a3)--(a4) (a3)--(a7) (a4)--(a8);
    \end{tikzpicture}
  \end{minipage}
}

\newcommand{\col}[1]{\begin{tabular}{@{}c@{}}#1\end{tabular}}

\newcommand\blfootnote[1]{%
  \begingroup
  \renewcommand\thefootnote{}\footnote{#1}%
  \addtocounter{footnote}{-1}%
  \endgroup
}

\numberwithin{equation}{section}

\begin{document}
\title{Composed varieties\\%
and reflection theorems
of Ohno--Nakagawa type%
}
\author{Evan M{.} O'Dorney\\%
  \em \small Center for Communications Research, Princeton, NJ%
}
\maketitle
% For MathJax for arXiv, etc.
\begin{abstract}
The Ohno-Nakagawa (O-N) reflection theorem is an unexpectedly simple identity relating the number of $\mathrm{GL}_2 \mathbb{Z}$-classes of binary cubic forms (equivalently, cubic rings) of two different discriminants $D$, $-27D$; it generalizes cubic reciprocity and the Scholz reflection theorem. In this paper, we present a new approach to proving and generalizing this theorem using Fourier analysis on the adelic cohomology $H^1(\mathbb{A}_K, M)$ of a finite Galois module, modeled after the celebrated Fourier analysis on $\mathbb{A}_K$ used in Tate's thesis. This method produces reflection theorems of O-N type from identities of lattice counts over local fields and yields new reflection theorems of O-N type for cubic forms and rings over arbitrary number fields, quartic forms and rings, and also quadratic forms counting by a peculiar invariant, the \emph{superdiscriminant} $a(b^2 - 4ac)$. O-N-type theorems are valuable both for their intrinsic beauty and for analytic applications, which include discriminant-reducing identities, symmetrization of functional equations, and subconvexity for Shintani zeta functions.
\end{abstract}

\blfootnote{%
  \paragraph{ORCID:} 0000-0002-7958-2060. \textbf{Email:} \url{emo916math@gmail.com}.
  
  
  \paragraph{I have no competing interests.} This research was partly supported by a Graduate Research fellowship from the National Science Foundation (grant \#DGE-1646566).
  
  \paragraph{Data availability statement:} This project has no associated data.
  
  \paragraph{AI statement:} ChatGPT 5.6 Sol was used as an aid in copyediting, literature search, mathematical error checking, and filling in details of proofs in this paper. All remaining errors are the responsibility of the author.
}

\section{Introduction}
One of the founding results in the theory of prehomogeneous vector spaces is the work of Shintani (\cite{Shintani}; see Taniguchi \cite{TaniguchiOnZeta_pre} for a readable modern introduction), who showed that the generating functions for the class numbers of cubic forms (equivalently, cubic rings) of given discriminant
\begin{equation}\label{eq:Shin1}
  \zeta^{\pm}(s) = \sum_{n \geq 1} \frac{h(\pm n)}{n^{s}}, \qquad h(n) = \!\!\sum_{\substack{f(x,y) = a x^3 + b x^2 y + c x y^2 + d y^3 \\ \disc f = n\text{, up to } \SL_2 \ZZ}} \frac{1}{\size{\Stab_{\SL_2 \ZZ} f}}
\end{equation}
enjoy a meromorphic continuation to the whole complex plane with simple poles at $s = 1$ and $s = 5/6$, with a functional equation
\begin{equation} \label{eq:Shin_fnl_eq}
  \begin{bmatrix}
    \zeta^+(1-s) \\ \zeta^-(1-s)
  \end{bmatrix}
  = 2^{-1}3^{3 s - 2} \pi^{-4s} \Gamma\(s - 1/6\)\Gamma(s)^2 \Gamma\(s + 1/6\)
  \begin{bmatrix}
    \sin 2\pi s & \sin \pi s \\ 3 \sin \pi s & \sin 2\pi s
  \end{bmatrix}
  \begin{bmatrix}
    \zeta^+_3(s) \\ \zeta^-_3(s)
  \end{bmatrix}
\end{equation}
where $\zeta^{\pm}_3(s)$ are the corresponding zeta functions for the dual lattice of integer-matrix forms:
\begin{equation}\label{eq:Shin3}
  \zeta^{\pm}_{3}(s) = \sum_{n \geq 1} \frac{h_{3}(\pm 27n)}{n^{s}}, \qquad
  h_{3}(n) = \!\!\!\!\sum_{\substack{f(x,y) = a x^3 + 3 b x^2 y + 3 c x y^2 + d y^3 \\ \disc f = n\text{, up to } \SL_2 \ZZ}} \frac{1}{\size{\Stab_{\SL_2 \ZZ} f}}
\end{equation}
Twenty-five years later, Ohno unexpectedly discovered numerically, and Nakagawa then proved, that the $\zeta_3^\pm$ are actually equal, term by term (up to an outlying constant), to the $\zeta^\mp$, yielding an \emph{extra functional equation} or \emph{reflection theorem}:
\begin{thm}[Ohno--Nakagawa \cite{Ohno,Nakagawa}]\label{thm:O-N} We have the identities
  \[
    \zeta_3^+(s) = \zeta^-(s) \textand \zeta_3^-(s) = 3\zeta^+(s).
  \]
  In other words, for every nonzero integer $D$,
  \begin{equation}
    \label{eq:O-N_cubic}
    h_3(-27 D) = \begin{cases}
      3 h(D), & D > 0 \\
      h(D), & D < 0.
    \end{cases}
  \end{equation}
\end{thm}
While a functional equation of the type \eqref{eq:Shin_fnl_eq} holds for any prehomogeneous vector space \cite{Shintani,SatoShintani}, the Ohno--Nakagawa (hereafter O-N) relation \eqref{eq:O-N_cubic} is much more mysterious and specific to a restricted class of algebraic group representations. It generalizes the Scholz reflection theorem as well as cubic reciprocity (see \cite[Examples 6.2--6.3]{OOnARemarkable}). Despite its elementary appearance, all known proofs of Theorem \ref{thm:O-N} \cite{Nakagawa,Marinescu,OOnARemarkable,Gao} use global class field theory.

Reflection theorems of this type are of value for multiple reasons. Their intrinsic beauty---an exact equality where we would expect only random fluctuations---points up certain group representations as worthier of study. They are also useful in assorted applications in analytic number theory. For instance, Ellenberg and Venkatesh \cite{EV} use reflection theorems of Scholz type to prove upper bounds on $\ell$-torsion in class groups of number fields, while Mih\u ailescu \cite{MihCat2} uses Leopoldt's
reflection theorem \cite{Leopoldt} to simplify a step of his monumental proof of the Catalan conjecture (that $8$ and $9$ are the only consecutive perfect powers).  A very general reflection theorem for Arakelov class
groups is due to Gras \cite{Gras}.

To date, we have found three families of composed varieties that admit \emph{natural duality}, a remarkable phenomenon leading to a reflection theorem of O-N type. We conjecture that the third of these can be generalized further (see Table \ref{tab:composed}). These composed varieties are geometric orbits cut out by an invariant $I$ of an algebraic group $\Gamma$ acting on a finite-dimensional representation $V$. Surprisingly, the group $\Gamma$ does not have to be reductive (see the first row of the table), nor does the representation $V$ have to be prehomogeneous (see the last two rows). There is a vast sea of representations beyond the Sato-Kimura classification of prehomogeneous vector spaces, and we wonder if many additional examples of natural duality are waiting to be found.

\begin{table}[htbp]
  \centering
  \begin{tabular}{c@{\ }c@{\ }c@{\ }c@{\ }c@{\ }c}
    Acting group $ \Gamma $ & Representation $V$ & Invariant $I$ & \col{ Objects\\parametrized} &
    \col{ Point\\stabilizer\\$ M $} & See \\ \hline
    $ \displaystyle \left\{\left[\begin{array}{cc}
      \lambda & 0            \\
      t       & \lambda^{-2}
    \end{array}\right]\right\}^{\strut} \subset \GL_2 $ & \col{ Quadratic forms, \\ $ \Sym^2(2) $} & \col{Superdiscriminant\\$a(b^2 - 4ac)$} & ? & $ C_2 $ & \col{Theorems\\ \ref{thm:O-N_quad_Z}, \ref{thm:O-N_quad}} \\
    $ \SL_2 $ & \col{ Cubic forms, \\ $ \Sym^3(2) $} & Discriminant & \col{Cubic rings / \\ $ 3 $-torsion in \\ quadratic rings} &  $ C_3 $\vspace{1ex} & \col{Theorems\\\ref{thm:O-N}, \ref{thm:O-N_traced_intro}, \ref{thm:O-N_traced}} \\
    $ \SL_3 $ & \col{ Pairs of ternary \\ quadratic forms, \\ $ \Sym^2(3)^{\oplus 2} $} & Cubic resolvent & \col{Quartic rings / \\ $ 2 $-torsion in \\ cubic rings} & $ C_2 \cross C_2 $ & \cite{OQuartic}\\ \col{$ \SL_n $ ($n \geq 3$ odd)\\(conjectural)} & \col{ Pairs of $n$-ary \\ quadratic forms, \\ $ \Sym^2(n)^{\oplus 2} $} & $n$-ic resolvent & \col{ $ 2 $-torsion in \\ $n$-ic rings } & $ C_2^{n-1} $ & \col{\cite{OQuartic},\\Conjecture\\1.10}
  \end{tabular}
  \caption{Composed varieties admitting a (proven or plausible) reflection identity of Ohno--Nakagawa type}
  \label{tab:composed}
\end{table}

The first row of the table is at once the simplest and the most unexpected representation to admit an O-N-type identity. Quadratic forms under the action of a maximal unipotent subgroup $\GG_a \subset \SL_2$, or what is the same thing, inhomogeneous quadratics $f(x) = a x^2 + b x + c$ under translations $x \mapsto x + t$, admit two algebraically independent invariants: the discriminant $D = b^2 - 4ac$ and the leading coefficient $a$. The natural generating function for these is a two-variable Dirichlet series $\zeta(s,w) = \sum_f \size{a(f)}^s \size{D(f)}^w$, which was studied by Shintani \cite{ShintaniQuadratic} (similar functions appear in \cite{SiegelFE} and \cite[\textsection2.3]{Wen}). However, it is not obvious from their methods that anything extraordinary should arise by specializing to the line $s = w$, that is, counting by the product $a D$ of the two fundamental invariants, which we call the ``superdiscriminant.'' (We have tried other combinations $a^u D^v$ and found no further reflection theorems.) This line \emph{does} appear as the fixed line of the Weyl reflection $(s,w) \mapsto (w,s)$, one of the $12$ functional equations of the closely related double Dirichlet series for the root system $A_2$ \cite{DiamantisGoldfeldShintani,KTW19,HiramotoShintani}, but at a cost of removing the interesting behavior at the prime $2$; see Remark \ref{rem:quad} for details.

\begin{thm}[\textbf{``Quadratic O-N''}] \label{thm:O-N_quad_Z}
  If $n$ is a nonzero integer, let $q(n)$ be the number of integer quadratic polynomials $f(x) = a x^2 + b x + c$ with
  \[
  I(f) = a(b^2 - 4a c) = n,
  \]
  up to the transformation $x \mapsto x + t$ ($t \in \ZZ$). Let $q_2(n)$, $q^+(n)$, and $q_2^+(n)$, respectively, be the number of these $f$ such that $2|b$ (for $q_2$), such that the roots of $f$ are real (for $q^+$), or which satisfy both conditions (for $q_2^+$). Then for all nonzero integers $n$,
  \begin{align*}
    q_2^+(4n) &= q(n) \\
    q_2(4n) &= 2 q^+(n).
  \end{align*}
\end{thm}
More generally, a similar reflection theorem holds over an arbitrary number field of class number $1$; see Theorem \ref{thm:O-N_quad} below.

The second row of Table \ref{tab:composed} is the situation of Ohno--Nakagawa. In this paper, we generalize it to an arbitrary base number field, resolving a conjecture of Dioses \cite[Conjecture 1.1]{Dioses}.
While over $\ZZ$ there is a natural duality between $1111$-forms and $1331$-forms, which naturally parametrize cubic rings and $3$-torsion in class groups of quadratic rings respectively (see \cite{B1,B2}), over $\OO_K$ there are additional integral models of cubic forms (see below) with no known corresponding algebraic object to parametrize (besides cubic rings or self-balanced ideals satisfying added conditions at the primes dividing $3$). Still, the reflection theorem holds:

\begin{thm}\label{thm:O-N_traced_intro}
  Let $K$ be a number field, and let $\aa$ and $\tt$ be ideals of $\OO_K$ with $\tt \mid (3)$. Consider
  \[
  \V_{\aa,\tt}(\OO_K) := \{f(x,y) = a x^3 + b x^2 y + c x y^2 + d y^3 : a \in \aa, b \in \tt, c \in \tt\aa^{-1}, d \in \aa^{-2}\},
  \]
  a $4$-dimensional $\OO_K$-lattice with an action of the group
  \[
  \G_\aa \coloneqq \SL(\OO_K \oplus \aa).
  \]
  For nonzero $D \in \tt^3\aa^{-2}$, define the class number
  \[
  h_{\aa,\tt}(D)
  = \sum_{\substack{f \in \G_{\aa} \backslash \V_{\aa,\tt}(\OO_K) \\ \disc f = D}} \frac{1}{\size{\Stab_{\G_\aa} f}},
  \]
  a weighted count of the forms in $\V_{\aa,\tt}(\OO_K)$ having discriminant $D$. Then we have the identity
  \begin{equation} \label{eq:O-N_traced_intro}
    h_{\aa\tt^{-3}, 3 \tt^{-1}}(-27D)
    = \frac{3^{\#\{v|\infty : D \in (K_v^\cross)^2\}}}{N_{K/\QQ}(\tt)} \cdot h_{\aa,\tt}(D).
  \end{equation}
\end{thm}

The special case $K = \QQ$, $\aa = 1$, $\tt = 1$ recovers Ohno-Nakagawa (Theorem \ref{thm:O-N}). It is worth noting that the counterpart of $\V_{\aa,\tt}$ in this theorem is $\V_{\aa\tt^{-3}, 3 \tt^{-1}}$ and not $\aa\V_{\aa, 3 \tt^{-1}} \isom \V_{\aa^{-1}, 3 \tt^{-1}}$, the dual lattice under the $\SL_2$-invariant alternating pairing on $V(K)$,
\[
  \<(a,b,c,d), (a',b',c',d')\> = a d' - \frac{1}{3} b c' + \frac{1}{3} c b' - d a'
\]
which plays a prominent role in Shintani's work (but will not appear here).

The third row of Table \ref{tab:composed} is a remarkable reflection theorem for quartic rings, which is proved by the same means as the preceding. However, the local computations are vastly more lengthy and have recently been published by the author in a book-length memoir \cite{OQuartic}.

In the fourth row of Table \ref{tab:composed}, we speculate that the analogous result holds for pairs of $n$-ary quadratic forms: see \cite[Conjecture 1.10]{OQuartic}. Such forms have been used to study $2$-torsion in $n$-ic rings \cite{W2xnxn,SwamAvg2Tors}. The reflection conjecture seems far from the reach of current methods.

\subsection{Implications}

The results in this paper belong to a long tradition of \emph{reflection theorems} going back to Scholz \cite{ScholzRefl}, who proved that the $3$-torsion subgroups of the class groups of the quadratic fields $\QQ(\sqrt{D})$ and $\QQ(\sqrt{-3D})$ differ in dimension by at most $1$. A
related theorem due to Leopoldt \cite{Leopoldt} relates different components of the $p$-torsion of the class group of a number field containing $\mu_p$ when decomposed under the Galois group of that field. Applications of such reflection theorems are
far-ranging: for instance, Ellenberg and Venkatesh \cite{EV} use reflection theorems of Scholz type to prove upper bounds on $\ell$-torsion in class groups of number fields, while Mih\u ailescu \cite{MihCat2} uses Leopoldt's
generalization to simplify a step of his monumental proof of the Catalan conjecture (that $8$ and $9$ are the only consecutive perfect powers). Through the
years, numerous reflection principles for different generalizations of ideal
class groups have come into print. A very general reflection theorem for Arakelov class
groups is due to Gras \cite{Gras}.

By plugging in fundamental discriminants $D$ into Theorem \ref{thm:O-N}, we get back Scholz's reflection theorem, as Nakagawa points out (\cite{Nakagawa}, Remark 0.9). Likewise, our Theorem \ref{thm:O-N_traced_intro} is a further generalization of Scholz's reflection theorem for number fields (see Ellenberg--Venkatesh \cite{EV}, Lemma 3.3).

While the discriminants $D$, $-27D$ in the O-N theorem are of the same order of magnitude, our methods also produce some \emph{discriminant-reducing identities} (see Section \ref{sec:disc_red}), useful tools in arithmetic statistics (see Bhargava, Taniguchi, and Thorne \cite{BTTImproved}) that reduce troublesome discriminant divisibility conditions to counting problems with lower-order discriminants.

The Shintani zeta functions $\zeta^\pm$ encapsulate the distribution of discriminants of cubic rings and fields, and their simple poles at $1$ and $5/6$ reflect the first- and second-order terms in counting such discriminants \cite{BST-2ndOrd,TT_rc}. Using the O-N relation, the functional equation can be diagonalized into a self-reflective form \cite[p.~1088]{Ohno}. This self-reflective functional equation was used by Hough and Lee \cite{HoughLee21} to obtain subconvexity bounds on Shintani zeta functions. Subconvexity is an invaluable property controlling zeta functions in the critical strip, used when stronger results such as the Riemann hypothesis are unavailable. The results in this paper point toward proofs of subconvexity bounds for several other families of Shintani-type zeta functions.

More broadly, the presence of natural duality on a particular space can be construed as motivation to study it amid the infinite sea of algebraic group representations. In particular, we wonder if the superdiscriminant $a D$ of an inhomogeneous quadratic has a deeper meaning to fill in the blank space in Table \ref{tab:composed}.

\subsection{Methods}

Several proofs of O-N are now in print \cite{Nakagawa,Marinescu,OOnARemarkable,Gao}, all of which consist of two main steps:
\begin{itemize}
  \item A ``local'' step to count the cubic rings in each cubic field, equivalently the $\GL_2(\ZZ)$-orbits in each $\GL_2(\QQ)$-orbit, and put the result in a usable form;
  \item A ``global'' step that uses global class field theory to understand the cubic fields that can have specified local orbit counts.
\end{itemize}
In this paper, the distinction between these steps is formalized and clarified. The global step, in particular, is streamlined into a \emph{local-global reflection engine} (Theorem \ref{thm:main_compose}) that can be reused to prove additional O-N-like reflection theorems over number fields $K$, given an appropriate \emph{local reflection theorem} over each completion of $K$.

The local-global reflection engine is proved by Fourier analysis on adelic cohomology: each reflection identity ultimately boils down to Poisson summation on a cohomology group $H^1(\AA_K, M)$, such as $H^1(\AA_K, \mu_n) \isom \AA_K^\cross/(\AA_K^\cross)^n$. This is a novel method inspired by Tate's thesis \cite{Tate_thesis}, which uses Poisson summation on the additive group of the adeles to prove functional equations for Hecke zeta functions. More recent precedent can be found in the work of Frei, Loughran, and Newton \cite{FLN18}, who use Poisson summation on the (multiplicative) idele group to count abelian extensions subject to infinitely many local conditions, and of Taniguchi--Thorne \cite{TT_oexp}, who compute Fourier duals of splitting-type indicators for binary forms over finite fields, finding cancellations that are just as mysterious as, but possibly unrelated to, the apparent coincidences that go into O-N. More recently, the author and Alberts \cite{BAEO_with_Corrigendum} extended the method of Poisson summation on adelic $H^1$ to count coclasses asymptotically by quite general discriminant-like invariants. Further refinements of this method underlie the recent counting results for nilpotent number field extensions by Alberts--Bucur \cite{AlbertsBucurMDS} and Alberts--Lemke Oliver--Wang--Wood \cite{ALOWW}.

The local-global reflection engine is applicable on a quite general class of representations of algebraic groups, which we call \emph{composed varieties.} The term ``composed'' refers to the presence of a composition law on the orbits, which corresponds to the group law of a Galois cohomology group $H^1(K,M)$. Our guiding example is the variety of binary cubic forms of discriminant $D$ with its action of $\Gamma = \SL_2$, where rational orbits correspond to $H^1(K,M)$ for an appropriate module $M$ of order $3$. This is one of Bhargava's \emph{higher composition laws} \cite{B2}.

\subsection{Outline of the paper}

In Section \ref{sec:composed}, we introduce composed varieties and prove the local-to-global reflection engine (Theorem \ref{thm:main_compose}) by Fourier analysis of the local and global Tate pairings. The following section develops the structure of the cohomology groups $H^1(K,M)$ needed for the local calculations. In Section \ref{sec:quadratic}, we prove quadratic O-N (Theorem \ref{thm:O-N_quad}), of which Theorem \ref{thm:O-N_quad_Z} is a special case. In Section \ref{sec:cubic}, we prove cubic O-N (Theorem \ref{thm:O-N_traced_intro}) and its corollaries. Section \ref{sec:non-natural} applies the same engine to more general local weightings, including discriminant-reduction identities and binary cubic forms over $\ZZ[1/N]$. We conclude with a function-field application.

Except in the final section, the global field is a number field; the local arguments are stated over the corresponding local fields whenever possible.

\subsection{Acknowledgements}
For fruitful discussions, I would like to thank (in no particular order):
Manjul Bhargava, Xiaoheng (Jerry) Wang, Fabian Gundlach, Levent Alp\"oge, Melanie Matchett Wood, Kiran Kedlaya, Alina Bucur, Benedict Gross, Sameera Vemulapalli, Frank Calegari, Frank Thorne, Brandon Alberts, Peter Sarnak, Jack Thorne, Theresa (Tess) Anderson, and Jennifer O'Dorney.

\subsection{Notation}\label{sec:notation}
The following conventions will be observed in the remainder of the paper.

We denote by $\NN$ the set of nonnegative integers.

If $P$ is a statement, then
\[
  \1_P = \begin{cases}
    1 & \text{$P$ is true} \\
    0 & \text{$P$ is false}.
  \end{cases}
\]
If $S$ is a set, then $\1_S$ denotes the characteristic function $\1_S(x) = \1_{x \in S}$.

An \emph{algebra} will always be commutative and of finite rank over a field, while a \emph{ring} or \emph{order} will be a finite-rank, torsion-free ring over a Dedekind domain, containing $1$. An order need not be a domain.

If $\OO_K$ is a Dedekind domain, $K$ is its field of fractions, $V$ is an $n$-dimensional vector space over $K$ and $A,B \subseteq V$ are two full-rank lattices, we denote by the \emph{index} $[A:B]$ the unique fractional ideal $\cc$ such that
\[
\Lambda^n A = \cc \Lambda^n B
\]
as lattices in the top exterior power $\Lambda^n V$. Note that if $A \supseteq B$ and 
\[
A/B \cong \OO_K/\cc_1 \oplus \cdots \oplus \OO_K/\cc_r
\]
as $\OO_K$-modules, then
\[
[A:B] = \cc_1\cc_2 \cdots \cc_r.
\]
If $L$ is a $K$-algebra, $\OO \subseteq L$ is an order, and $\aa \subseteq L$ is a fractional ideal, the index $[\OO : \aa]$ is called the \emph{norm} of $\aa$ and will be denoted by $N_\OO(\aa)$ or, when the context is clear, by $N(\aa)$. The norm is multiplicative for $\OO$-ideals that are invertible (i.e.\ locally principal), but not for general $\OO$-ideals.

If $K$ is a field, we denote by $K^\alg$ an algebraic closure of $K$.

If $S$ is a finite set, we let $\Sym(S)$ denote the set of permutations of $S$; thus $S_n = \Sym(\{1,\ldots,n\})$. If $\size{S} = \size{T}$, and if $g \in \Sym(S)$, $h \in \Sym(T)$ are elements, we say that $g$ and $h$ are \emph{conjugate} if there is a bijection between $S$ and $T$ under which they correspond. Likewise when we say that two subgroups $G \subseteq \Sym(S)$, $H \subseteq \Sym(T)$ are conjugate.

We will use the semicolon to separate the coordinates of an element of a product of rings. For instance, in $\ZZ \cross \ZZ$, the nontrivial idempotents are $(1;0)$ and $(0;1)$.

If $n$ is a positive integer, then $\zeta_n$ denotes a primitive $n$th root of unity in $\QQ^\alg$, while $\bar\zeta_n$ denotes the $n$th root of unity
\[
\bar\zeta_n = \(1; \zeta_n; \zeta_n^2; \ldots; \zeta_n^{n-1}\) \in (\QQ^\alg)^n.
\]

Throughout the proofs of the local reflection theorems (Sections \ref{sec:quadratic}--\ref{sec:cubic}), we will fix a local field $K$, its valuation $v = v_K$ (which has value group $\ZZ$), its residue field $k_K$ of order $q$, and a uniformizer $\pi = \pi_K$. Often the choice of uniformizer is immaterial, but not always (for example, a tamely ramified cubic extension will be described without loss of generality as $R = K[\sqrt[3]{\pi}]$.) We put
\[
  e = \begin{cases}
    v_K(2) & \text{in the quadratic case} \\
    v_K(3) & \text{in the cubic case};
  \end{cases}
\]
note that $e$ is either the absolute ramification index of $K$ or $0$. We fix an extension of $v_K$ to $K^\alg$. We let $\mm_K$ denote the maximal ideal of $K$, and likewise let $\mm_{K^\alg}$ be the maximal ideal of the corresponding valuation ring $\OO_{K^\alg}$; note that $\mm_{K^\alg}$ is not finitely generated. We also allow $v = v_K$ to be applied to elements of $K^\alg$, with this fixed extension, scaled so that its restriction to $K$ has value group $\ZZ$. We use the absolute value bars $\size{\bullet}$ for the corresponding multiplicative absolute value, whose normalization will be left undetermined.

If $R/\OO_K$ is a ring, we denote by $R^{N=1}$ the subgroup of units of norm $1$ down to $K$. Implicitly, the group operation is multiplication, so that $R^{N=1}[n]$, for instance, denotes the $n$th roots of unity of norm $1$ (of which there may be more than $n$ if $R$ is not a domain).

\section{Composed varieties}
\label{sec:composed}
It has long been noted that orbits of certain algebraic group actions on varieties over a field $K$ parametrize rings of low rank over $K$, which can sometimes also be identified with the cohomology of small Galois modules over $K$. The aim of this section is to explain all this in a level of generality suitable for our applications.

\begin{defn}
  Let $K$ be a field. A \emph{composed variety} over $K$ consists of the data $\bigl(\Gamma, V, X, X_0\bigr)$, where
  \begin{itemize}
    \item $V$ is a finite-dimensional, faithful representation of a connected linear algebraic group $\Gamma$ over $K$;
    \item $X \subseteq V$ is a closed subvariety whose points $X(K^{\alg})$ over the algebraic closure $K^{\alg}$ form a single $\Gamma(K^{\alg})$-orbit, such that the stabilizer $M = \Stab_\Gamma(x_0)$ is finite, abelian, and \'etale of order not divisible by $\ch K$;
    \item $X_0 = \Gamma(K)x_0$ is an orbit of a $K$-rational point $x_0 \in X$.
  \end{itemize}
\end{defn}

\begin{rem}
  We have endeavored to make this definition as general as possible. We require that $X$ be closed in order to get finite class numbers (compare Example \ref{ex:Gm}). That the stabilizer $M$ be abelian is essential at several steps. To lift the hypotheses that $M$ be \'etale and of order prime to $\ch K$ entails replacing Galois cohomology by flat cohomology, an interesting extension that we do not pursue here although the groundwork was laid by \v{C}esnavi\v{c}ius \cite{Ces_Poitou,Ces_Topo}.
\end{rem}

The term \emph{composed} is derived from Gauss composition of binary quadratic forms and the higher composition laws of the work of Bhargava and his successors, from which we derive many of our examples. In particular, the $\Gamma(K)$-orbits often admit a group structure, which can be explained as a cohomology group:

\begin{prop}~\label{prop:composed} Let $\bigl(\Gamma, V, X, X_0\bigr)$ be a composed variety. Then:
  \begin{enumerate}[$($a$)$]
    \item\label{comp:psi} There is a natural injection
    \[
      \psi : \Gamma(K)\backslash X(K) \hookrightarrow H^1(K,M)
    \]
    by which the orbits $\Gamma(K)\backslash X(K)$ parametrize some subset of the Galois cohomology group $H^1(K, M)$;
    \item\label{comp:0} We have $\psi(X_0) = 1$;
    \item\label{comp:stab} The $\Gamma(K^{\alg})$-stabilizer of every $x \in X(K)$ is canonically isomorphic to $M$ as a Galois module;
    \item\label{comp:chg_bspnt} If $X_0$ is replaced by a different base orbit $X_1 = \Gamma(K) x_1$, then the associated parametrization $\psi_{X_1} : \Gamma(K)\backslash X(K) \to H^1(K,M)$ differs only by translation:
    \[
      \psi_{X_1}(Y) = \psi(Y)\psi(X_1)^{-1},
    \]
    under the isomorphism between the stabilizers $M$ established in the previous part.
    \item\label{comp:base_chg} If $K'/K$ is a field extension (not necessarily algebraic), then $\bigl(\Gamma', V', X', X_0'\bigr)$ is also a composed variety, where $\Gamma', V', X'$ are formed by extension of scalars and $X_0' = \Gamma(K') x_0$. The corresponding parametrizations $\psi, \psi'$ are related by the commutative diagram
    \[
      \begin{tikzcd}
        \Gamma(K)\bs X(K) \arrow[r, hook, "\psi"] \arrow[d]
          & H^1(K,M) \arrow[d, "\Res"] \\
        \Gamma'(K')\bs X'(K') \arrow[r, hook, "\psi'"]
          & H^1(K',M)
      \end{tikzcd}
    \]
  \end{enumerate}  
\end{prop}
\begin{proof}
  \begin{enumerate}[$($a$)$]
    \item It is well known that, more generally, if $X$ is a homogeneous space for an algebraic group $\Gamma$ with a $K$-rational point $x_0$ and smooth stabilizer $M$, there is a long exact sequence of pointed sets
    \begin{equation} \label{eq:LES_hom}
      M(K) \longrightarrow
      \Gamma(K) \longrightarrow
      X(K) \stackrel{\psi}{\longrightarrow}
      H^1(K, M) \longrightarrow H^1(K, \Gamma).
    \end{equation}
    Here we recall the construction of the connecting map $\psi$ for future reference. Let $x \in X(K)$ be given. Since $X$ is a $\Gamma(K^{\alg})$-orbit, we can find $\gamma \in \Gamma(K^{\alg})$ such that $\gamma(x_0) = x$. There are $|M|$ such $\gamma$, differing by right-multiplication by $M$; since $\ch K \nmid |M|$, these $\gamma$ are defined over the separable closure $K^{\sep}$. Fix one such $\gamma$. For any $g \in \Gal(K^{\sep}/K)$, $g(\gamma)$ also takes $x_0$ to $x$, so the product $\gamma^{-1} \cdot g(\gamma)$
    lies in $\Stab_{\Gamma(K^{\sep})} x_0 = M$. Define a map $\sigma_x : \Gal(K^{\sep}/K) \to M$ by
    \[
    \sigma_x(g) = \gamma^{-1} \cdot g(\gamma).
    \]
    It is routine to verify the following:
    \begin{enumerate}[(i)]
      \item $\sigma_x$ satisfies the cocycle condition $\sigma_x(gh) = \sigma_x(g) \cdot g(\sigma_x(h))$ and hence defines a coclass in $H^1(K,M)$.
      \item Given $x$, the choice of $\gamma$ is unique up to right-multiplication by $m \in M$. If we replace $\gamma$ by $\gamma m$, then $\sigma_x$ changes by multiplying by the coboundary $g \mapsto m^{-1} g(m)$.
      \item If $x$ is replaced by a different orbit representative $\alpha x$, where $\alpha \in \Gamma(K)$, the cocycle $\sigma_x$ is unchanged.
      \item If the basepoint $x_0$ is replaced by a different orbit representative $\alpha x_0$ for some $\alpha \in \Gamma(K)$, the cocycle $\sigma_x$ is unchanged, up to identifying $M$ with $\Stab_{\Gamma(K^{\sep})} (\alpha x_0) = \alpha M \alpha^{-1}$ in the obvious way.
    \end{enumerate}
    So we get a well-defined map
    \begin{align*}
      \psi : \Gamma(K)\backslash X(K) &\to H^1(K,M) \\
      \Gamma(K) x &\mapsto \sigma_x.
    \end{align*}
    It remains to show that $\psi$ is injective. Suppose that $x_1, x_2 \in X(K)$ map to equivalent cocycles $\sigma_{x_1}$, $\sigma_{x_2}$. Let $\gamma_i \in \Gamma(K^{\sep})$ be an associated transformation that maps $x_0$ to $x_i$. By right-multiplying $\gamma_1$ by an element of $M$, as above, we can remove any coboundary discrepancy and assume that $\sigma_{x_1} = \sigma_{x_2}$ as cocycles. That is, for every $g \in \Gal(K^{\sep}/K)$,
    \[
    \gamma_1^{-1} \cdot g(\gamma_1) = \gamma_2^{-1} \cdot g(\gamma_2),
    \]
    which can also be written as
    \[
    g(\gamma_2 \gamma_1^{-1}) = \gamma_2 \gamma_1^{-1}.
    \]
    Thus, $\gamma = \gamma_2 \gamma_1^{-1}$ is Galois stable and hence defined over $K$. We have $\gamma(x_1) = x_2$, establishing that these points lie in the same $\Gamma(K)$-orbit, as desired.
    \item If $x = x_0$, we can take $\gamma = \id \in \Gamma(K^{\sep})$, and it is immediate that $\sigma_x(g) = \id$ for all $g \in \Gal(K^{\sep}/K)$.
    \item If $\gamma(x_0) = x$, then the $\Gamma(K^{\sep})$-stabilizer of $x$ is of course $\gamma M \gamma^{-1}$. We claim that the obvious map
    \begin{align*}
      M \to \gamma M \gamma^{-1} \\
      \mu \mapsto \gamma \mu \gamma^{-1}
    \end{align*}
    is an isomorphism of Galois modules. We compute, for $g \in \Gal(K^{\sep}/K)$,
    \[
    g\(\gamma \mu \gamma^{-1}\) = g(\gamma) g(\mu) g(\gamma)^{-1}
    = \gamma \sigma_x(g) g(\mu) \sigma_x(g)^{-1} \gamma^{-1}
    = \gamma g(\mu) \gamma^{-1},
    \]
    since $M$ is abelian. This establishes the isomorphism. By a similar calculation, the identification of stabilizers is independent of $\gamma$ and is thus canonical.
  \end{enumerate}
  Parts \ref{comp:chg_bspnt} and \ref{comp:base_chg} are routine.
\end{proof}

While $\psi$ is always injective, it need not be surjective, as we will see by examples in the following section, motivating the following definitions.
\begin{defn}~\label{defn:fullHasse}
  \begin{enumerate}[$($a$)$]
    \item A composed variety $X$ is \emph{full} if its cohomology parametrization $\psi$ is surjective, that is, $X$ includes a $\Gamma(K)$-orbit for every coclass $\alpha \in H^1(K,M)$.
    \item If $K$ is a global field, a composed variety is \emph{Hasse} if it satisfies the following version of the Hasse local-global principle: For $\alpha \in H^1(K,M)$, if the localization $\alpha_v \in H^1(K_v,M)$ at each place $v$ lies in the image of the local parametrization
    \[
    \psi_v : \Gamma(K_v)\backslash X(K_v) \to H^1(K_v, M),
    \]
    then $\alpha$ also lies in the image of the global parametrization $\psi$.
  \end{enumerate}
\end{defn}

\begin{prop}\label{prop:fin}
  A composed variety over a finite field is full.
\end{prop}
\begin{proof}
  This is evident from the exact sequence \eqref{eq:LES_hom}, since $H^1(K, \Gamma) = 0$ by Lang's theorem.
\end{proof}

\subsection{Localization of global orbit counts}
\subsubsection{Notations}
In this section, we will be interested in composed varieties over local fields $K_v$ and global fields $K$, which we call \emph{local} and \emph{global composed varieties}, respectively. We write $\V^K$ for the set of places of $K$ and $\V^K_{\mathrm f}$ for the set of nonarchimedean places.

We will often consider a local composed variety $(\Gamma, V, X, X_0)$ equipped with an $\OO_v$-lattice $\Lambda_v$. Then we let $\Gamma_{\Lambda_v}$ be the subgroup of $\Gamma(K_v)$ fixing $\Lambda_v$. Often $\Lambda_v$ is chosen specially, for instance so that $\Gamma_{\Lambda_v}$ is a maximal compact subgroup of $\Gamma(K_v)$, but we do not require this. If $v$ is an archimedean place, these notions do not make sense and we arbitrarily set $\Lambda_v = V(K_v)$, so $\Gamma_{\Lambda_v} = \Gamma(K_v)$.

Likewise, we will often consider a global composed variety $(\Gamma, V, X, X_0)$ equipped with an \emph{adelic lattice} $\Lambda = \prod_v \Lambda_v$. Here $\Lambda_v \subseteq V(K_v)$ is an $\OO_v$-lattice at each nonarchimedean place and $\Lambda_v = V(K_v)$ at each archimedean place, and, for some (equivalently, every) basis $x_1,\ldots,x_n$ of $V$, we have $\Lambda_v = \OO_v\<x_1,\ldots,x_n\>$ for all but finitely many $v$. Adelic lattices are an essential tool for understanding Shintani zeta functions over global fields \cite{WrightAdelicZeta, DW2, TaniguchiOnZeta_pre}. We let $\Gamma_{\Lambda}$ be the subgroup of $\Gamma(\AA_K)$ preserving $\Lambda$, and write $\Gamma_{\Lambda}(K) = \Gamma_{\Lambda} \intsec \Gamma(K)$.

\subsubsection{Spreading out}

First, we show that the behavior of a global composed variety at most places $v$ is well controlled.

\begin{lem}\label{lem:spread}
  Let $(\Gamma, V, X, X_0)$ be a composed variety over a global field $K$, and let $\Lambda = \prod_v \Lambda_v$ be an adelic lattice in $V$. For almost all places $v$ of $K$, the following occurs:
  \begin{enumerate}[$($a$)$]
    \item\label{spread:good_red} $\Gamma$, $X$, and $M$ have good reduction at $v$, and their reductions have the same number of geometric irreducible components as the global varieties;
    \item\label{spread:H1ur} There are $|H^0(K_v, M)|$-many $\Gamma(K_v)$-orbits of $X(K_v)$ that intersect $\Lambda_v$, namely those parametrized by the \emph{unramified} cohomology $H^1_\ur(K_v, M)$.
    \item\label{spread:Gv} Each such orbit intersects $\Lambda_v$ in a single orbit of the subgroup $\Gamma_{\Lambda_v}$.
  \end{enumerate}
\end{lem}

\begin{proof}
  Choose a basis $(x_1, \ldots, x_n)$ of $V$. We first claim that, for almost all places $v$, the following hold:
  \begin{enumerate}[(i)]
    \item\label{sprc:basis} $(x_1, \ldots, x_n)$ is an $\OO_v$-basis of $\Lambda_v$;
    \item\label{sprc:imm} the closed immersion $\Gamma \hookrightarrow \GL(V) \isom \GL_n$, the action of $\Gamma$ on $X$, the point $x_0$, and the stabilizer $M$ all extend over $\OO_v$; the resulting models of $\Gamma$ and $X$ are smooth and geometrically irreducible, and the model of $M$ is finite \'etale of degree $|M|$;
    \item\label{sprc:orb} the orbit map
    \[
      \rho : \Gamma \longrightarrow X, \qquad \gamma \longmapsto \gamma x_0,
    \]
    extends to a finite \'etale surjection over $\OO_v$, necessarily an $M$-torsor.
  \end{enumerate}
  
  This is shown by spreading out the schemes and morphisms involved, since the conditions of being smooth, \'etale, or finite \'etale are open conditions. As to the irreducibility condition in \ref{sprc:imm}, it is well known that, for $Y$ of finite presentation over $S$, the locus of $s \in S$ for which $Y_s$ is geometrically irreducible is constructible \cite[Th\'eor\`eme 9.7.7]{EGA_IV}, so since the generic fiber is geometrically irreducible, so are almost all the special fibers. This proves \ref{spread:good_red}.
  
  We now fix a place $v$ satisfying \ref{sprc:basis}--\ref{sprc:orb}. Note that for the chosen integral model of $\Gamma$,
  \begin{equation} \label{eq:GLv}
    \Gamma_{\Lambda_v} = \Gamma(\OO_v) = \Gamma(K_v) \intsec \GL_n(\OO_v).
  \end{equation}

  Since the ring $\OO_v^\ur$ is strictly Henselian, every finite \'etale torsor over $\OO_v^\ur$ is trivial. Thus, for every $x \in X(\OO_v^\ur)$, the fiber of $\rho$ at $x$ has an $\OO_v^\ur$-point, that is, we can find $\gamma_x \in \Gamma(\OO_v^\ur)$ such that
  \[
    \gamma_x x_0 = x.
  \]
  In particular,
  \[
    X(\OO_v^\ur) = \Gamma(\OO_v^\ur)x_0
    \qquad\text{and}\qquad
    X(k_v^{\alg}) = \Gamma(k_v^{\alg})\bar x_0.
  \]
  Since $M$ is finite \'etale over $\OO_v$, its geometric points are defined over $\OO_v^\ur$, so $M$ descends to a Galois module $\bar M$ over $k_v$ and
  \[
    H^1_\ur(K_v, M) \isom H^1(k_v, \bar M).
  \]

  Now decompose
  \[
    X(\OO_v) = X(K_v) \intsec \Lambda_v
  \]
  into $\Gamma(\OO_v)$-orbits. If $x \in X(\OO_v)$, choose $\gamma_x \in \Gamma(\OO_v^\ur)$ with $\gamma_x x_0 = x$. The associated cocycle is unramified, so every such orbit determines a coclass in $H^1_\ur(K_v, M)$.

  Suppose that two $\Gamma(\OO_v)$-orbits represented by $x, y \in X(\OO_v)$ determine the same coclass. By Proposition \ref{prop:composed}, there is a $\gamma \in \Gamma(K_v)$ with $\gamma x = y$. However, the elements of $\Gamma(K_v^\ur)$ taking $x$ to $y$ are precisely
  \[
    \gamma_y \mu \gamma_x^{-1}, \qquad \mu \in M,
  \]
  and each of these lies in $\Gamma(\OO_v^\ur)$. Hence $\gamma \in \Gamma(\OO_v^\ur) \intsec \Gamma(K_v) = \Gamma(\OO_v)$, so $x$ and $y$ were already in the same $\Gamma(\OO_v)$-orbit.
  This proves \ref{spread:Gv} once we know that every unramified coclass occurs.

  Since $\Gamma$ is connected, its good reduction $\bar\Gamma$ is connected, and by Proposition \ref{prop:fin}, the composed variety $\bigl(\bar\Gamma, k_v^n, \bar X, \bar x_0\bigr)$ is full, that is, every coclass in $H^1(k_v, \bar M)$ is represented by a $\bar\Gamma(k_v)$-orbit on $\bar X(k_v)$. Since $X$ is smooth over $\OO_v$, every point of $\bar X(k_v)$ lifts to a point of $X(\OO_v)$. Thus every coclass in $H^1_\ur(K_v, M)$ occurs among the integral orbits.

  We have therefore proved that the $\Gamma(K_v)$-orbits meeting $\Lambda_v$ are exactly those parametrized by $H^1_\ur(K_v, M)$ and that each meets $\Lambda_v$ in a single $\Gamma_{\Lambda_v}$-orbit. Finally, for the finite Galois module $\bar M$ over the finite field $k_v$,
  \[
    |H^1(k_v, \bar M)| = |H^0(k_v, \bar M)| = |H^0(K_v, M)|,
  \]
  so there are exactly $|H^0(K_v, M)|$ such $\Gamma(K_v)$-orbits. This proves \ref{spread:H1ur} and completes the proof. \qedhere
\end{proof}

\subsubsection{Local and global lattice counts}

\begin{defn} \label{def:loc_lat_cnt}
  If $\bigl(\Gamma, V, X, X_0\bigr)$ is a local composed variety with lattice $\Lambda_v$ and $x_v \in X(K_v)$, the \emph{local lattice count} $\ell_v(x_v)$ is the number of translates $\gamma_v\Lambda_v$, $\gamma_v \in \Gamma(K_v)$, that contain $x_v$.
\end{defn}
\begin{rem}
  Equivalently, after mapping $\gamma_v \mapsto \gamma_v^{-1}$, $\ell_v(x_v)$ is the number of left cosets
  \[
    \Gamma_{\Lambda_v}\gamma_v \in \Gamma_{\Lambda_v}\bs\Gamma(K_v)
  \]
  such that $\gamma_v x_v \in \Lambda_v$. Note that $\ell_v(x_v)$ depends only on the orbit $\Gamma(K_v)x_v$.
\end{rem}

\begin{defn} \label{def:glo_lat_cnt}
  If $\bigl(\Gamma, V, X, X_0\bigr)$ is a global composed variety with adelic lattice $\Lambda$ and $x = (x_v)_v \in X(\AA_K)$, the \emph{global lattice count} $\ell(x)$ is the number of translates $\gamma\Lambda$, $\gamma \in \Gamma(\AA_K)$, that contain $x$.
\end{defn}
\begin{rem}
  Equivalently, after mapping $\gamma \mapsto \gamma^{-1}$, $\ell(x)$ is the number of left cosets
  \[
    \Gamma_\Lambda\gamma \in \Gamma_\Lambda\bs\Gamma(\AA_K)
  \]
  such that $\gamma x \in \Lambda$. Note that $\ell(x)$ depends only on the orbit $\Gamma(\AA_K)x$.
\end{rem}

\begin{lem}\label{lem:cl_num_glo_loc}
  For any $x_v \in X(K_v)$, the local lattice count $\ell_v(x_v)$ is finite. For $x = (x_v)_v \in X(\AA_K)$, the global lattice count $\ell(x)$ is finite and is given by
  \begin{equation}\label{eq:cl_num_glo_loc}
    \ell(x) = \prod_{v \in \V^K} \ell_v(x_v).
  \end{equation}
\end{lem}
\begin{proof}
  At an archimedean place, our convention $\Lambda_v = V(K_v)$ gives $\ell_v(x_v) = 1$, so suppose that $v$ is nonarchimedean. Fix $x_v \in X(K_v)$ and consider the orbit map
  \[
    \rho_{x_v} : \Gamma \longrightarrow X, \qquad \gamma \longmapsto \gamma x_v.
  \]
  It is a finite map, hence proper. Therefore
  \[
    S_v = \{\gamma \in \Gamma(K_v) : \gamma x_v \in \Lambda_v\}
       = \rho_{x_v}^{-1}\bigl(X(K_v) \intsec \Lambda_v\bigr)
  \]
  is compact. Since $\Gamma_{\Lambda_v}$ is open in $\Gamma(K_v)$, the compact set $S_v$ meets only finitely many left cosets of $\Gamma_{\Lambda_v}$. Thus the local lattice count $\ell_v(x_v)$ is finite.

  Now let $x \in X(\AA_K)$. For almost all $v$, we have $x_v \in \Lambda_v$, and Lemma \ref{lem:spread}\ref{spread:Gv} gives $\Gamma(K_v)x_v \intsec \Lambda_v = \Gamma_{\Lambda_v}x_v$, implying $\ell_v(x_v) = 1$ for almost all $v$. Finally, both the adelic group and the adelic lattice are restricted products of their local counterparts, so choosing an adelic translate containing $x$ is equivalent to choosing, independently at each place, one of the finitely many local translates containing $x_v$. This proves the finiteness of $\ell(x)$ as well as \eqref{eq:cl_num_glo_loc}.
\end{proof}

For $\alpha \in H^1(K, M)$, define
\[
  \ell(\alpha) =
  \begin{cases*}
    \ell(x) & if $x \in X(K)$ corresponds to the coclass $\alpha$, \\
    0 & if no such $x$ exists,
  \end{cases*}
\]
and define $\ell_v(\alpha_v)$ analogously for $\alpha_v \in H^1(K_v, M)$.

\begin{lem}\label{lem:Hasse}
  If $X$ is Hasse, then for all $\alpha \in H^1(K, M)$,
  \[
    \ell(\alpha) = \prod_{v \in \V^K} \ell_v(\alpha_v).
  \]
\end{lem}
\begin{proof}
  If the right-hand side is zero, then some localization $\alpha_v$ is not represented by an orbit in $X(K_v)$, so $\alpha$ cannot be represented by a global orbit and the left-hand side is also zero. If the right-hand side is nonzero, the Hasse property implies that $\alpha$ is represented globally, and the equality follows from Lemma \ref{lem:cl_num_glo_loc}.
\end{proof}

At times we will want to count local or global points of $X$ with a prescribed weighting.
\begin{defn}
  If $(\Gamma, V, X, X_0)$ is a composed variety over a local field $K_v$ with a chosen $\OO_v$-lattice $\Lambda_v \subseteq V(K_v)$, a \emph{local weighting} is a function $w_v : \Lambda_v \intsec X(K_v) \to \CC$ invariant under $\Gamma_{\Lambda_v}$.
\end{defn}
\begin{rem}
  At an archimedean place, by our convention $\Lambda_v = V(K_v)$ and $\Gamma_{\Lambda_v} = \Gamma(K_v)$, so a local weighting is simply a function on the $\Gamma(K_v)$-orbits of $X(K_v)$.
\end{rem}

\begin{defn}
  If $(\Gamma, V, X, X_0)$ is a composed variety over a global field $K$ with a chosen adelic lattice $\Lambda$, a \emph{global weighting} is a function $w : \Lambda \intsec X(\AA_K) \to \CC$ which can be expressed as a uniformly convergent product of local weightings,
  \begin{equation}\label{eq:glo_wtg}
    w(x) = \prod_{v \in \V^K} w_v(x_v).
  \end{equation}
\end{defn}

\begin{defn}
  If $w_v$ is a local weighting on $(\Gamma, V, X, X_0)$, define the \emph{weighted local lattice count} at $x_v \in X(K_v)$ by
  \[
    \ell_v^{w_v}(x_v) = \sum_{\substack{\gamma \in \Gamma_{\Lambda_v}\bs\Gamma(K_v) \\ \gamma x_v \in \Lambda_v}}
      w_v(\gamma x_v).
  \]
  If $\alpha_v \in H^1(K_v, M)$ is represented by $x_v$, set $\ell_v^{w_v}(\alpha_v) = \ell_v^{w_v}(x_v)$; otherwise set $\ell_v^{w_v}(\alpha_v) = 0$.
\end{defn}
\begin{rem}
  The sum is finite by Lemma \ref{lem:cl_num_glo_loc} and is independent of the choice of representatives for the left cosets by the $\Gamma_{\Lambda_v}$-invariance of $w_v$. If $w_v = 1$, we recover the local lattice count $\ell_v(x_v)$ of Definition \ref{def:loc_lat_cnt}.
\end{rem}

\begin{defn}
  If $w$ is a global weighting on $(\Gamma, V, X, X_0)$, define the \emph{weighted global lattice count} at $x \in X(K)$ by
  \[
    \ell^w(x) = \sum_{\substack{\gamma \in \Gamma_\Lambda\bs\Gamma(\AA_K) \\ \gamma x \in \Lambda}} w(\gamma x).
  \]
  If $\alpha \in H^1(K, M)$ is represented by $x$, set $\ell^w(\alpha) = \ell^w(x)$; otherwise set $\ell^w(\alpha) = 0$.
\end{defn}
\begin{rem}
  The sum is finite by Lemma \ref{lem:cl_num_glo_loc} and is independent of the choice of representatives for the left cosets by the $\Gamma_\Lambda$-invariance of $w$. If $w = 1$, we recover the global lattice count $\ell(x)$ of Definition \ref{def:glo_lat_cnt}.
\end{rem}

\begin{lem}\label{lem:Hasse_wtd}
  If $X$ is Hasse and $w = \prod_v w_v$ is a global weighting, then for every
  $\alpha \in H^1(K, M)$,
  \[
    \ell^w(\alpha) = \prod_{v \in \V^K} \ell_v^{w_v}(\alpha_v).
  \]
\end{lem}
\begin{proof}
  If some localization $\alpha_v$ is not represented by a point of $X(K_v)$,
  both sides vanish. Otherwise the Hasse property implies that $\alpha$ is
  represented globally. The equality then follows by factoring the defining
  adelic sum into its local sums, exactly as in Lemma \ref{lem:cl_num_glo_loc}.
\end{proof}

\begin{lem}\label{lem:finite_total_cl_num}
  There are only finitely many coclasses $\alpha \in H^1(K, M)$
  for which $\ell(\alpha) \neq 0$.
\end{lem}
\begin{proof}
  If $\ell(\alpha) \neq 0$, then by Lemma \ref{lem:spread}, the coclass $\alpha$ is unramified at all places $v$ outside a fixed finite set $S$. Thus $\alpha$ belongs to the Selmer group $H^1_S(K, M)$ of coclasses unramified outside $S$, which is well known to be finite (by the Hermite--Minkowski theorem, and its analogue for tame extensions of function fields; or alternatively by the discreteness of $H^1(K, M)$ inside $H^1(\AA_K, M)$, which will be used below).
\end{proof}

Thus the following is well defined.
\begin{defn}
  The \emph{(total) weighted class number} of $X$ is
  \[
    h^w(X) = \sum_{\alpha \in H^1(K, M)} \ell^w(\alpha).
  \]
  For the trivial weighting $w = 1$, we write $h(X) = h^1(X)$ and call it the \emph{class number} of $X$.
\end{defn}

\subsection{From total class numbers to orbit counts}

Our motivation for studying lattice counts is that their total over coclasses can be used to count orbits of integral forms of $\Gamma$ on $X$. For $\gamma \in \Gamma(\AA_K)$, let
\begin{equation}
  X_\gamma = X(K) \intsec \gamma^{-1} \Lambda \textand \Gamma_\gamma = \Gamma(K) \intsec \Gamma_{\gamma^{-1} \Lambda}.
\end{equation}
Then $\Gamma_\gamma$ acts on $X_\gamma$. For $\delta \in \Gamma_\Lambda$ and $\kappa \in \Gamma(K)$, we have
\begin{equation}
  \Gamma_{\delta \gamma \kappa} = \kappa^{-1} \Gamma_\gamma \kappa \textand X_{\delta \gamma \kappa} = \kappa^{-1} X_\gamma.
\end{equation}
In particular, the group $\Gamma_\gamma$, up to conjugacy, and the corresponding orbit set $\Gamma_\gamma \bs X_\gamma$ depend only on the class of $\gamma$ in the adelic class set $\Gamma_\Lambda \bs \Gamma(\AA_K) / \Gamma(K)$, whose cardinality is often called the \emph{class number} of $\Gamma$ and is well known to be finite. As it turns out, the total weighted class number counts orbits of $\Gamma_\gamma$ for all classes of $\gamma$:

\begin{lem}\label{lem:total_cl_num}
  Let $w$ be a global weighting on $X$. Then
  \begin{equation}\label{eq:total_cl_num}
    \frac{h^w(X)}{\size{H^0(K, M)}}
    =
    \sum_{\gamma \in \Gamma_\Lambda \bs \Gamma(\AA_K) / \Gamma(K)}
    \ \sum_{x \in \Gamma_\gamma \bs X_\gamma}
    \frac{w(\gamma x)}
    {\size{\Stab_{\Gamma_\gamma} x}}.
  \end{equation}
  
\end{lem}
\begin{proof}
  Since the cohomology parametrization of Proposition \ref{prop:composed} is injective and $\ell^w(\alpha) = 0$ for coclasses not represented by points of $X(K)$, the numerator on the left is
  \[
    h^w(X)
    =
    \sum_{x \in \Gamma(K) \bs X(K)}
    \ \sum_{\substack{\beta \in \Gamma_\Lambda \bs \Gamma(\AA_K) \\
                        \beta x \in \Lambda}}
    w(\beta x).
  \]
  Now subdivide $\Gamma_\Lambda \bs \Gamma(\AA_K)$ into double cosets
  $\Gamma_\Lambda \bs \Gamma(\AA_K) / \Gamma(K)$. Fix a representative $\gamma$ of one such double coset. The left cosets of $\Gamma_\Lambda$ contained in
  $\Gamma_\Lambda \gamma \Gamma(K)$ are naturally parametrized by
  $\Gamma_\gamma \bs \Gamma(K)$: indeed,
  \[
    \Gamma_\Lambda \gamma \kappa_1 = \Gamma_\Lambda \gamma \kappa_2
    \quad\Longleftrightarrow\quad
    \kappa_2 \kappa_1^{-1} \in \Gamma(K) \intsec \gamma^{-1} \Gamma_\Lambda \gamma
    = \Gamma_\gamma.
  \]
  Hence the contribution of this double coset is
  \begin{equation}\label{eq:x_dbl_coset}
    \sum_{x \in \Gamma(K) \bs X(K)}
    \ \sum_{\substack{\kappa \in \Gamma_\gamma \bs \Gamma(K) \\
                        \kappa x \in X_\gamma}}
    w(\gamma \kappa x).
  \end{equation}
  We now rearrange this as a sum over $\Gamma_\gamma$-orbits in $X_\gamma$. Fix $y \in X_\gamma$, and let $x$ be a representative for its $\Gamma(K)$-orbit. The terms of \eqref{eq:x_dbl_coset} which map to the $\Gamma_\gamma$-orbit of $y$ are represented by those $\kappa \in \Gamma(K)$ for which $\kappa x \in \Gamma_\gamma y$. After replacing $\kappa$ within its $\Gamma_\gamma$-coset, we may require $\kappa x = y$. If $\kappa_0 x = y$, all such $\kappa$ are exactly
  \[
    \kappa_0 \Stab_{\Gamma(K)} x,
  \]
  and two of them determine the same left coset in $\Gamma_\gamma \bs \Gamma(K)$ exactly when they differ by an element of $\Stab_{\Gamma_\gamma} y$. Thus the number of terms of \eqref{eq:x_dbl_coset} mapping to the orbit of $y$ is
  \[
    \frac{\size{\Stab_{\Gamma(K)} y}}
         {\size{\Stab_{\Gamma_\gamma} y}}.
  \]
  By Proposition \ref{prop:composed}\ref{comp:stab}, the $\Gamma$-stabilizer of every
  $y \in X(K)$ is isomorphic to $M$ as a Galois module, so
  \[
    \size{\Stab_{\Gamma(K)} y} = \size{H^0(K, M)}.
  \]
  Thus the contribution of the chosen double coset $\Gamma_\Lambda \gamma \Gamma(K)$ is
  \[
    \size{H^0(K, M)}
    \sum_{y \in \Gamma_\gamma \bs X_\gamma}
    \frac{w(\gamma y)}{\size{\Stab_{\Gamma_\gamma} y}}.
  \]
  Summing over the double cosets and dividing by $\size{H^0(K, M)}$ gives
  \eqref{eq:total_cl_num}.
\end{proof}

\begin{rem}
  If $\Gamma$ has class number $1$, there is just one term on the right-hand side of \eqref{eq:total_cl_num}, and the formula simplifies to
  \begin{equation}
    \frac{h^w(X)}{\size{H^0(K, M)}}
    = \sum_{x \in \Gamma_{\id} \bs X_{\id}}
    \frac{w(x)}
    {\size{\Stab_{\Gamma_{\id}} x}},
  \end{equation}
  which recognizably has the form of a class number as in \eqref{eq:Shin1}.
\end{rem}

\subsection{Fourier analysis of the local and global Tate pairings}
\label{sec:Fourier}

We now introduce the main innovative technique of this paper: Fourier analysis of local and global Tate duality. In structure we are indebted to Tate's celebrated thesis \cite{Tate_thesis}, in which he
\begin{enumerate}
  \item\label{it:Tate_local} constructs a perfect pairing on the additive group of a local field $K$, taking values in the unit circle $\CC^{N=1}$, and thus furnishing a notion of Fourier transform for $\CC$-valued $L^1$ functions on $K$;
  \item\label{it:Tate_product} derives thereby a pairing and Fourier transform on the adele group $\AA_K$ of a global field $K$;
  \item\label{it:Tate_Poisson} proves that the discrete subgroup $K \subseteq \AA_K$ is a self-dual lattice and that the Poisson summation formula 
  \begin{equation}
    \sum_{x \in K} f(x) = \sum_{x \in K} \widehat f(x)
  \end{equation}
  holds for all $f$ satisfying reasonable integrability conditions.
\end{enumerate}

In this paper, we work not with the additive group $K$ but with a Galois cohomology group $H^1(K, M)$. The needed theoretical result is \emph{Poitou-Tate duality}, a nine-term exact sequence of which the middle three terms are of main interest to us:
\begin{equation}\label{eq:Poitou-Tate}
(\text{finite kernel}) \longto H^1(K, M) \longto \sideset{}{'}{\bigoplus}_v H^1(K_v, M) \longto H^1(K, M')^\vee \longto (\text{finite cokernel}).
\end{equation}
This can be interpreted as saying that $H^1(K, M)$ and $H^1(K, M')$ map to dual lattices in the respective adelic cohomology groups, where, if $m$ is the exponent of $M$, we write $M' = \Hom(M, \mu_m)$ for the Tate dual:
\[
H^1(\AA_K, M) = \sideset{}{'}{\bigoplus}_v H^1(K_v, M) \textand 
H^1(\AA_K, M') = \sideset{}{'}{\bigoplus}_v H^1(K_v, M').
\]
For a local field $K$, the cup product followed by the local invariant map gives the perfect Tate pairing
\[
H^1(K, M) \cross H^1(K, M')
\longto H^2(K, \mu_m)
\xrightarrow{\inv_K} \frac{1}{m}\ZZ/\ZZ.
\]
Composing with $t \mapsto e^{2\pi i t}$, we regard this pairing as taking values in $\mu_m(\CC)$ and write it multiplicatively as $\langle \alpha, \beta \rangle$. The adelic pairing is the product of the local pairings:
\[
\<\{\alpha_v\}, \{\beta_v\}\> = \prod_v \<\alpha_v, \beta_v\>
\in \mu_m(\CC).
\]
Now, for any Schwartz--Bruhat function $g : H^1(\AA_K, M) \to \CC$ (in this paper, $g$ will always be compactly supported and locally constant), we have Poisson summation
\[
\sum_{\alpha \in H^1(K, M)} g(\alpha) = c_M \sum_{\beta \in H^1(K, M')} \widehat g(\beta)
\]
for some constant $c_M$ which we think of as the covolume of $H^1(K,M)$ as a lattice in the adelic cohomology (adjusted by the kernel and cokernel in \eqref{eq:Poitou-Tate}).

We apply Poisson summation to the weighted local lattice-count functions introduced above.
\begin{defn} \label{defn:dual}
  Let $K_v$ be a local field. For $i = 1,2$, let
  $\bigl(\Gamma^{(i)}, V^{(i)}, X^{(i)}, X_0^{(i)}\bigr)$ be composed varieties over
  $K_v$ with Tate-dual stabilizers $M^{(1)}$ and $M^{(2)}$, equipped with lattices
  $\Lambda_v^{(i)}$ and local weightings $w_v^{(i)}$. We say that the two weightings are \emph{dual} with \emph{duality constant} $c_v$ if
  \begin{equation}\label{eq:general_local_refl}
    \widehat{\ell}_v^{w_v^{(1)}} = c_v \cdot \ell_v^{w_v^{(2)}},
  \end{equation}
  where the normalization of the Fourier transform is given by
  \[
    \widehat{f}(\beta)
    = \frac{1}{\size{H^0(K_v, M^{(1)})}}
      \sum_{\alpha \in H^1(K_v, M^{(1)})}
      f(\alpha) \langle \alpha, \beta \rangle.
  \]
  An identity of the form \eqref{eq:general_local_refl} is called a
  \emph{local reflection theorem}. If the weightings $w_v^{(1)} = w_v^{(2)} = 1$ are trivial, we say that
  the two composed varieties with chosen lattices are \emph{naturally dual}.
\end{defn}

\begin{thm}[\textbf{local-to-global reflection engine}]
  \label{thm:main_compose}
  For $i = 1,2$, let
  $\bigl(\Gamma^{(i)}, V^{(i)}, X^{(i)}, X_0^{(i)}\bigr)$ be Hasse composed
  varieties over a global field $K$ whose stabilizers $M^{(1)}$ and $M^{(2)}$ are Tate duals.
  Equip them with adelic lattices $\Lambda^{(i)}$ and global weightings
  $w^{(i)} = \prod_v w_v^{(i)}$. Suppose that for every place $v$ the local weightings are dual with constant $c_v$. Then $c_v = 1$ for almost all $v$, and we obtain a global reflection identity:
  \begin{equation}\label{eq:main_compose}
    \frac{h^{w^{(2)}}(X^{(2)})}{\size{H^0(K, M^{(2)})}}
    = \left(\prod_v c_v^{-1}\right)
      \frac{h^{w^{(1)}}(X^{(1)})}{\size{H^0(K, M^{(1)})}}.
  \end{equation}
\end{thm}
\begin{proof}
  For $x \in H^1(K, M^{(i)})$, let
  \[
    \ell^{(i)}(x) = \prod_v \ell_v^{w_v^{(i)}}(x_v).
  \]
  For almost all $v$, by Lemma \ref{lem:spread}, $\ell_v^{(i)} = \1_{H^1_\ur(K_v, M^{(i)})}$, so $\ell^{(i)}$ is a Schwartz--Bruhat function, indeed locally constant of compact support. Moreover, for almost all $v$,
  \[
    \widehat{\ell}_v^{(1)} = \widehat{\1}_{H^1_\ur(K_v, M^{(1)})}
    = \1_{H^1_\ur(K_v, M^{(2)})} = \ell_v^{(2)},
  \]
  so almost all duality constants $c_v$ are $1$.
  By Lemma \ref{lem:Hasse_wtd},
  \[
    h^{w^{(i)}}(X^{(i)})
    = \sum_{\alpha \in H^1(K, M^{(i)})} \ell^{(i)}(\alpha).
  \]
  Since the adelic Tate pairing is the product of the local Tate pairings,
  the Fourier transform respects the product decomposition, and the local reflection
  identities give
  \[
    \widehat{\ell}^{(1)} = \left(\prod_v c_v\right) \ell^{(2)}.
  \]
  Poisson summation therefore yields
  \[
    h^{w^{(2)}}(X^{(2)})
    = c_{M^{(2)}} \left(\prod_v c_v^{-1}\right)
      h^{w^{(1)}}(X^{(1)}),
  \]
  where the constant $c_{M^{(2)}}$ depends only on the Galois module $M^{(2)}$. It can be pinned down by appealing to the Greenberg--Wiles formula, which can be thought of as Poisson summation on the characteristic function of an adelic open box $\prod_v L_v \subset H^1(\AA_K, M)$. The value of the constant (which arises from a global Euler characteristic computation) is
  \[
    c_{M^{(2)}}
    = \frac{\size{H^0(K, M^{(1)})}}{\size{H^0(K, M^{(2)})}},
  \]
  establishing \eqref{eq:main_compose}.
\end{proof}

\subsection{Examples}
Composed varieties are plentiful. In this section, $K$ is any field not of one of finitely many bad characteristics for which the exposition does not make sense.

\begin{examp}\label{ex:bin_cubic}
  Let $V = \Sym^3(K^2)$ and let $X_{D} \subset V$ be the closed variety of binary cubic forms $f$ over $K$ with fixed nonzero discriminant $D$. This has an algebraic action of $\SL_2$, which is transitive over $K^\alg$ (using that $\PGL_2(K^\alg)$ is triply transitive on $\PP^1$), and there is a convenient basepoint
  \[
  f_0(X,Y) = X^2 Y - \frac{D}{4} Y^3.
  \]
  The point stabilizer $M$ is isomorphic to $\ZZ/3\ZZ$, cyclically permuting the roots of $f_0$, but the nontrivial elements of $M$ are only defined over $K(\sqrt{D})$; that is, $M \isom \{0, \sqrt{D}, -\sqrt{D}\}$ as sets with Galois action. Coupled with the appropriate higher composition law (Theorem \ref{thm:hcl_cubic_ring}), this recovers the parametrization of cubic \'etale algebras with fixed quadratic resolvent by $H^1(K,M)$ in Proposition \ref{prop:Kummer_cubic}. To see that it is the same parametrization, note that a $\gamma \in \Gamma(K^\alg)$ that takes $f_0$ to $f$ is determined by where it sends the rational root $[1:0]$ of $f_0$, so the three $\gamma$'s are permuted by $\Gal(K^\alg/K)$ just like the three roots of $f$. In particular, $X_{D}$ is full.
\end{examp}
\begin{rem}
  Note that in this framework, we must choose $\SL_2$ and not $\GL_2$, which does not preserve the discriminant, and has a nonabelian point stabilizer $S_3$ (or even larger for the untwisted action: see \eqref{eq:twisted_action} for a discussion of the twist), corresponding to the six automorphisms of $K \cross K \cross K$. In the statement of a theorem like O-N (Theorem \ref{thm:O-N}), one can freely replace the group by a larger or smaller one of finite index, as the class numbers on either side are merely scaled.
\end{rem}
\begin{examp}
  Continuing with the sequence of known ring parametrizations, we might study the level variety $X$ of pairs of ternary quadratic forms with fixed discriminant $D_0$, which parametrize quartic rings \cite{B3}. This has one orbit over $K^\alg$ under the action of the group $\Gamma = \SL_2 \cross \SL_3$; unfortunately, the point stabilizer is isomorphic to the alternating group $A_4$, which is not abelian.
  
  So we narrow the group, which widens the ring of invariants and allows us to take a smaller closed level set $X$. We let $\Gamma = \SL_3$ alone act on pairs $(A, B)$ of ternary quadratic forms, which preserves the resolvent
  \[
    f(X,Y) = 4 \det\(A X + B Y\),
  \]
  a binary cubic form. We let $X$ be the closed variety of $(A,B)$ for which $f = f_0$ is a fixed separable polynomial. These parametrize quartic \'etale algebras $L$ over $K$ whose cubic resolvent $R$ is fixed. There is a natural base orbit $(A_0, B_0)$ whose associated $L \isom K \cross R$ has a linear factor. The point stabilizer $M \isom C_2 \cross C_2$, with the three non-identity elements permuted by $\Gal(K^\alg / K)$ in the same manner as the three roots of $f_0$. We have reconstructed the parametrization of quartic \'etale algebras with fixed cubic resolvent by $H^1(K,M)$ in \cite[Proposition~5.13]{QuarticExercises1}. In particular, $X$ is full.
\end{examp}
\begin{examp}
  The action of $\SL_2$ on binary quartic forms factors faithfully through $\PSL_2$. Let $X$ be the closed level set of binary quartic forms with fixed invariants $I$ and $J$, with the induced action of $\PSL_2$. The orbits of this space have been found useful for parametrizing $2$-Selmer elements of the elliptic curve $E : y^2 = x^3 + x I + J$ \cite{BhSh_BinQ}. The point stabilizer is $M \isom \ZZ/2\ZZ \cross \ZZ/2\ZZ$ with the Galois-module structure $E[2]$. This space $X$ embeds into the space of the preceding example via a map which we call the \emph{Wood embedding} after its prominent role in Wood's work \cite{WoodBQ}:
  \begin{align*}
    f &\mapsto (A,B) \\
    a x^4 + b x^3 y + c x^2 y^2 + d x y^3 + ey^4 &\mapsto
    \(
    \begin{bmatrix}
      & & 1/2 \\
      & -1 & \\
      1/2 & &
    \end{bmatrix},
    \begin{bmatrix}
      a & b/2 & c/3 \\
      b/2 & c/3 & d/2 \\
      c/3 & d/2 & e
    \end{bmatrix}\).
  \end{align*}
  In general, $X$ is \emph{not} full. For instance, over $K = \RR$, if $E$ has full $2$-torsion, there are only three kinds of binary quartics over $\RR$ with positive discriminant (positive definite, negative definite, and those with four real roots) which cover three of the four elements in $H^1(\RR, E[2])$. Two of these three (positive definite, four real roots) form the subgroup of elements whose corresponding $E$-torsor $z^2 = f(x,y)$ is soluble at $\infty$: these are the ones we retain when studying $\Sel_2 E$. The fourth element of $H^1(\RR, E[2])$ yields \'etale algebras whose $(A,B)$ has
  \[
  A = \begin{bmatrix}
    1/2 & & \\
    & 1 & \\
    & & 1/2
  \end{bmatrix},
  \]
  a conic with no real points. However, over global fields, the Hasse--Minkowski theorem implies that $X$ is Hasse.
\end{examp}

\begin{examp}\label{ex:Gm}
  For a \emph{non-example}, let $\Gamma = \GG_m$ act on the line $V = K$ by
  $\lambda(x) = \lambda^n x$. The geometric orbit $V - \{0\}$ has stabilizer
  $\mu_n$, but it is not closed, and the local lattice counts are infinite.

  To retain the stabilizer $\mu_n$ while making the ambient representation
  faithful, let $\Gamma = \GG_m$ act on $V = K^3$ by $\lambda(x, y, z) = (\lambda^n x, \lambda^{-n} y, \lambda z)$ and take
  \[
    X = \{(x,y,z) : xy = 1,\ z = 0\}.
  \]
  Then $X$ is closed, the representation on $V$ is faithful, and the stabilizer
  of every geometric point of $X$ is $\mu_n$. Thus $(\Gamma,V,X,X_0)$ is a
  composed variety, and Hilbert's Theorem 90 shows that it is full.
\end{examp}

Because of the flexibility afforded by faithful linear representations and closed homogeneous subvarieties, it is reasonable to suppose that any finite $K$-Galois module $M$ appears as the point stabilizer of some full composed variety over $K$. In the case of $M$ having trivial action, such a variety was constructed by Gundlach (\cite{Gundlach}, Theorem 2.3). Given a pair $X^{(1)}$, $X^{(2)}$ of such composed varieties with Tate dual stabilizers, any large enough lattices $\Lambda^{(i)}$ will admit weightings $w^{(i)}$ that are mutually dual, thereby getting reflection theorems from Theorem \ref{thm:main_compose}. What is more noteworthy is when a pair of lattices are \emph{naturally} dual at all finite places, especially when \emph{all} generic $\Gamma^{(1)}$-orbits $X^{(1)}$ having a $K$-point have a naturally dual orbit within $V^{(2)}$. This is the case for the reflection theorems of O-N type listed in Table \ref{tab:composed}.

\begin{rem}
  In the cubic and quartic cases, the lattices $\Lambda^{(i)}$ are \emph{dual} under an identification of the ambient space $ V $ with its dual $ V^* $ (which are isomorphic representations, up to an outer automorphism of $ \Gamma = \SL_3 $ in the quartic case). But in the quadratic case, while $V$ has a unique invariant $I = a D$, the dual representation $ V^* $ decomposes into $K^\alg$-orbits according to a different invariant $J = a/D^2$ and has no closed orbits besides the origin, raising the question of whether the alignment of the Fourier and additive duals may be coincidental in general.
\end{rem}

Many composed varieties have \emph{no} natural dual, as exemplified by the following.
\begin{examp}
  The group $\Gamma = \SO_2$ of rotations over $K = \QQ$ preserving the quadratic form $I(x,y) = x^2 + xy + y^2$ acts on the two-dimensional space $V$ of binary cubic forms of the shape
  \[
    f(x,y) = a x^3 - 3b x^2 y + (-3a - 3b)x y^2 - a y^3,
  \]
  which are exactly those symmetric under the threefold shift $x \mapsto y \mapsto -x-y \mapsto x$. This representation is used by Bhargava and Shnidman \cite{BhSh} to parametrize \emph{cyclic} cubic rings, that is, those with an automorphism of order $3$. The geometric orbits are the loci where $I(a,b) = a^2 + a b + b^2$ is fixed (note that $\disc f = 81 I^2$), which are closed and of finite stabilizer $M \isom \mu_3$ whenever $I \neq 0$, so we get composed varieties $X_I$ for all nondegenerate fibers of $I$ having a $K$-point. As may be checked, such $X_I$ are full, and there is a natural base orbit consisting of forms totally split over $\QQ$. However, natural duality fails in cases where $I$ is not a local unit.
  
  For example, let $p \equiv 1 \mod 3$ and take the fiber $I = p$, which has a rational point because $p$ splits in $\QQ(\sqrt{-3})$. Its discriminant is $81p^2$, hence has $p$-adic valuation $2$. Let $\Lambda_p$ be the lattice of forms integral at $p$. By Kummer theory, the cohomology group $H^1(\QQ_p, M) = H^1(\QQ_p, \ZZ/3\ZZ)$ is isomorphic to $\QQ_p^\cross/(\QQ_p^{\cross})^3$, which has dimension $2$ over $\FF_3$. A function $g$ on $H^1(\QQ_p, M)$ may be written as a matrix
  \[
  \begin{bmatrix}
    \multicolumn{1}{c|}{g(1)} & g(\alpha) & g(\alpha^2) \\ \hline
    g(\beta) & g(\alpha\beta) & g(\alpha^2\beta) \\
    g(\beta^2) & g(\alpha\beta^2) & g(\alpha^2\beta^2)
  \end{bmatrix}
  \]
  in which the trivial coclass and the unramified cohomology $\<\alpha\>$ are marked off by dividers.
  
  The six \emph{ramified} cohomology elements each have one integral orbit, corresponding to the maximal order in the corresponding cubic extension of $\QQ_p$; the three \emph{unramified} cohomology elements---the trivial coclass for $\QQ_p^3$, and the other two for the degree-$3$ unramified field extension in its two orientations---all have \emph{no} integral orbits, because the three orders
  \[
  \{(x_1, x_2, x_3) \in \ZZ_p^3 : x_i \equiv x_j \mod p\}
  \]
  are all asymmetric under the threefold automorphism of $\ZZ_p^3$. So we get a local lattice-count function
  \[
  \ell_p = 
  \begin{bmatrix}
    \multicolumn{1}{c|}{0} & 0 & 0 \\ \hline
    1 & 1 & 1 \\
    1 & 1 & 1
  \end{bmatrix}
  \]
  whose Fourier transform
  \[
  \widehat\ell_p =
  \begin{bmatrix}
    \multicolumn{1}{c|}{2} & -1 & -1 \\ \hline
    0 & 0 & 0 \\
    0 & 0 & 0
  \end{bmatrix}
  \]
  has mixed signs and thus cannot be the local lattice-count function of any composed variety. We might try to remedy this by changing the basepoint, but this results in non-real values of $\widehat\ell_p$. Similar obstructions to natural duality appear in typical examples drawn from the composed varieties parametrizing rings with automorphisms found by Gundlach \cite{Gundlach}.
\end{examp}

\section{The structure of \texorpdfstring{$H^1(K, M)$}{H¹(K, M)}}

\subsection{Realizing \texorpdfstring{$H^1(K, M)$}{H¹(K, M)} as a multiplicative subgroup}

Let $K$ be a field not of characteristic $p$, and let $M$ be a cyclic $G_K$-module of order $p$. In this section we begin to study the structure of the group $H^1(K, M)$, especially in the case that $K$ is a $p$-adic field. Much of the following generalizes the work of Del Corso--Dvornicich \cite{DCD}, with suitable adaptations of method.

Let $M' = \Hom(M, \mu_p)$ be the Tate dual of $M$. The Galois actions on $M$, $M'$, and $\mu_p$ can be described as characters
\[
  \chi_M, \chi_{M'}, \chi_\cyc : G_K \to \FF_p^\cross,
\]
the last of these being the familiar cyclotomic character. We have $\chi_M \cdot \chi_{M'} = \chi_\cyc$. The action of $G_K$ on $M'$ factors through the Galois group of a finite extension $E/K$, giving an isomorphism
\[
  \chi_{M'} : \Gal(E/K) \isom \mu_d,
\]
where $d = [E : K]$, a divisor of $p - 1$. Fix an abelian group isomorphism $\xi : M \isom \mu_p$ and note that $\xi$ is also an isomorphism of $E$-modules. By Kummer theory,
\[
  H^1(E, M) \isom H^1(E, \mu_p) \isom E^\cross/(E^\cross)^p.
\]
Our first proposition characterizes the image of the restriction map $\Res : H^1(K, M) \to H^1(E, M)$, which is injective because $\gcd(d, p) = 1$:
\begin{prop}[\textbf{Kummer theory for $H^1(K, M)$}] \label{prop:E_isotyp}
  Let $K$ be a field, $\ch K \neq p$. Let $M$ be a $K$-Galois module of order $p$, and let $M' = \Hom(M, \mu_p)$ be its Tate dual. The image of $H^1(K, M)$ in $H^1(E, M) \isom E^\cross/(E^\cross)^p$ is the $M'$-isotypical component, that is, the subgroup of $\alpha \in E^\cross/(E^\cross)^p$ such that for all $g \in \Gal(E / K)$,
  \begin{equation}
    g(\alpha) \equiv \alpha^{\chi_{M'}} \mod (E^\cross)^p.
  \end{equation}
\end{prop}
\begin{rem}
  The term ``isotypical component'' is apt because $E^\cross / (E^\cross)^p$ is a representation of $\Gal(E/K)$ of characteristic $p$, which is prime to $d = \size{\Gal(E/K)}$ so we have complete reducibility into irreducibles. This will be useful in the next subsection.
\end{rem}
\begin{proof}
  The case $p = 2$ is trivial since $E = K$, and the case $p = 3$ follows from the classical solution to the cubic: see \cite[Proposition 5.1]{OEtale}. These are the only cases that will be used in this paper. Nevertheless, we prove the general case for posterity's sake.
  
  Let $\sigma \in Z^1(K, M)$ be a cocycle. Let $\alpha \in E^\cross$ be a Kummer element for $\Res^E_K \sigma$. By definition, this means that there is a choice of $p$th root $\sqrt[p]{\alpha} \in K^\alg$ such that for all $h \in H \coloneqq G_E$,
  \begin{equation} \label{eqx:Kum_p}
    h\bigl(\sqrt[p]{\alpha}\bigr) = \xi\bigl(\sigma(h)\bigr) \cdot \sqrt[p]{\alpha}.
  \end{equation}
  (The choice of $p$th root can change when $\sigma$ is adjusted by a coboundary.)
  For $g \in G_K$, letting $j \in \ZZ$ be a lift of $\chi_{M'}(g)^{-1}$, we claim that $g(\alpha)^j$ is another valid Kummer element for $\Res^E_K \sigma$. We do this by showing by routine calculation that the $p$th root
  \[
    \frac{g\bigl(\sqrt[p]{\alpha}\bigr)^j}{\xi\bigl(\sigma(g)\bigr)}
  \]
  of $g(\alpha)^j$ may replace $\sqrt[p]{\alpha}$ in \eqref{eqx:Kum_p}.
  {The proof is a calculation:
  \begin{align*}
    h\left(\frac{g\bigl(\sqrt[p]{\alpha}\bigr)^j}{\xi\bigl(\sigma(g)\bigr)}\right)
    &= \frac{hg\bigl(\sqrt[p]{\alpha}\bigr)^j} {h \xi\bigl(\sigma(g)\bigr)} \\
    &= \frac{g\Bigl(g^{-1}hg\bigl(\sqrt[p]{\alpha}\bigr)\Bigr)^j}
    {\xi\bigl(h \sigma(g)\bigr)} \\
    &= \frac{g\Bigl(\xi\bigl(\sigma(g^{-1} h g)\bigr) \cdot \sqrt[p]{\alpha}\Bigr)^j}
    {\xi\bigl(h \sigma(g)\bigr)} \\
    &= \xi\bigl(j \cdot \chi_\cyc(g) \cdot \sigma(g^{-1} h g) - h \sigma(g) \bigr) g\bigl(\sqrt[p]{\alpha}\bigr)^j \\
    &= \xi\bigl(g \sigma(g^{-1} h g) - h \sigma(g) \bigr) g\bigl(\sqrt[p]{\alpha}\bigr)^j \\
    &= \xi\Bigl(\bigl(\sigma(gg^{-1}hg) - \sigma(g)\bigr) - \bigl(\sigma(hg) - \sigma(h)\bigr)\Bigr) g\bigl(\sqrt[p]{\alpha}\bigr)^j \\
    &= \xi\bigl(\sigma(h)\bigr) \cdot \frac{g\bigl(\sqrt[p]{\alpha}\bigr)^j}{\xi\bigl(\sigma(g)\bigr)}.
  \end{align*}}%
  Since the Kummer element $\alpha$ of $\Res^E_K \sigma$ is unique up to $p$th powers, we derive that $\alpha$ belongs to the $\chi_{M'}$-isotypical component.
  
  Conversely, suppose that $\alpha \in E^\cross$ is a representative of a member of the $\chi_{M'}$-isotypical component of $E^\cross/(E^\cross)^p$. Let
  \[
    I = \{i \in \ZZ : i^d \equiv 1 \mod p\}.
  \]
  For each $i \in I$, there is a $\beta_i \in E^\cross$ such that
  \begin{equation} \label{eqx:beta}
    g_i(\alpha)^i = \beta_i^p \alpha,
  \end{equation}
  where $g_i \in \Gal(E/K)$ is the unique element with $\chi_{M'}(g_i)^{-1} = i$. Although $I$ is infinite, we can choose $\beta_i$ for a representative of each of the $d$ relevant congruence classes and set
  \[
    \beta_{i + kp} = g_i(\alpha)^{k} \beta_i, \quad k \in \ZZ.
  \]
  Now for $g \in G_K$, define
  \begin{equation} \label{eqx:sigma}
    \sigma(g) = \xi^{-1}\left(\frac{g\bigl(\sqrt[p]{\alpha}\bigr)^i}{\beta_i \cdot \sqrt[p]{\alpha}}\right),
  \end{equation}
  where $i \in \ZZ$ is a lift of $\chi_{M'}(g)^{-1}$. Expanding out, we find that $\sigma$ is a cocycle if and only if the $\beta_i$ satisfy the relation
  \begin{equation} \label{eqx:beta_cocycle}
    \beta_{ij} = \beta_i \cdot g_i(\beta_j)^i,
  \end{equation}
  Combining the relevant instances of \eqref{eqx:beta} shows merely that the $p$th powers of both sides of \eqref{eqx:beta_cocycle} are equal. Let $\gamma_{i,j}$ be the discrepancy in \eqref{eqx:beta_cocycle}:
  \[
    \gamma_{i,j} = \frac{\beta_{ij}}{\beta_i \cdot g_i(\beta_j)^i}.
  \]
  We can check that $\gamma_{i,j}$ depends only on the congruence classes of $i$ and $j$ modulo $p$ and defines a cocycle $\gamma \in H^2\bigl(\mu_d, M\bigr)$, which vanishes since $\gcd(d,p) = 1$. Hence there are constants $\epsilon_i \in \mu_p(E)$ such that
  \[
    \gamma_{i,j} = \frac{\epsilon_{ij}}{\epsilon_i \cdot g_i(\epsilon_j)^i},
  \]
  and upon dividing the $\beta_i$ by the $\epsilon_i$, we satisfy the cocycle condition \eqref{eqx:beta_cocycle}. So $\sigma \in Z^1(K, M)$ is a cocycle. 
  Setting $i = 1$ in \eqref{eqx:sigma}, we see that $\sigma$ has associated Kummer element $\alpha$.
\end{proof}

\subsection{Isotypical components of the multiplicative group modulo \texorpdfstring{$p$}{p}th powers}

Now we specialize to the case that $K$ is a finite extension of $\QQ_p$. To the author's knowledge, no one has computed the representation structure of the groups $E^\cross / (E^\cross)^p$ appearing in Proposition \ref{prop:E_isotyp}. However, Del Corso and Dvornicich \cite{DCD} compute the representation structure of the  closely related group $F^\cross / (F^\cross)^p$, where $F$ is the compositum of all cyclic extensions of $K$ of degree dividing $p - 1$. 

\begin{lem} \label{lem:misc_bij_p-1}
Let $K$ be a finite extension of $\QQ_p$, and let $F$ be the compositum of all cyclic extensions of $K$ of degree dividing $p - 1$. The following groups are canonically isomorphic:
\begin{enumerate}[(a)]
  \item\label{misc:mod} The Galois modules $M$ over $K$ of order $p$, up to isomorphism, under tensor product;
  \item\label{misc:rep} The irreducible representations $\chi$ of $\Gal(F/K)$, under tensor product;
  \item\label{misc:elt} $K^\cross / (K^\cross)^{p-1}$.
\end{enumerate}
\end{lem}
\begin{proof}
  A Galois module $M$ is given by a map $\tilde\chi : G_K \to \Aut(C_p) \isom \FF_p^\cross$. Such a $\chi$ necessarily factors through a finite cyclic quotient $\Gal(E/K)$ of order dividing $p - 1$, and hence factors through $\Gal(F/K)$, establishing a map from \ref{misc:mod} to \ref{misc:rep}. For the inverse, since $\mu_{p-1} \subset K$, the irreducible representations $\chi$ of $\Gal(F/K)$ are one-dimensional and can thus be viewed as Galois modules. Note that since $\GL_1(\FF_p) = \FF_p^\cross$ is abelian, two representations $\chi_i : \Gal(F/K) \to \FF_p^\cross$ are isomorphic if and only if they are equal.
  
  Finally, the relation of \ref{misc:elt} to the other two comes from Kummer theory, since $\mu_{p-1} \subseteq K$.
\end{proof}

Now suppose that, as in \cite{DCD}, $F$ is the compositum of all cyclic extensions of a $p$-adic field $K$ of degree $p - 1$. By Kummer theory, such extensions are of the form $K[\sqrt[p-1]{\alpha}]$ for $\alpha \in K^\cross/\bigl(K^\cross\bigr)^{p - 1}$. The group $K^\cross/\bigl(K^\cross\bigr)^{p - 1}$ has the structure $C_{p-1} \cross C_{p-1}$, with a basis $\{u, \pi_K\}$ where $u$ is a lift of a generator of $k_K^\cross$. Hence $F = K[\sqrt[p-1]{u}, \sqrt[p-1]{\pi_K}]$, and the pairing
\begin{equation} \label{eq:Kum_pairing}
  \begin{tabular}{c@{}c@{}c@{}c@{}c}
    $\Gal(F/K)$ & ${} \cross {}$ & $K^\cross/\bigl(K^\cross\bigr)^{p - 1}$ & ${} \to {}$ & $\mu_{p-1}$ \\
    $(g$ & $,$ & $\alpha)$ & ${} \mapsto {}$ & $\dfrac{g(\sqrt[p-1]{\alpha})}{\sqrt[p-1]{\alpha}}$
  \end{tabular}
\end{equation}
is perfect. Denote by $\chi_\alpha$ the one-dimensional representation of $\Gal(F/K)$ over $\FF_p$ given by pairing with $\alpha$.

\begin{nota}
If a Galois module $M$, a representation $\chi$, and an element $\alpha$ correspond under the bijections of Lemma \ref{lem:misc_bij_p-1}, we say that $\alpha$ is the \emph{Kummer element} of $M$ or of $\chi$, and we write $M = M_\alpha$, $\alpha = \alpha_M$, $\chi = \chi_\alpha$. The Kummer element is uniquely characterized by the relation 
\begin{equation}
  g|_{M} = \chi(g) = \frac{g(\sqrt[p-1]{\alpha})}{\sqrt[p-1]{\alpha}} \mod p
\end{equation}
for all $g \in G_K$ (or all $g \in \Gal(F/K)$).
\end{nota}

\begin{lem} \label{lem:E_F}
If $K \subseteq E \subseteq F$ is an intermediate field, then the natural map $E^\cross / (E^\cross)^p \to F^\cross / (F^\cross)^p$ is an embedding whose image is the sum of the $\chi$-isotypical components $W_\chi$ of $F^\cross / (F^\cross)^p$ for which $\chi : \Gal(F/K) \to \FF_p^\cross$ factors through $\Gal(E/K)$.
\end{lem}
\begin{proof}
  Put $G = \Gal(F/K)$, $H = \Gal(F/E)$, and $r = \size{H} = [F : E]$. Since $r \mid \size{G} = (p - 1)^2$, we have $p \nmid r$. If $a \in E^\cross$ becomes a $p$th power in $F^\cross$, say $a = b^p$, then
  \[
    a^r = N_{F/E}(b)^p,
  \]
  so $[a]^r = 1$ in $E^\cross/(E^\cross)^p$. Since this group has exponent $p$ and $p \nmid r$, we have $[a] = 1$. Thus the natural map is injective.

  Let $W = F^\cross/(F^\cross)^p$. The image of $E^\cross/(E^\cross)^p$ in $W$ is contained in $W^H$. Conversely, if $[b] \in W^H$, then
  \[
    [N_{F/E}(b)] = \prod_{h \in H} h[b] = [b]^r.
  \]
  Choose $s \in \ZZ$ with $r s \equiv 1 \mod p$. Then the image of $N_{F/E}(b)^s \in E^\cross$ in $W$ is $[b]$, so the image is exactly $W^H$.

  Since $p \nmid \size{G}$, the $\FF_p[G]$-module $W$ is semisimple; and since $G$ is abelian of exponent dividing $p - 1$, all of its irreducible representations over $\FF_p$ are one-dimensional. Hence
  \[
    W^H = \bigoplus_{\chi|_H = 1} W_\chi.
  \]
  The condition $\chi|_H = 1$ is equivalent to $\chi$ factoring through $G/H = \Gal(E/K)$, proving the claim.
\end{proof}

The following result deserves to be better known.
\begin{prop} \label{prop:p-1st_root_of_-p}
  The Kummer element $\alpha_{\mu_p}$ of the cyclotomic module $\mu_p$ is $-p$.
\end{prop}
\begin{rem}
  In particular, $\QQ_p[\zeta_p] = \QQ_p[\sqrt[p-1]{-p}]$ \cite[Ch.~4, \textsection2, Exercise 15]{BorSha66}.
\end{rem}
\begin{proof}
  It suffices to prove this relation for $K = \QQ_p$, as passing to an extension does not change the Kummer element. Let $L = \QQ_p[\sqrt[p-1]{-p}]$, a totally ramified extension with uniformizer $\pi_L = \sqrt[p-1]{-p}$. Consider the polynomial
  \[
    f(x) = \frac{1 - (1 + \pi_L x)^p}{\pi_L p}.
  \]
  Since $p \nmid \binom{p}{k}$ for $0 < k < p$, we see that $f(x) \in \OO_L[x]$ and
  \[
    f(x) \equiv x^p - x \mod \sqrt[p-1]{-p}.
  \]
  Each element of $\FF_p$ is a simple root of $x^p - x$. So by Hensel's lemma, $f$ has $p$ distinct roots; for each $a \in \FF_p$, there is a unique $p$th root of unity $\zeta_{p,a} \in L$ satisfying the congruence
  \[
    \zeta_{p,a} \equiv 1 + \pi_L a \mod \pi_L^2.
  \]
  In particular, $L$ has all the $p$th roots of unity and thus equals $\QQ_p[\zeta_p]$.
  
  To complete the proof, we must show that the two Galois actions on $L = \QQ_p[\zeta_p]$ coincide, for to say that $\mu_p$ has Kummer element $-p$ is to claim that, for $r \in \mu_{p-1}$, the automorphism $\sigma_r$ of $L$ given by
  \[
    \sigma_r(\pi_L) = r \cdot \pi_L
  \]
  sends each $\zeta \in \mu_p$ to $\zeta^r$ (exponent reduced modulo $p$). We have
  \[
    \sigma_r(\zeta_{p,a}) \equiv \sigma_r(1 + \pi_L a) \equiv 1 + r \pi_L a \equiv (1 + \pi_L a)^r \equiv \zeta_{p,a}^r \mod \pi_L^2,
  \]
  and hence
  \[
     \sigma_r(\zeta_{p,a}) = \zeta_{p,a}^r
  \]
  since the $p$th roots of unity are pairwise incongruent modulo $\pi_L^2$.
\end{proof}

Let $F$ be a $p$-adic field. The structure of $F^\cross/(F^\cross)^p$ is well known. Following Nguyen-Quang-Do \cite{Nguyen}, we filter $F^\cross$ by the subgroups
\[
  U^i(F) = \{1 + \pi_F^i a : a \in \OO_F\},
\]
where we put $U^0(F) = \OO_F^{\cross}$ and $U^{-1}(F) = F^\cross$. Write
\[
  \gr^i U(F) = U^i(F)/U^{i+1}(F),
\]
so that
\[
  \gr^i U(F) \isom \begin{cases}
    \ZZ & i = -1, \\
    k_F^\cross & i = 0, \\
    k_F & i \geq 1.
  \end{cases}
\]
Let $\mathcal U^i(F)$ be the image of $U^i(F)$ in $F^\cross/(F^\cross)^p$, and put
\[
  \gr^i \mathcal U(F)
  = \mathcal U^i(F)/\mathcal U^{i+1}(F)
  = U^i(F) / \Bigl(U^{i+1}(F) \cdot \bigl(U^i(F) \intsec (F^\cross)^p\bigr)\Bigr).
\]
The following well-known result gives the structure of these graded pieces:

\begin{prop}[\textbf{Shafarevich basis}] \label{prop:Sh_basis}
  Let $F \supseteq \QQ_p$ be a local field. The groups $\gr^i \mathcal U(F)$, for $i \geq -1$, are as follows:
  \begin{enumerate}[(a)]
    \item\label{Sh:uzr} $\gr^{-1} \mathcal U(F) \isom C_p$, generated by the class of a uniformizer $\pi_F$.
    \item\label{Sh:generic} If
    \[
      0 < i < \frac{p}{p-1} e_F, \quad p \nmid i,
    \]
    then the natural map $\gr^i U(F) \to \gr^i \mathcal U(F)$ is an isomorphism.
    \item\label{Sh:intimate} If $i = \frac{p}{p-1} e_F$ and $\mu_p \subseteq F$, then $\gr^i \mathcal U(F) \isom C_p$.
    \item\label{Sh:0} In all other cases, $\gr^i \mathcal U(F) = 0$.
  \end{enumerate}
\end{prop}
\begin{proof} See \cite[Proposition 6]{DCD}.
\end{proof}

If we concatenate $\FF_p$-bases for all the nonzero $\gr^i \mathcal U(F)$, arbitrarily lifted, we obtain an $\FF_p$-basis for $F^\cross/(F^\cross)^p$, commonly called the \emph{Shafarevich basis}. As the three types of elements in the Shafarevich basis have no standard names, we propose the following:
\begin{itemize}
  \item The uniformizer $\pi_F$;
  \item The $e_F f_F$ \emph{generic units} $1 + \pi_F^i a$ arising from part \ref{Sh:generic} of Proposition \ref{prop:Sh_basis}, where $0 < i < \frac{p}{p-1} e_F$ with $p \nmid i$, and where $a$ runs over an $\FF_p$-basis for $k_F \isom \FF_{p^{f_F}}$;
  \item The \emph{intimate unit} appearing in part \ref{Sh:intimate} of Proposition \ref{prop:Sh_basis} when $F$ has a primitive $p$th root of unity $\zeta_p$. It spans the final nonzero unit layer and is a Kummer element for the unique unramified extension of $F$ of degree $p$. Although not needed here, one can show that a formula for the intimate unit is
  \[
    u_{\mathrm{int}} = 1 + \beta p (1 - \zeta_p)
  \]
  where $\beta \in \OO_F$ is an element such that $\tr_{k_F / \FF_p} \bar\beta \neq 0$.
\end{itemize}

We now study the Galois action on the graded pieces of $F^\cross/(F^\cross)^p$.
\begin{prop}\label{prop:Gal_on_gr}
As a representation of $\Gal(F/K)$ over $\FF_p$, the structure of $\gr^j \mathcal U(F)$ is as follows:
\begin{enumerate}[(a)]
  \item\label{GalGr:uzr} $\gr^{-1} \mathcal U(F) \isom \chi_1$, the trivial representation.
  \item\label{GalGr:generic} For $0 < j < pe_K$ with $p \nmid j$,
  \[
    \gr^{j} \mathcal U(F) \isom \bigoplus_{i = 0}^{p - 2} \chi_{\pi_K^j u^i}^{\oplus f_K}.
  \]
  \item\label{GalGr:intimate} $\gr^{pe_K} \mathcal U(F) \isom \mu_p \isom \chi_{-p}$.
\end{enumerate}
\end{prop}
\begin{proof}
  This is essentially \cite[Proposition 7]{DCD}. In translating the notation of \cite{DCD} to our setting, note that their $W^{(j)}$ is our $\gr^j \mathcal U(F)$ when $j \geq 1$, while their $W^{(0)}$ is our $\gr^{-1} \mathcal U(F)$; and that $W_{m n}^{(j)}$, for $m, n \in \FF_p^\cross$, is the isotypical component of $W^{(j)}$ corresponding to the representation $\chi_{m,n} : (\sigma \mapsto m, \tau \mapsto n)$, where $[\sigma,\tau]$ is the basis of $\Gal(F/K)$ dual to $[\pi_K,u]$ for a suitable normalization of $u$. To speak about such a dual basis, one must fix a root $r = \zeta_{p-1}$ (which choice is avoided in our formulation), and the relevant translation of notation is
  \[
    \chi_{m,n} = \chi_{\pi_K^{\log_r m} u^{\log_r n}}.
  \]
  For part \ref{GalGr:intimate}, the formulation of \cite{DCD} is that $W^{(pe_K)} = W^{(pe_K)}_{kh}$, where $k$ and $h$ are defined in terms of the action of $\Gal(F/K)$ on $\mu_p$. Thus $\gr^{pe_K} \mathcal U(F) \isom \mu_p$, and the isomorphism $\mu_p \isom \chi_{-p}$ follows from Proposition \ref{prop:p-1st_root_of_-p}.
\end{proof}

It is clear by now that an important invariant of $M$ is the valuation $\theta(M) = v_K(\alpha_M) \in \ZZ/(p-1)\ZZ$, which we call the \emph{module offset} of $M$. The following properties are trivial but worth recording:
\begin{prop}\label{prop:module_offset}\
  \begin{enumerate}[(a)]
    \item\label{off:+} $\theta(M \tensor N) = \theta(M) + \theta(N)$ and $\theta\bigl(\Hom(M,N)\bigr) = \theta(N) - \theta(M)$;
    \item\label{off:mu_p} $\theta(\mu_p) = e$ and $\theta(M') = e - \theta(M)$;
    \item\label{off:unr} A module $M$ is unramified if and only if $\theta(M) = 0$.
  \end{enumerate}
\end{prop}

\subsection{Level spaces}
Following Nguyen-Quang-Do \cite[Definition 1.1]{Nguyen}, we make the following definition:
\begin{defn}
The \emph{defect} $d(\delta)$ of an element $[\delta] \in F^\cross / (F^\cross)^p$ is the maximal value of $v_F(\delta - 1)$ for representatives $\delta \in \OO_F^\cross$ of the coset $[\delta]$. If no representative in $\OO_F^\cross$ exists, we put $d(\delta) = -1$. If $\sigma \in H^1(K,M)$ is a coclass, we define its defect $d(\sigma) = d(\delta)$ to be the defect of the corresponding Kummer element $\sigma$ under the bijections of Proposition \ref{prop:E_isotyp} and Lemma \ref{lem:E_F}.
\end{defn}
\begin{prop}\label{prop:defect} The defect $d(\delta)$ of an element $[\delta] \in F^\cross / (F^\cross)^p$ lies in $\{-1, 0, 1, \ldots, pe\} \union \{\infty\}$.
\end{prop}
\begin{proof}
  Immediate from Proposition \ref{prop:Sh_basis}.
\end{proof}
\begin{defn}
The \emph{level} of a coclass $\sigma \in H^1(K, M)$ is defined by
\[
  \ell(\sigma) = \floor{\frac{d(\delta)}{p}}.
\]
For $i \in \ZZ$, define the \emph{level space}
\[
  \L_i = \{\sigma \in H^1(K, M) : \ell(\sigma) \geq i\}.
\]
\end{defn}
\begin{prop}\label{prop:levels_new} The level $\ell(\sigma)$ and the level spaces $\L_i$ have the following properties:
  \begin{enumerate}[(a)]
    \item\label{lev:range} $\ell(\sigma) \in \{-1, 0, 1, \ldots, e\} \union \{\infty\}$.
    \item\label{lev:defect} The level $\ell = \ell(\sigma)$ of a coclass $\sigma \in H^1(K, M)$ uniquely determines the defect $d = d(\sigma)$ in the following way:
    \begin{enumerate}[(i)]
      \item\label{def:-1} If $\ell = -1$, then $d = -1$.
      \item\label{def:generic} If $0 \leq \ell < e$, then $d$ is the unique solution in the interval $p\ell < d < p(\ell + 1)$ to the congruence
      \begin{equation} \label{eq:dp1}
        d \equiv v_{K}(\alpha_{M'}) \mod p - 1.
      \end{equation}
      \item\label{def:e} If $\ell = e$, then $d = pe$.
      \item\label{def:oo} If $\ell = \infty$, then $d = \infty$.
    \end{enumerate}
    \item\label{lev:subgp} For $i \in \ZZ$, the level space $\L_i$ is a subgroup of $H^1(K, M)$.
    \item\label{lev:ur} $\L_{e+1} = \{1\}$, and $ \L_e = H^1_{\ur}(K,M)$.
    \item\label{lev:size_Li} For $ 0 \leq i \leq e $,
    \[
      \size{\L_i} = q^{e - i} \size{H^0(K, M)}.
    \]
    \item\label{lev:size_all} $ \L_{-1} $ is the whole of $ H^1(K,M) $, and
    \[
      \size{H^1(K,M)} = q^e \size{H^0(K, M)} \cdot \size{H^0(K, M')}.
    \]
  \end{enumerate}
\end{prop}
\begin{proof}
  \begin{enumerate}[(a)]
    \item Immediate from Proposition \ref{prop:defect}.
    \item Cases \ref{def:-1}, \ref{def:e}, and \ref{def:oo} are immediate from Proposition \ref{prop:defect}. In case \ref{def:generic}, let $\sigma$ have defect $d$, where $p\ell \leq d < p(\ell + 1)$. Then $\sigma$ has nonzero projection to the $M'$-isotypical component of $\gr^d \mathcal U(F)$. Note that $d = p\ell$ is excluded because $\gr^d \mathcal U(F)$ vanishes there by Proposition \ref{prop:Sh_basis}\ref{Sh:0}. By Proposition \ref{prop:Gal_on_gr}\ref{GalGr:generic}, we must have
    \[
    \chi_{\pi_K^d u^h} \isom M'
    \]
    for some $h \in \ZZ$, establishing \eqref{eq:dp1}. Among the $p - 1$ integers strictly between $p\ell$ and $p(\ell + 1)$, exactly one satisfies \eqref{eq:dp1}.
    \item Simply note that $\L_i$ is the preimage in $H^1(K,M)$ of a multiplicative subgroup of $F^\cross/(F^\cross)^p$, namely
    \begin{itemize}
      \item $\mathcal U^{-1}(F)$ if $i \leq -1$
      \item $\mathcal U^{pi}(F)$ if $i \geq 0$.
    \end{itemize}
    \item The first claim follows from part \ref{lev:range}. For the second claim, note that if $\ell(\sigma) = e$, then $\delta$ has the maximal possible finite defect $pe$ and thus lies in the intimate space $\mathcal U^{pe}(F) \isom \gr^{pe} \mathcal U(F)$. By Proposition \ref{prop:Gal_on_gr}\ref{GalGr:intimate}, $M' \isom \mu_p$, so $M$ has the trivial action and $\delta$ is a Kummer element for the unique unramified degree $p$ extension $L/K$ (see \cite[Lemmas 5 and 8]{DCD}).
    \item The case $i = e$ follows from the preceding part, since
    \[
      \size{H^1_\ur(K,M)} = \size{H^0(K, M)}.
    \]
    Thus it suffices to prove for $0 \leq i < e$ that
    \begin{equation} \label{eqx:LLq}
      \size{\L_i / \L_{i+1}} = q.
    \end{equation}
    The nonidentity cosets of $\L_i / \L_{i+1}$ consist of coclasses of level $i$. By part \ref{lev:defect}, these coclasses all have a common defect $d$. For this $d$, by Proposition \ref{prop:Gal_on_gr}\ref{GalGr:generic},
    \[
      \size{\L_i / \L_{i+1}} = \size{\bigl(\gr^d \mathcal U(F)\bigr)_{M'}} = \size{C_p^{f_K}} = q,
    \]
    establishing \eqref{eqx:LLq}.
    \item Finally, by Proposition \ref{prop:Gal_on_gr}\ref{GalGr:uzr},
    \[
      \L_{-1} / \L_0 \isom \bigl(\gr^{-1} \mathcal U(F)\bigr)_{M'} \isom \begin{cases}
        C_p & M' \text{ has trivial action} \\
        0 & \text{otherwise.}
      \end{cases} \qedhere
    \]
  \end{enumerate}
\end{proof}
\begin{rem}
  By Tate duality, part \ref{lev:size_all} can also be written
  \[
    \frac{\size{H^0(K, M)} \cdot \size{H^2(K, M)}}{\size{H^1(K, M)}} = p^{-[K : \QQ_p]},
  \]
  thus recovering an Euler characteristic result.
\end{rem}

\subsection{The Tate pairing}
Now that we understand $H^1(K, M)$, we must understand the Tate pairing
\[
  H^1(K, M) \cross H^1(K, M') \to \mu_p
\]
by which the Fourier transform in Theorem \ref{thm:main_compose} is defined. In particular, we must understand how the Tate pairing behaves on level spaces.

\begin{lem}\label{lem:Tate_Hilbert}
  The Tate pairing between $H^1(K, M)$ and $H^1(K, M')$ is the restriction of the Hilbert pairing on $H^1(F, \mu_p) \isom F^\cross / (F^\cross)^p$.
\end{lem}
\begin{proof}
  In the commutative diagram
  \begin{equation} \label{eq:Tate_Hilbert}
  \begin{tikzcd}[column sep=large]
    {} & H^1(K, M) \cross H^1(K, M') \arrow[r, "\cup"] \arrow[d, "\Res \cross \Res"]
       & H^2(K, \mu_p) \arrow[d, "\Res"] \arrow[r, "e^{2\pi i\,\inv_K}", "\sim"']
       & \mu_p \arrow[d, equal] \\
    F^\cross/(F^\cross)^p \cross F^\cross/(F^\cross)^p \arrow[r, "\sim"]
       & H^1(F, \mu_p) \cross H^1(F, \mu_p) \arrow[r, "\cup"]
       & H^2(F, \mu_p) \arrow[r, "e^{2\pi i\,\inv_F}", "\sim"']
       & \mu_p
  \end{tikzcd}
  \end{equation}
  the top row is the Tate pairing over $K$, the bottom row the Hilbert pairing. The commutativity of the right-hand square is justified as follows: we have
  \[
    \inv_F \circ \res = [F : K] \cdot \inv_K
  \]
  (see Serre \cite{SerreLF}, Prop.~XIII.7), but multiplication by $[F : K] = (p - 1)^2 \equiv 1 \mod p$ preserves $\mu_p$. This yields the desired identification of pairings.
\end{proof}
\begin{rem}
In particular, the Hilbert pairing on $F^\cross / (F^\cross)^p$ is perfect between the $M$- and $M'$-isotypical components for all $M$, and in fact vanishes on all other pairs of components, as one can deduce from the Galois equivariance of the Hilbert symbol.
\end{rem}

\begin{lem}\label{lem:Tate_level}
  For $ -1 \leq i \leq e + 1 $, with respect to the Tate pairing between $ H^1(K,M) $ and $ H^1(K, M') $,
  \[
  \L_i(M)^\perp = \L_{e-i}(M').
  \]
\end{lem}
\begin{proof}
  For $i = -1$ or $i = e + 1$, the result is trivial. Suppose $0 \leq i \leq e$. Under the Hilbert pairing on $F^\cross/(F^\cross)^p$, the standard orthogonality relation for the unit filtration is
  \begin{equation}\label{eq:unit_Hilbert_orth}
    \bigl(\mathcal U^j(F)\bigr)^\perp = \mathcal U^{pe-j+1}(F), \qquad 1 \leq j \leq pe,
  \end{equation}
  where the unit groups are understood modulo $p$th powers; see Nguyen-Quang-Do \cite[Proposition 3.4.1.1]{Nguyen}. For $0 < i < e$, Proposition \ref{prop:Sh_basis}\ref{Sh:0} gives
  \[
    \mathcal U^{pi}(F) = \mathcal U^{pi+1}(F),
  \]
  so \eqref{eq:unit_Hilbert_orth} gives
  \[
    \bigl(\mathcal U^{pi}(F)\bigr)^\perp
      = \mathcal U^{p(e-i)}(F).
  \]
  At the endpoints, $\mathcal U^0(F)=\mathcal U^1(F)$ modulo $p$th powers, while
  $\bigl(\mathcal U^{pe}(F)\bigr)^\perp=\mathcal U^1(F)=\mathcal U^0(F)$.
  Taking the $M$- and $M'$-isotypical components and using Lemma \ref{lem:Tate_Hilbert} now yields
  $\L_i(M)^\perp=\L_{e-i}(M')$.
\end{proof}

\begin{cor}\label{cor:levels}
The characteristic function $L_i$ of the level space $\L_i$ has Fourier transform given by
\begin{equation}\label{eq:hat 1 Li}
  \widehat L_i = \begin{cases}
    q^e \size{H^0(K, M')} L_{e + 1} & i = -1 \\
    q^{e - i} {L_{e-i}} & 0 \leq i \leq e \\
    \dfrac{1}{\size{H^0(K, M)}}L_{-1} & i = e + 1.
  \end{cases}
\end{equation}
\end{cor}
\begin{proof}
  Since, for any subgroup $Y \subseteq H^1(K, M)$, we have
  \[
    \widehat{\1}_Y = \frac{\size{Y}}{\size{H^0(K, M)}} \1_{Y^\perp},
  \]
  the result follows from Proposition \ref{prop:levels_new}\ref{lev:ur}--\ref{lev:size_all} and Lemma \ref{lem:Tate_level}.
\end{proof}

\subsection{The corresponding \'etale algebra}

Attached to any coclass $\sigma \in H^1(K, M)$ is a degree $p$ \'etale algebra $L$ over $K$, together with the structure that identifies its $\Aut(M)$-resolvent with the \'etale algebra encoding the Galois action on $M$; see \cite[Theorem 4.5(c)]{OEtale}. For a cyclic module of order $3$, this correspondence admits an especially explicit description. If $T$ is the quadratic \'etale algebra corresponding to $M$ and $T'$ is the quadratic \'etale algebra corresponding to its Tate dual $M'$, ordinary Kummer theory over $T'$ identifies $H^1(K, M)$ with a norm-one quotient of $T'^\cross$, and the three coordinates of $L$ are obtained from the three conjugates of a cube root of the Kummer element. We record the resulting formula from \cite[Proposition 5.1]{OEtale}.

\begin{prop}[\cite{OEtale}, Proposition 5.1]\label{prop:Kummer_cubic} Let $K$ be a field with $\ch K \neq 2, 3$, and let $ M = C_3^{D}$ be a $K$-Galois module with underlying group $ C_3 $ with Galois action factoring through the quadratic extension $K[\sqrt{D}]$. Let $T' = K[\sqrt{-3D}]$. Then we have a group isomorphism
  \[
  H^1(K, M) \isom T'^{N=1}/(T'^{N=1})^3
  \]
  in which $\delta \in T'^{N=1}/(T'^{N=1})^3$ corresponds to the cubic extension $L/K$ that is the image of the $K$-linear map
  \begin{align*}
    \kappa : K \cross T' &\to (K^\alg)^3 \\
    (a,\xi) &\mapsto \(a + \tr_{(K^\alg)^2/K^\alg} \xi \omega \sqrt[3]{\delta}\)_{\omega \in ((K^\alg)^2)^{N=1}[3]},
  \end{align*}
  where $\sqrt[3]{\delta} \in (K^\alg)^2$ is chosen to have norm $1$, and $\omega$ ranges through the set
  \[
  ((K^\alg)^2)^{N=1}[3] = \{(1;1), (\zeta_3; \zeta_3^2), ( \zeta_3^2; \zeta_3)\}
  \]
  of cube roots of $1$ in $(K^\alg)^2$ of norm $1$.
\end{prop}

\begin{prop}[\textbf{explicit Kummer theory for cubic algebras}]\label{prop:Kummer_resolvent_cubic}
  Let $R$ be a quadratic \'etale algebra over $K$ ($\ch K \neq 3$), and let
  \[
  L = K + \kappa(R)
  \]
  be the cubic algebra of resolvent $R' = R \odot K[\mu_3]$ (the Tate dual of $R$) corresponding to an element $\delta \in R^{N=1}$ in Proposition \ref{prop:Kummer_cubic}, where
  \[
  \kappa(\xi) = \(\tr_{(K^\alg)^2/K^\alg} \xi \omega \sqrt[3]{\delta}\)_{\omega \in \((K^\alg)^2\)^{N=1}[3]} \in (K^\alg)^3
  \]
  so $\kappa$ maps $R$ bijectively onto the traceless plane in $L$. Then the index form of $L$ is given explicitly by
  \begin{equation}\label{eq:Kummer_rsv_cubic}
    \begin{aligned}
      f : L/K &\to \Lambda^2(L/K) \\
      \kappa(\xi)& \mapsto 3\sqrt{-3}\delta\xi^3 \wedge 1,
    \end{aligned}
  \end{equation}
  where we identify
  \[
  \Lambda^2 L/K \isom \Lambda^3 L \isom \Lambda^2 R' \isom R'/K \isom \sqrt{-3} \cdot R/K
  \]
  using the fact that $R'$ is the discriminant resolvent of $L$.
\end{prop}
\begin{proof}
  Direct calculation, after reducing to the case $K = K^\alg$, $\delta = 1$, $L = K \cross K \cross K$.
\end{proof}

\begin{lem}\label{lem:disc}
  Let $\sigma \in H^1(K, M)$ be a nontrivial coclass. Then the degree $p$ \'etale algebra $L$ corresponding to $\sigma$ is a field, and its discriminant ideal can be computed from the defect $d = d(\sigma)$ as follows:
  \begin{equation}
    \Disc(L/K) = \begin{cases}
      \bigl(\pi_K^{p e + p - 1}\bigr), & d = -1, \\
      \bigl(\pi_K^{p e + p - 1 - d}\bigr), & 0 \leq d < p e, \\
      (1), & d \in \{p e, \infty\}.
    \end{cases}
  \end{equation}
\end{lem}
\begin{proof}
  Suppose first that $L$ is not a field; we claim that $\sigma = 1$. Let $\psi_L : G_K \to \Hol(M)$ be the corresponding homomorphism to the holomorph of $M$ as in \cite[Theorem 4.5]{OEtale}. If $L$ is not a field, then the image of $\psi_L$ is a nontransitive subgroup $J$ of $\Hol(M)$. If $J = \{\id\}$, then $\sigma = 1$. Otherwise, let $\lambda(x) = gx + t$ be a nonidentity element of $J$ such that $g$ generates the image of $J$ in the cyclic group $\Aut(M) \isom \FF_p^\cross$. If $g = 1$, then $\lambda$ is a $p$-cycle and $J$ is transitive, a contradiction. So $\lambda$ has a unique fixed point. By postconjugation (which corresponds to adjusting $\alpha$ by a coboundary), we may assume that $t = 0$. Now $J$ contains the maps $\lambda^r(x) = g^r x$. If $J$ has any other element $\mu(x) = g'x + t'$ with $t' \neq 0$, then $g' = g^r$ for some $r$, and composing with $\lambda^{-r}$, we get a $p$-cycle, a contradiction. So the image of $\psi_L$ is contained in $\Aut M$ and thus $\sigma = 1$. Hence $L$ is a field.
  
  By \cite[Lemma 8]{DCD},
  \[
    v_K\bigl(\Disc(L/K)\bigr) = \begin{cases}
      \dfrac{v_F\bigl(\Disc(FL/F)\bigr)}{p - 1} + p - 2, & v_F\bigl(\Disc(FL/F)\bigr) > 0, \\
      0, & v_F\bigl(\Disc(FL/F)\bigr) = 0. 
    \end{cases}
  \]
  But $FL = F[\sqrt[p]{\delta}]$ has the Kummer element $\delta \in F^\cross / (F^\cross)^p$ corresponding to $\sigma$. Its discriminant can be computed in terms of the defect of $\delta$; see \cite[Lemmas 5--7]{DCD}.
\end{proof}

\subsection{Aligning notation in the tame case}

If $ K $ is a \emph{tame} local field, that is, $ \ch k_K \neq p $, we define the defect and level of a coclass $\sigma$ by
\[
  d(\sigma) = \ell(\sigma) = \begin{cases*}
    -1 & if $\sigma$ is ramified \\
    0 & if $\sigma \neq 1$ is unramified \\ 
    \infty & if $\sigma = 1$.
  \end{cases*}
\]
Note that $\sigma$ is unramified if and only if the corresponding Kummer element $\alpha \in E^\cross/(E^\cross)^p$ is a unit. It is not hard to show that Propositions \ref{prop:defect} and  \ref{prop:levels_new}, Lemma \ref{lem:Tate_level}, Corollary \ref{cor:levels}, and \ref{lem:disc} hold as stated, where we interpret $e$ as $e = v_K(p) = 0$. For instance, the nontrivial case $i = 0$ of Lemma \ref{lem:Tate_level} boils down to
\[
  H^1_\ur(K, M)^\perp = H^1_\ur(K, M'),
\]
a standard result \cite[Theorem 7.2.15]{NSW}. This setup allows us to avoid treating the tame case separately in the proofs of local reflection theorems.

\section{Quadratic forms by superdiscriminant}
\label{sec:quadratic}

We begin with the simplest Galois module $M \isom \ZZ/2\ZZ$, over a field $K$ of characteristic not $2$.

There are many full composed varieties whose point stabilizer is of order $2$, and the one we take is, to say the least, one of the more unexpected. The group $\GL_2$ acts on the space
\[
  V = \Sym^2(2) = \{ a x^2 + b x y + c y^2 : a,b,c \in \GG_a \}
\]
of binary quadratic forms in the natural way. Let $\Gamma$ be the algebraic subgroup, defined over $\ZZ$, of elements of a peculiar form:
\[
  \left\{
  \begin{bmatrix}
    u & 0 \\
    t & u^{-2}
  \end{bmatrix}
  : u \in \GG_m, t \in \GG_a \right\}.
\]
Abstractly, this group is a certain semidirect product of $\GG_a$ by $\GG_m$. As is not too hard to verify, the restriction of the action to $\Gamma$ has a single polynomial invariant, the \emph{superdiscriminant}
\[
  I := a D = a(b^2 - 4a c),
\]
and the nondegenerate geometric orbits are level sets for $I$.
Because $\Gamma$ fixes the point $[x : y] = [1 : 0]$, it is rather natural to write forms in $V$ inhomogeneously as $f(x) = a x^2 + b x + c$, as was done in stating Theorem \ref{thm:O-N_quad_Z}.

For simplicity, we work over a number field with class number $1$, although the same local reflection theorem yields over a general number field a theorem summing over the class group.

\begin{rem}
  The group $\Gamma$ is not reductive, that is, does not fit into the classical Dynkin-diagram parametrization for Lie groups. Non-reductive groups are decidedly in the minority within the whole context of using orbits to parametrize arithmetic objects, but they have occurred before: Altu\u{g}, Shankar, Varma, and Wilson \cite{ASVW} count $D_4$-fields using orbits of pairs of ternary quadratic forms under a certain nonreductive subgroup of $\GL_2 \cross \SL_3$.
\end{rem}

\begin{rem} \label{rem:quad}
  The zeta functions for inhomogeneous quadratics, though less famous than their homogeneous cubic analogues, have had a long history. Shintani \cite{ShintaniQuadratic}, building on related work of Siegel \cite{SiegelFE}, proved meromorphic continuation to $\CC^2$ with a Klein four group of functional equations (see also \cite{IbukiyamaSaitoII}). As show by Diamantis and Goldfeld \cite{DiamantisGoldfeldShintani}, Shintani's double zeta functions are also linear combinations of Mellin transforms of metaplectic Eisenstein series, and thereby they fit into the framework of Weyl group multiple Dirichlet series. They noted that this structure predicts additional ``hidden'' functional equations not visible from the prehomogeneous vector space standpoint. Specifically, the part of the Shintani series corresponding to odd discriminants $D$ can be expressed in terms of the Weyl group multiple Dirichlet series of the root system $A_2$ \cite[\textsection2.3]{Wen}, which admits a dihedral group of $12$ functional equations \cite{BlomerDouble, HiramotoShintani}. The story is unchanged over number fields \cite{TaniguchiDiscCubic, KTW19}: one must separate the primes dividing $\infty$ and $2$ into an exceptional set $S$ to get a Dirichlet series admitting the full $12$ functional equations. As pointed out by Blomer \cite[Appendix A.3; note that the arXiv version does not include the appendices]{KTW19}; one of the extra functional equations, the swap $(s,w) \mapsto (w,s)$, is a symmetry of the Dirichlet series itself related to quadratic reciprocity. Thus, Theorem \ref{thm:O-N_quad_local} strengthens this functional equation along the reflection line $s = w$ to allow a more natural treatment at $2$. It is therefore natural to ask whether this can be extended to a full group of $12$ functional equations for Shintani's original zeta functions (and Taniguchi's number field generalization), both the even and odd terms, thus realizing the hidden symmetries anticipated by Diamantis--Goldfeld.
\end{rem}

We build a composed variety structure $(\Gamma, V, X_I, X_{I,0})$ on the orbit
\[
  X_I = \{f \in V : I(f) = I\},
\]
taking $X_{I,0}$ to be the orbit of the basepoint
\[
  f_0 = I x^2 - \frac{1}{4I}
\]
of discriminant $1$. Then the rational orbits are parametrized by $D = b^2 - 4ac \in K^\cross/(K^\cross)^2$ consistent with the parametrization of their splitting fields via Kummer theory.

Now we introduce the lattices used in the local and global counts. If $K_v$ is a nonarchimedean local field and $\tau \in \OO_v$ divides $2$, then take the local lattice
\[
  \Lambda_{v,\tau} = \{ a x^2 + b x + c : a,c \in \OO_v, b \in \tau \OO_v \};
\]
it is easily seen that the subgroup of $\Gamma(K_v)$ preserving $\Lambda_{v,\tau}$ is precisely $\Gamma(\OO_v)$. If $K$ is a number field and $\tau \in \OO_K$ divides $2$, these local lattices 
\[
\Lambda_\tau = \prod_{v \mid \infty} V(K_v) \cross \prod_{v \nmid \infty} \Lambda_{v,\tau}V(K_v).
\]
We write
\[
  X_I(\Lambda_\tau) = X_I(K) \intsec \Lambda_\tau.
\]
For this set to be nonempty, we must have $\tau^2 \mid I$.

Our first local reflection theorem says that each of these lattices has a natural dual.

\begin{thm}[\textbf{``Local Quadratic O-N''}]\label{thm:O-N_quad_local}
Let $K$ be a nonarchimedean local field, $\ch K \neq 2$. For $0 \neq I \in \OO_K$ and $\tau \mid 2$, let
$\ell_{v,\tau,I}$ denote the local lattice-count function on $H^1(K,\ZZ/2\ZZ)$
attached to $X_I$ and the lattice $\Lambda_{v,\tau}$. Then
\begin{equation} \label{eq:O-N_quad_local}
  \widehat{\ell}_{v,\tau,I}
  = \size{\OO_K/\tau\OO_K} \cdot \ell_{v,2\tau^{-1},4\tau^{-4}I}.
\end{equation}
Equivalently, the two lattices are naturally dual in the sense of Definition
\ref{defn:dual}, with duality constant $\size{\OO_K/\tau\OO_K}$.
\end{thm}

Before proving this theorem, we deduce the global quadratic reflection theorem from it using the local-global reflection engine.

\begin{thm}[\textbf{``Quadratic O-N''}] \label{thm:O-N_quad}
  Let $K$ be a number field of class number $1$. For $\tau \in \OO_K$ a divisor of $2$ and $0 \neq I \in \tau^2 \OO_K$, let $X_I(\Lambda_\tau)$ be the set of binary quadratic forms $f = a x^2 + b x y + c y^2$ over $\OO_K$ with $\tau \mid b$ and
  \[
  I(f) = a(b^2 - 4a c) = I,
  \]
  a set with an action of the algebraic group
  \[
  \Gamma(\OO_K) =
  \left\{
  \begin{bmatrix}
    u & 0 \\
    t & u^{-2}
  \end{bmatrix}
  : u \in \OO_K^\cross, t \in \OO_K \right\}.
  \]
  Then for $\tau \mid 2$ and $0 \neq I \in \tau^2 \OO_K$,
  \[
  \sum_{f \in \Gamma\(\OO_K\)\bs X_{4\tau^{-4}I}(\Lambda_{2\tau^{-1}})}
  \frac{1}{\size{\Stab_{\Gamma(\OO_K)} f}}
  = \frac{2^{r_1(K)+r_2(K)}}{\size{N_{K/\QQ}(\tau)}}
  \sum_{\substack{f \in \Gamma\(\OO_K\)\bs X_I(\Lambda_\tau) \\
      \disc f > 0 \text{ at every real place}}}
  \frac{1}{\size{\Stab_{\Gamma(\OO_K)} f}}
  \]
  where $r_1(K)$ and $r_2(K)$ are the numbers of real and complex places of $K$, respectively.
\end{thm}

\begin{proof}
  The variety $X_I$ is full, hence Hasse. The adelic lattices obtained from
  $\Lambda_\tau$ and $\Lambda_{2\tau^{-1}}$ have class number one for $\Gamma$:
  the additive factor has strong approximation, while the multiplicative factor
  has the same class number as $K$, which is one by hypothesis. Hence Lemma
  \ref{lem:total_cl_num} identifies the normalized total class numbers in
  Theorem \ref{thm:main_compose} with the sums on either side of the
  theorem statement.

  At every finite place, Theorem \ref{thm:O-N_quad_local} gives the required
  local reflection theorem. The local duality constants are $c_v = 1$ for $v \nmid 2 \infty$, while
  for $v \mid 2$ we get constants $c_v = [\OO_v : \tau\OO_v]$, whose product is
  $\size{N_{K/\QQ}(\tau)}$. At a real place $v$, the weighting on the $V^{(1)}$ (right) side is
  $\1_{\{1\}}$, selecting the trivial class in $H^1(K_v,\ZZ/2\ZZ)$, equivalently,
  stipulating positive discriminant. Its Fourier transform is $\frac12$ times the
  constant weighting on the $V^{(2)}$ (left) side, so $c_v = 1/2$. At a complex place the
  cohomology group is trivial and the normalization by $\size{H^0(K_v,M)} = 2$
  likewise gives $c_v = 1/2$. Hence
  \[
    \left(\prod_v c_v\right)^{-1} = \frac{2^{r_1(K)+r_2(K)}}{\size{N_{K/\QQ}(\tau)}},
  \]
  Theorem \ref{thm:main_compose}, in which the product of the $c_v$ appears inversely,
  gives the stated factor.
\end{proof}
\begin{proof}[Proof of Theorem \ref{thm:O-N_quad_Z}]
  We specialize further to the case $K = \QQ$. We replace $\Gamma(\ZZ)$ by its index-$2$ subgroup, the group $\ZZ$ of translations. This merely doubles all orbit counts, and it acts freely on quadratics with nonzero discriminant, so we can suppress all mention of stabilizers for the clean statement given in Theorem \ref{thm:O-N_quad_Z}.
\end{proof}

\subsection{Examples}
Over $\ZZ$, it is not hard to compute all quadratics $f$ of a fixed superdiscriminant $I$. The leading coefficient $a$ must be a divisor of $I$ (possibly negative). Then, by replacing $x$ by $x + t$ where $t$ is an integer nearest to $-b/(2a)$, we can assume that $b$ lies in the window $-|a| < b \leq |a|$. We then seek $b$ in this window such that
\[
c = \frac{ab^2 - I}{4a^2}
\]
is an integer; in other words, we solve the quadratic congruence
\begin{equation} \label{eq:quad_ex}
  b^2 \equiv \frac{I}{a} \mod 4a,
\end{equation}
which is close to how Shintani originally defined the zeta functions for quadratic forms \cite{ShintaniQuadratic}.
\begin{examp}\label{examp:15}
  There are five quadratics of superdiscriminant $15$:
  \[
  \begin{tabular}{l|cccc}
    $f(x)$ & $q$ & $q^+$ & $q_2$ & $q_2^+$ \\ \hline
    ${-x^2} + x - 4 $ & $\checkmark$ & & &    \\
    $15x^2 + x      $ & $\checkmark$ & $\checkmark$ & & \\
    $15x^2 - x      $ & $\checkmark$ & $\checkmark$ & & \\
    $15x^2 + 11x + 2$ & $\checkmark$ & $\checkmark$ & & \\
    $15x^2 - 11x + 2$ & $\checkmark$ & $\checkmark$ & & 
  \end{tabular}
  \]
  So we get the totals
  \[
  q(15) = 5 \textand
  q^+(15) = 4.
  \]
  There are $18$ quadratics of superdiscriminant $60$:
  \[
  \begin{tabular}{l|cccc}
    $f(x)$ & $q$ & $q^+$ & $q_2$ & $q_2^+$ \\ \hline
    ${x^2 - 15        }$ & $\checkmark$ & $\checkmark$ & $\checkmark$ &  $\checkmark$ \\
    ${-x^2 - 15       }$ & $\checkmark$ & & $\checkmark$ & \\
    ${-3x^2 + 2x - 2  }$ & $\checkmark$ & & $\checkmark$ & \\
    ${-3x^2 - 2x - 2  }$ & $\checkmark$ & & $\checkmark$ & \\
    ${-4x^2 + x - 1   }$ & $\checkmark$ & & & \\
    ${-4x^2 - x - 1   }$ & $\checkmark$ & & & \\
    ${15x^2 + 2x      }$ & $\checkmark$ & $\checkmark$ & $\checkmark$ &  $\checkmark$ \\
    ${15x^2 - 2x      }$ & $\checkmark$ & $\checkmark$ & $\checkmark$ &  $\checkmark$ \\
    ${15x^2 + 8x + 1  }$ & $\checkmark$ & $\checkmark$ & $\checkmark$ &  $\checkmark$ \\
    ${15x^2 - 8x + 1  }$ & $\checkmark$ & $\checkmark$ & $\checkmark$ &  $\checkmark$ \\
    ${60x^2 + x       }$ & $\checkmark$ & $\checkmark$ & & \\
    ${60x^2 - x       }$ & $\checkmark$ & $\checkmark$ & & \\
    ${60x^2 + 31x + 4 }$ & $\checkmark$ & $\checkmark$ & & \\
    ${60x^2 - 31x + 4 }$ & $\checkmark$ & $\checkmark$ & & \\
    ${60x^2 + 41x + 7 }$ & $\checkmark$ & $\checkmark$ & & \\
    ${60x^2 - 41x + 7 }$ & $\checkmark$ & $\checkmark$ & & \\
    ${60x^2 + 49x + 10}$ & $\checkmark$ & $\checkmark$ & & \\
    ${60x^2 - 49x + 10}$ & $\checkmark$ & $\checkmark$ & &
  \end{tabular}
  \]
  In particular,
  \[
  q_2^+(60) = 5 = q(15) \textand q_2(60) = 8 = 2\cdot 4 = 2q^+(15)
  \]
  in accordance with Theorem \ref{thm:O-N_quad_Z}. \ignore{From the same theorem, we derive, without computation, that
  \[
  q_2^+(240) = q(60) = 18 \textand q_2(240) = 2q^+(60) = 26.
  \]}
\end{examp}
\begin{examp} \label{ex:QR}
  In general, the solubility of \eqref{eq:quad_ex} can be expressed in terms of Legendre symbols. For a superdiscriminant $I = p_1p_3$, where $p_1 \equiv 1$ (mod $4$) and $p_3 \equiv 3$ (mod $4$) are primes, the combination $a = p_3$, $D = p_1$ is feasible if and only if the congruence
  \[
  b^2 \equiv p_1 \mod 4p_3
  \]
  has a solution, which happens exactly when $\( \dfrac{p_1}{p_3} \) = 1$. Examining all factorizations of $I$ and $4I$, we find that
  \[
  q^+(p_1p_3) = 5 + \( \frac{p_1}{p_3} \) \textand
  q_2(4p_1p_3) = 10 + 2 \( \frac{p_3}{p_1} \).
  \]
  Thus our reflection theorem recovers the quadratic reciprocity law
  \[
  \( \frac{p_1}{p_3} \) = \( \frac{p_3}{p_1} \).
  \]
\end{examp}

\begin{qn}
  Does there exist a proof of Theorem \ref{thm:O-N_quad_Z} using no tools more advanced than quadratic reciprocity in $\ZZ$? In the prime-to-$2$ case $2 \nmid I$, the corresponding swap symmetry in the $S = \{\infty,2\}$ theory of Kim--Tsuzuki--Wakatsuki is deduced from Blomer's double Dirichlet series by quadratic reciprocity \cite[Appendix A.3]{KTW19}. Such a proof would extend that reciprocity argument through the wild dyadic contribution rather than removing it.
\end{qn}

\subsection{Proof of local quadratic O-N}

Theorem \ref{thm:O-N_quad_local} is an equality of functions on
\[
  H^1(K, C_2) \isom K^\cross / (K^\cross)^2.
\]
We filter $H^1(K, C_2)$ by the level spaces $\L_i$ of Proposition \ref{prop:levels_new}. Since $p = 2$, there are simplifications (for instance, the congruence \eqref{eq:dp1} is vacuous), and the level spaces $\L_{-1}, \L_0, \ldots, \L_{e+1}$ are
\[
\L_i = \begin{cases}
  K^\cross/(K^\cross)^2, & i = -1 \\
  \{[\alpha] \in \OO_K^\cross/(\OO_K^\cross)^2 : \alpha \in 1 + \pi^{2e-2i} \OO_K\}, & 0 \leq i \leq e \\
  \{1\}, & i = e+1.
\end{cases}
\]
The module $M = C_2$ is self-dual, and the Tate pairing on $H^1(K, C_2)$ is none other than the quadratic Hilbert pairing on $K^\cross/(K^\cross)^2$, by Lemma \ref{lem:Tate_Hilbert}.

Let $L_i$ be the characteristic function of $\L_i$. By Corollary \ref{cor:levels}, the Fourier transform of each $L_i$ is a scalar multiple of $L_{e - i}$. We will prove Theorem \ref{thm:O-N_quad_local} by writing the local lattice-count function as a linear combination of $L_i$ and showing that the coefficients of $L_i$ on each side of \ref{eq:O-N_quad_local} match.

\begin{proof}[Proof of Theorem \ref{thm:O-N_quad_local}]
We proceed by explicitly computing the local lattice-count function $\ell_{v,\tau,I}$, which sends each $[D] \in K^\cross/(K^\cross)^2 \isom H^1(K, \ZZ/2\ZZ)$ to the number of cosets $[\gamma] \in \Gamma(\OO_K)\bs \Gamma(K)$ such that $\gamma v_0 \in \Lambda_{v,\tau}$, where $v_0$ is an arbitrary vector in $X_I(K)$ and $D(v_0) = D$. Let $t = v(\tau)$ and $e = v(2)$; we have
$e > 0$ exactly when $K$ is $2$-adic, and
\[
  0 \leq t \leq e.
\]

A coset $[\gamma]$ is specified by two pieces of information. First is the valuation $v(u)$ of the diagonal element $u$; this is equivalent to specifying $v(D)$ and $v(a)$, where
\[
  \gamma v_0 = a x^2 + b x y + c y^2 \textand D = b^2 - 4 a c.
\]
Second, we specify $t$ modulo the appropriate integral sublattice. If (as we may assume) $v(u) = 0$, then $t$ must be specified modulo $1$, which is the same as specifying $b$ modulo $2a$. So the problem of computing $\ell_{v,\tau,I}$ devolves onto computing how many $b \in \tau\OO_K$, up to translation by $2a$, yield an integral value for
\[
  c = \frac{b^2 - D}{4a};
\]
that is, we must solve the quadratic congruence (a local version of \eqref{eq:quad_ex})
\begin{equation} \label{eq:qfc}
  b^2 \equiv D \mod 4a.
\end{equation}
The answer, in general, depends on how close $D$ is to being a square in $K$, which is closely connected to the level of $D$ as an element of $K^\cross/(K^\cross)^2$.

We claim that, if we fix $v(I)$ and $v(a)$ (and hence $v(D)$), then the contribution of all solutions of \eqref{eq:qfc} to $\ell_{v,\tau,I}$ can be expressed as a linear combination of the $L_i$. The basic idea, which will recur for other reflection theorems, is to group the solutions into \emph{families} that have a constant number of solutions over some subset $S \subseteq K^\cross/(K^\cross)^2$. The subset $S$ will be called the \emph{support} of the family, and the number of solutions for each $D \in S$ will be called the \emph{thickness} of the family.

We divide the possible pairs $\bigl(v(a), v(D)\bigr)$ into three cases, called \emph{zones}, delimited by linear inequalities, depicted in Figure \ref{fig:zones_quad}. (The $v(D)$ axis shows only values with $v(D) \geq 2 t$ due to the restriction $b \in \tau\OO_K$.) The families in each zone have a unified description.

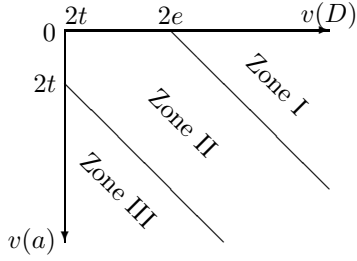
\begin{figure}
  \setlength{\unitlength}{1em}
  \[
  \begin{picture}(13,9)(-2,-8.2)
    \put(0,0){\vector(1,0){10}}
    \put(0,0){\vector(0,-1){8}}
    \put(0,-2){\line(1,-1){6}}
    \put(4,0){\line(1,-1){6}}
    \put(-2,-9){\makebox(2,2)[r]{$v(a)$~}}
    \put(-2,-3){\makebox(2,2)[r]{$2t$~}}
    \put(-2,-1){\makebox(2,2)[r]{$0$\vphantom{$^2$}~}}
    \put( 0,0){\makebox(2,2)[bl]{$2t$\vphantom{|}}}
    \put( 3,0){\makebox(2,2)[b]{$2e$\vphantom{|}}}
    \put( 9,0){\makebox(2,2)[b]{$v(D)$\vphantom{\big|}}}
    \put(2,-1){\rotatebox{-45}{\makebox(7,0){Zone II}}}
    \put(6,0){\rotatebox{-45}{\makebox(5.7,0){Zone I}}}
    \put(0,-4){\rotatebox{-45}{\makebox(5.7,0){Zone III}}}
  \end{picture}
  \]
  \caption{Three zones used in the local counting of quadratic forms by superdiscriminant}
  \label{fig:zones_quad}
\end{figure}

If $v(D) \geq v(4a)$ (Zone I), then \eqref{eq:qfc} simplifies to $4a|b^2$, that is,
\[
  v(b) \geq \ceil{e + \frac{1}{2}v(a)}.
\]
Since we are counting values of $b$ modulo $2a$, the number of solutions is
\[
  q^{\( e + v(a)\) - \ceil{e + \frac{1}{2} v(a)}} = q^{\floor{\frac{1}{2}v(a)}}.
\]
We thus get a family with this thickness, supported on either $\L_0$ or $\L_{-1} \setminus \L_{0}$ according as $v(D)$ is even or odd.

If $v(D) < v(4a)$, then $b^2$ must be actually able to cancel at least the leading term of $D$ to get any solutions. In particular, $v(D)$ must be even. Let $\tilde{D} = D / \pi^{v(D)}$, and let $\tilde{b} = b/\pi^{\frac{1}{2} v(D)}$, so $\tilde{b}$ must be a unit satisfying
\begin{equation} \label{eq:qfc2}
  \tilde b^2 \equiv \tilde D \mod \frac{4a}{D}.
\end{equation}
Let $m = v(4a/D)$. If $m \geq 2e+1$ (Zone III), a unit is a square modulo $\pi^m$ only if it is a square outright, so we get a family supported just on the trivial class $1 \in K^\cross/(K^\cross)^2$. Otherwise, we have $1 \leq m \leq 2e$ (Zone II), and the support is $L_{\ceil{m/2}}$. The corresponding thicknesses are easy to compute. The $\tilde b$ satisfying \eqref{eq:qfc2} form a fiber of the group homomorphism
\[
  \phi = \bullet^2 \colon \( \OO_K / \pi^{v(2a) - \frac{1}{2}v(D)} \)^\cross
  \to \( \OO_K / \pi^{v(4a) - v(D)} \)^\cross,
\]
and the cokernel of this homomorphism has size $[\L_0 : \L_i]$, so the thickness is
\begin{align*}
  \size{\ker \phi} &= [\L_0 : \L_i] \cdot \frac{\Size{\(\OO_K / \pi^{v(2a) - \frac{1}{2}v(D)} \)^\cross}}{\Size{\(\OO_K / \pi^{v(4a) - v(D)} \)^\cross}} \\
    &= [\L_0 : \L_i] \cdot \frac{\( 1 - \frac{1}{q} \) q^{v(2a) - \frac{1}{2}v(D)}}
  {\( 1 - \frac{1}{q} \) q^{v(4a) - v(D)}} \\
  &= [\L_0 : \L_i] \cdot q^{\frac{1}{2} v(D) - e} \\
  &= \begin{cases*}
    q^{\frac{1}{2} v(D) - e + \ceil{m/2}} = q^{\floor{v(a)/2}} & (Zone II) \\
    2q^{\frac{1}{2} v(D)} & (Zone III).
  \end{cases*}
\end{align*}

By way of illustration, we tabulate the contributions to $\ell_{v,\tau,I}$ when $e = 2$ in Table \ref{tab:quad_ex}.
\begin{table}
  \[
\newcommand{\rline}[1]{\multicolumn{1}{c|}{#1}}
\begin{tabular}{r|ccccccccc}
  $\downarrow v(a); v(D) \rightarrow$ & $0$ & $1$ & $2$ & $3$ & $4$ & $5$ & $6$ & $7$ & $8$ \\ \hline
  $0$ & ${L_2}$ & & ${L_1}$ & \rline{} & ${L_0}$ & ${L_{-1}} - {L_0}$ & ${L_0}$ & ${L_{-1}} - {L_0}$ & ${L_0}$ \\ 
  \cline{2-2}\cline{6-6}
  $1$ & \rline{${L_3}$} & & ${L_1}$ & & \rline{${L_0}$} & ${L_{-1}} - {L_0}$ & ${L_0}$ & ${L_{-1}} - {L_0}$ & ${L_0}$ \\
  \cline{3-3}\cline{7-7}
  $2$ & ${L_3}$ & \rline{} & $q{L_2}$ & & $q{L_1}$ & \rline{} & $q{L_0}$ & $q\({L_{-1}} - {L_0}\)$ & $q{L_0}$ \\
  \cline{4-4}\cline{8-8}
  $3$ & ${L_3}$ & & \rline{$q{L_3}$} & & $q{L_1}$ & & \rline{$q{L_0}$} & $q\({L_{-1}} - {L_0}\)$ & $q{L_0}$ \\
  \cline{5-5}\cline{9-9}
  $4$ & ${L_3}$ & & $q{L_3}$ & \rline{} & $q^2{L_2}$ & & $q^2{L_1}$ & \rline{} & $q^2{L_0}$ \\
  \cline{6-6}\cline{10-10}
\end{tabular}
\]
\caption{Families of quadratics over a $2$-adic local field with absolute ramification index $e = 2$. The staircase-shaped lines delimit Zones I, II, and III (from top to bottom)}
\label{tab:quad_ex}
\end{table}
We have shown the subdivision of the table into three \emph{zones} given by the inequalities:
\begin{itemize}
  \item Zone I: $v(D) \geq v(4a)$
  \item Zone II: $v(a) \leq v(D) < v(4a)$
  \item Zone III: $v(D) < v(a)$.
\end{itemize}
The graphs of these zones are shown in Figure \ref{fig:zones_quad} (note the necessary condition that $v(D) \geq 2t$).

The feature to be noted is that, under the transformation
\[
  t \leftrightarrow e - t, \quad v(a) \leftrightarrow v(D) - 2t,
\]
the shape of Zone II is flipped about a diagonal line and Zones I and III are interchanged. This will be the basis for our proof of Theorem \ref{thm:O-N_quad_local}; but with slight irregularities owing to the floor functions in the formulas and the fact that $\L_{-1} \setminus \L_0$, instead of $\L_{-1}$, appears as a support.

There are two ways to finish the proof. One is to establish a bijection of families, as outlined in the previous paragraph, so that $L_i$ and $L_{e-i}$ are interchanged as supports and all the thicknesses correspond appropriately. Such an approach will be used for cubic O-N in Section \ref{sec:wild_bij}. The other is to verify the local reflection computationally, by means of a generating function. We present the second method. Here the expressions are simple enough to fit on the page; with the aid of a computer, the same method is adaptable to proving other local reflection theorems.

Let
\[
  F_{t,e}(Z) = \sum_{n \geq 0} \ell_{v,\tau,\pi^n} Z^n,
\]
a formal power series whose coefficients are functions of $D \in K^\cross/(K^\cross)^2$. The desired theorem can be written as
\[
  \widehat F_{t,e} = q^{e-t} F_{e-t,e}
\]
where the Fourier transform is taken with respect to $D$ only, and the equality of coefficients of $Z^k$ encodes the claim that Theorem \ref{thm:O-N_quad_local} holds for $v(I) = k$.

We write $F = F_{\zoneI} + F_{\zoneII} + F_{\zoneIII}$, where $F_{\X}$ is the contribution coming from Zone $\X$ in the preceding analysis. Writing $i = v(a)$ and $d = v(D)$, we proceed to compute each $F_{\X}$ in turn:
\begin{align*}
  F_{\zoneI} &= \sum_{i \geq 0} \sum_{d \geq 2e+i} \begin{cases}
    q^{\floor{i/2}} L_0 Z^{i+d}, & d \text{ even} \\
    q^{\floor{i/2}} (L_{-1} - L_0) Z^{i+d}, & d \text{ odd}
  \end{cases}\\
  &= Z^{2e} \sum_{i \geq 0} \sum_{d \geq i} \begin{cases}
    q^{\floor{i/2}} L_0 Z^{i+d}, & d \text{ even} \\
    q^{\floor{i/2}} (L_{-1} - L_0) Z^{i+d}, & d \text{ odd}.
  \end{cases} \\
\intertext{Splitting $i = 2i_f + i_p$, where $0 \leq i_p \leq 1$, and likewise $d = 2d_f + d_p$, we get}
  F_{\zoneI} &= Z^{2e} \sum_{i_p = 0}^1 \sum_{i_f \geq 0} \left(
    \sum_{d \geq 2i_f + i_p} q^{i_f} (-1)^d L_0 Z^{d + 2i_f + i_p} + \sum_{d_f \geq i_f} q^{i_f} L_{-1} Z^{2d_f + 2i_f + i_p + 1} \right) \\
  &= Z^{2e} \sum_{i_p = 0}^1 \sum_{i_f \geq 0} \left(\frac{q^{i_f}(-1)^{i_p} L_0 Z^{4i_f + 2i_p}}{1 + Z} +
    \frac{q^{i_f} Z^{4i_f + i_p + 1}L_{-1}}{1 - Z^2} \right) \\
  &= Z^{2e} \sum_{i_p = 0}^1  \left( \frac{(-1)^{i_p} Z^{2i_p} L_0}{(1+Z)(1 - qZ^4)} + \frac{Z^{i_p + 1} L_{-1}}{(1 - Z^2)(1 - qZ^4)} \right) \\
  &= Z^{2e}\(\frac{(1-Z^2)}{(1+Z)(1 - qZ^4)} L_0 + \frac{Z(1 + Z)}{(1-Z^2)(1 - qZ^4)} L_{-1} \) \\
  &= \frac{Z^{2e}(1-Z)}{1 - qZ^4} L_0 + \frac{Z^{2e+1}}{(1-Z)(1 - qZ^4)} L_{-1}.
\end{align*}
For Zone II, which appears only when $e > 0$, we simplify the sum
\[
  F_{\zoneII} = \sum_{i\geq 0} \sum_{\substack{i \leq d < i + 2e \\ d \geq 2t \\ d\text{ even}}} q^{e + \floor{\frac{i}{2}} - \frac{d}{2}} L_{e + \floor{\frac{i}{2}} - \frac{d}{2}} Z^{i+d}
\]
by grouping terms with the same level $L_j$. We have $j = e + \floor{\frac{i}{2}} - \frac{d}{2}$, so the values of $i$ and $j$ determine $d$. The condition $d \geq 2t$ reduces to $i \geq 2(j + t - e)$; the other condition $i \leq d < i + 2e$ is automatically satisfied if $1 \leq j \leq e-1$, while if $j = 0$ or $j = e$, we must have $i$ odd or $i$ even respectively. For $1 \leq j \leq e-1$, the $L_j$-piece of $F_{\zoneII}$ is therefore
\begin{align*}
&\sum_{i \geq \max\{0, 2(j + t - e)\}} q^{\floor{i/2}} L_j Z^{i + 2\( e + \floor{\frac{i}{2}} - 2j\)} \\
&= \sum_{i_p = 0}^1 \sum_{i_f \geq \max\{0, j + t - e\}} q^{i_f} Z^{4i_f + i_p + 2e - 2j} L_j \\
&= \left( \sum_{i_p = 0}^1 Z^{i_p} \right) \left( \sum_{i_f \geq \max\{0, j + t - e\}} \( qZ^4\)^{i_f} \right) Z^{2e - 2j} L_j \\
&= \frac{(1+Z)\( qZ^4\)^{\max\{0,j+t-e\}}}{1 - qZ^4} \cdot Z^{2e - 2j}L_j.
\end{align*}
For $j = 0$ and $j = e$, since $i_p$ can only take one of its two values, the initial factor $1+Z$ is to be replaced by $Z$ and $1$ respectively, so
\[
  F_{\zoneII} = \frac{Z^{2e+1}}{1 - qZ^4} L_0 + \sum_{1 \leq j \leq e-1} \frac{(1+Z)\( qZ^4\)^{\max\{0,j+t-e\}}}{1 - qZ^4} \cdot Z^{2e - 2j}L_j + \frac{\( qZ^4 \)^t}{1 - qZ^4} L_e.
\]
Finally, Zone III presents no added difficulties:
\begin{align*}
  F_{\zoneIII} &= \sum_{\substack{d\geq 2t \\ d\text{ even}}} \sum_{i \geq d + 1} 2q^{d/2} L_{e+1} Z^{i+d} \\
  &= 2 \sum_{\substack{d\geq 2t \\ d\text{ even}}} q^{d/2} \cdot \frac{Z^{2d+1}}{1 - Z} \cdot L_{e+1} \\
  &= \frac{2q^tZ^{4t+1}}{(1 - Z)(1 - qZ^4)} L_{e+1}.
\end{align*}
Summing up, we get for $e \geq 1$ (the case $e = 0$ can be handled similarly)
\begin{align*}
  F_{t,e} &= F_{\zoneI} + F_{\zoneII} + F_{\zoneIII} \\
  &= \frac{Z^{2e+1}}{(1-Z)(1 - qZ^4)} L_{-1} + \frac{Z^{2e}(1-Z)}{1 - qZ^4} L_0
  + \frac{Z^{2e+1}}{1 - qZ^4} L_0 + \sum_{1 \leq j \leq e-1} \frac{(1+Z)\( qZ^4\)^{\max\{0,j+t-e\}}}{1 - qZ^4} \cdot Z^{2e - 2j}L_j \\ & \quad {} + \frac{\( qZ^4 \)^t}{1 - qZ^4} L_e + \frac{2q^tZ^{4t+1}}{(1 - Z)(1 - qZ^4)} L_{e+1} \\
  &= \frac{Z^{2e+1}}{(1-Z)(1 - qZ^4)} L_{-1} + \frac{Z^{2e}}{1 - qZ^4} L_0 + \sum_{1 \leq j \leq e-1} \frac{(1+Z)\( qZ^4\)^{\max\{0,j+t-e\}}}{1 - qZ^4} \cdot Z^{2e - 2j}L_j \\ & \quad {} + \frac{\( qZ^4 \)^t}{1 - qZ^4} L_e + \frac{2Z \( qZ^4 \)^t}{(1 - Z)(1 - qZ^4)} L_{e+1}.
\end{align*}
We now note that for $-1 \leq j \leq e+1$,
\[
  \coef_{L_j} F_{t,e} = \frac{\size{L_j}}{2 q^t} \coef_{L_{e-j}} F_{e-t,e},
\]
so
\[
  \widehat F_{t,e} = q^{e-t} F_{e-t,e}
\]
as desired.
\end{proof}

\section{Cubic Ohno-Nakagawa}
\label{sec:cubic}
\subsection{Statements of results}
The space $V(K)$ of binary cubic forms over a local or global field $K$ can be equipped with many invariant lattices. Let $V_{\OO_K}$ be the lattice of binary cubic forms with trivial Steinitz class; these can be written as
\[
V_{\OO_K} = \{a x^3 + b x^2 y + c x y^2 + d y^3 : a,b,c,d \in \OO_K\}.
\]
In this subsection we abbreviate the form $a x^3 + b x^2 y + c x y^2 + d y^3$ to $(a,b,c,d)$. A theorem of Osborne classifies all lattices $L \subseteq V_{\OO_K}$ that are $\GL_2(\OO_K)$-invariant and \emph{primitive}, in the sense that $\pp^{-1}L \nsubseteq V(\OO_K)$ for all finite primes $\pp$ of $\OO_K$ (any lattice can be scaled by a unique fractional ideal to become primitive):
\begin{prop}[Osborne \cite{Osborne}, Theorem 2.2] \label{prop:Osborne}
  A primitive $\GL_2(\OO_K)$-invariant lattice in $V(\OO_K)$ is given by any combination of the primitive $\GL_2(\OO_{K,\pp})$-invariant lattices in the completions $V(\OO_{K,\pp})$, which are:
  \begin{enumerate}[(a)]
    \item\label{it:3} If $\pp|3$, the lattices $\Lambda_{\pp,i} = \{(a,b,c,d) : b \equiv c \equiv 0 \bmod \pp^i \}$, for $0 \leq i \leq v_\pp(3)$;
    \item\label{it:2} If $\pp|2$ and $N_{K/\QQ}(\pp) = 2$, the five lattices
    \begin{align*}
      \Lambda_{\pp,1} &= V(\OO_{K,\pp}), \\
      \Lambda_{\pp,2} &= \{(a, b, c, d) \in V(\OO_{K,\pp}) : a + b + d \equiv a + c + d \equiv 0 \mod \pp \} \\
      \Lambda_{\pp,3} &= \{(a, b, c, d) \in V(\OO_{K,\pp}) : a + b + c \equiv b + c + d \equiv 0 \mod \pp \} \\
      \Lambda_{\pp,4} &= \{(a, b, c, d) \in V(\OO_{K,\pp}) : b + c \equiv 0 \mod \pp \} \\
      \Lambda_{\pp,5} &= \{(a, b, c, d) \in V(\OO_{K,\pp}) : a \equiv d \equiv b + c \mod \pp \},
    \end{align*}
    \item For all other $\pp$, the maximal lattice $V(\OO_{K,\pp})$ only.
  \end{enumerate}
\end{prop}
For $\pp | 2$, the latter four lattices are not $\GL_2$-invariant as integral forms of the affine space $V$, because they lose their $\GL_2$-invariance as soon as we extend scalars so that the residue field has more than two elements. By contrast, if $\pp|3$, the $\GL_2$-invariance of the space $\Lambda_{\pp,i}$ can be established purely formally. We establish a local reflection theorem on this lattice, which we will call the space of \emph{$\pp^i$-traced} forms due to a connection to the relative trace map on the associated cubic ring (see below).

Although Osborne deals only with the case of $V(\OO_K)$, his method generalizes easily to the lattice
\[
  V(\OO_K,\aa) = \{a x^3 + b x^2 y + c x y^2 + d y^3 : a \in \aa, b \in \OO_K, c \in \aa^{-1}, d \in \aa^{-2}\}
\]
that comprises the index forms
\[
  f : M \to \Lambda^2 M, \qquad M \coloneqq \OO_K \oplus \aa
\]
of cubic rings of Steinitz class $\aa$ in Theorem \ref{thm:hcl_cubic_ring}. Here 
\[
  \GL(M) =
  \left\{\begin{bmatrix}
    a_{11} & a_{12} \\
    a_{21} & a_{22}
  \end{bmatrix} \in \GL_2(K): a_{ij} \in \aa^{j-i}\right\}
\]
naturally acts on both $M$ and $\Lambda^2 M$, thus affecting $V(\OO_K, \aa)$ via a twisted action
\begin{equation}\label{eq:twisted_action}
\left(\begin{bmatrix}
  a_{11} & a_{12} \\
  a_{21} & a_{22}
\end{bmatrix}
\mathop{.\vphantom{I}} f\right)(x,y) = \frac{1}{a_{11} a_{22}-a_{12} a_{21}} f(a_{11} x + a_{21} y, a_{12} x + a_{22} y).
\end{equation}
(Compare Cox \cite{potf}, p.~142 and Wood \cite{WGauss}, Theorem 1.2.) For our purposes, the denominator of \eqref{eq:twisted_action} can be ignored since we only need the action of the group $\SL(M)$, which preserves the discriminant $D \in \aa^{-2}$ of the form. Here and below, the discriminant $\Disc(\OO)$ of a cubic ring over $\OO_K$ means the discriminant ideal of its trace pairing. Under the parametrization of Theorem \ref{thm:hcl_cubic_ring}, if the index form has discriminant $D \in \aa^{-2}$, then
\[
  \Disc(\OO) = D\aa^2.
\]

For instance, over $K = \ZZ$ there are ten primitive invariant lattices, comprising five types at $2$ and two types at $3$. The O-N-like reflection theorems relating these ten lattices were computed by Ohno and Taniguchi \cite{10lat}; over a number field, such an analysis can be more easily folded into the paradigm of splitting type selectors in Section \ref{sec:non-natural}. On the other hand, the behavior at $3$ generalizes elegantly:

\begin{thm}[\textbf{``Local cubic O-N''}] \label{thm:O-N_cubic_local}
  Let $K$ be a nonarchimedean local field, $\ch K \neq 2, 3$, and let $X_D$ be the
  composed variety of binary cubic forms of nonzero discriminant $D$ under
  $\Gamma=\SL_2$. For $\alpha \in K^\cross$ and $\tau \mid 3$, let
  \[
    \Lambda_{v,\alpha,\tau}
    = \{a\alpha x^3+b\tau x^2y+c\alpha^{-1}\tau xy^2+d\alpha^{-2}y^3:
       a,b,c,d\in\OO_K\}.
  \]
  Write $\ell_{v,\alpha,\tau,D}$ for the local lattice-count function attached
  to $X_D$ and $\Lambda_{v,\alpha,\tau}$. Then
  \begin{equation}\label{eq:x_cubic_dual_t}
    \widehat{\ell}_{v,1,\tau,D}
    = \size{\OO_K/\tau\OO_K}\,
      \ell_{v,1,3\tau^{-1},-27\tau^{-6}D},
  \end{equation}
  and, after the corresponding change of basis,
  \begin{equation}\label{eq:x_cubic_dual_a_t}
    \widehat{\ell}_{v,\alpha,\tau,D}
    = \size{\OO_K/\tau\OO_K}\,
      \ell_{v,\alpha\tau^{-3},3\tau^{-1},-27D}.
  \end{equation}
  Thus, when $\Lambda_{v,\alpha,\tau}$ is taken as the first lattice and
  $\Lambda_{v,\alpha\tau^{-3},3\tau^{-1}}$ as the second, the local duality constant
  of Definition \ref{defn:dual} is $\size{\OO_K/\tau\OO_K}$.
\end{thm}

The two formulations \eqref{eq:x_cubic_dual_t}, \eqref{eq:x_cubic_dual_a_t} are equivalent under replacing $D$ by $\alpha^2 D$ and changing basis by the twisted action \eqref{eq:twisted_action} of the matrices
\[
  \begin{bmatrix}1 & \\ & \alpha^{-1}\end{bmatrix}, \qquad
  \begin{bmatrix}1 & \\ & \alpha^{-1} \tau^{3}\end{bmatrix} \quad \in \GL_2(K)
\]
on the binary cubic forms counted on the left and right sides of \eqref{eq:x_cubic_dual_t} respectively.
The formulation \eqref{eq:x_cubic_dual_t} is the one we will prove, but \eqref{eq:x_cubic_dual_a_t} has the needed form of a local reflection theorem to apply at each place to get Theorem \ref{thm:O-N_traced_intro}:
\begin{proof}[Proof of Theorem \ref{thm:O-N_traced_intro}]
  Fix a Steinitz class $\aa$. Apply Theorem \ref{thm:main_compose} to the
  adelic lattice whose finite components are the lattices
  $\Lambda_{v,\alpha_v,\tau_v}$ of Theorem \ref{thm:O-N_cubic_local}, where $\alpha_v$ and $\tau_v$ generate the completed ideals $\aa_v$ and $\tt_v$, respectively. At the infinite
  places the local cohomology is trivial, and the ratio of the normalization factors in \eqref{eq:main_compose} gives the stated archimedean constant. Lemma \ref{lem:total_cl_num}
  then converts the normalized total class numbers into the desired weighted orbit counts. Since $\SL_2$ has class number one (strong approximation) in the number field case, we get one term on each side of the equation as claimed.
\end{proof}

If $\OO$ is a ring of finite rank over a Dedekind domain $\OO_K$, define its \emph{trace ideal} $\tr(\OO)$ to be the image of the trace map $\tr_{\OO/\OO_K} : \OO \to \OO_K$. Note that $\tr(\OO)$ is an ideal of $\OO_K$ and, since $1 \in \OO$ has trace $n = \deg(\OO/\OO_K)$, it is a divisor of the ideal $(n)$. In particular, if $n$ is a unit in $\OO_K$, this notion is uninteresting. Let $\tt$ be an ideal of $\OO_K$ dividing $(n)$. We say that the ring $\OO$ is \emph{$\tt$-traced} if $\tr(\OO) \subseteq \tt$.

By Theorem \ref{thm:hcl_cubic_ring}, we can parametrize cubic orders $\OO$ by their Steinitz class $\aa$ and index form
\[
\Phi(x\xi + y\eta) = (a x^3 + b x^2 y + c x y^2 + d y^3)(\xi \wedge \eta)
\]
relative to a decomposition $\OO = \OO_K \oplus \OO_K\xi \oplus \aa\eta$, where $a \in \aa$, $b \in \OO_K$, $c \in \aa^{-1}$, and $d \in \aa^{-2}$. Then a short computation using the multiplication table from Theorem \ref{thm:hcl_cubic_ring} shows that, if $(1, \xi, \eta)$ is a normal basis, then $\tr(\xi) = b$ and $\tr(\eta) = -c$, so $\tr(\OO) = \<3, b, \aa c\>$. Thus the based $\tt$-traced rings over $\OO_K$ are parametrized by the rank-$4$ lattice of cubic forms
\[
\V_{\aa,\tt}(\OO_K) := \{a x^3 + b x^2 y + c x y^2 + d y^3 : a \in \aa, b \in \tt, c \in \tt\aa^{-1}, d \in \aa^{-2}\},
\]
on which $\GL(\OO \oplus \aa)$ acts by the twisted action \eqref{eq:twisted_action}. For instance, if $\OO_K = \ZZ$, $\aa = (1)$, and $\tt = (3)$, this is the lattice of integer-matrix cubic forms considered in the introduction. Our goal in this section is to prove a generalization for all number fields $K$ and spaces $V(\OO_K, \aa, \tt)$.

We can rewrite Theorem \ref{thm:O-N_traced_intro} in terms of cubic rings:
\begin{thm}[\textbf{O-N for traced cubic rings}]\label{thm:O-N_traced}
  Let $K$ be a number field, and let $\aa$, $\tt$ be ideals of $\OO_K$ with $\tt \mid 3$. Define $h_{\aa,\tt}(D)$ to be the number of $\tt$-traced cubic rings over $\OO_K$ with Steinitz class $\aa$ and discriminant ideal $D\aa^2$, each ring weighted by the reciprocal of its number of automorphisms:
  \[
  h_{\aa,\tt}(D)
  = \sum_{\substack{\Disc \OO = D\aa^2 \\ \tt\text{-traced}}} \frac{1}{\size{\Aut_K \OO}}.
  \]
  Then we have the global reflection theorem
  \begin{equation} \label{eq:O-N_traced}
    h_{\aa\tt^{-3}, 3 \tt^{-1}}(-27D)
    = \frac{3^{\#\{v|\infty : D \in (K_v^\cross)^2\}}}{N_{K/\QQ}(\tt)} \cdot h_{\aa,\tt}(D).
  \end{equation}
\end{thm}

We can also rewrite our results in terms of Shintani zeta functions.

\begin{defn}
Let $K$ be a number field. If $\OO/\OO_K$ is a cubic ring of nonzero discriminant, the \emph{signature} $\sigma(\OO)$ of $\OO$ is the Kummer element $\alpha \in (K\tensor_\QQ \RR)^\cross / ((K\tensor_\QQ \RR)^\cross)^2$ corresponding to the quadratic resolvent of $\OO$. That is, it takes the value $\alpha_v = +1$ or $-1$ at each real place $v$ of $K$ according as $\OO_v \isom \RR \cross \RR \cross \RR$ or $\RR \cross \CC$, and $\alpha_v = 1$ at each complex place.
\end{defn}

\begin{defn}\label{defn:Shintani}
Given a number field $K$, a signature $\sigma \in (K\tensor_\QQ \RR)^\cross / ((K\tensor_\QQ \RR)^\cross)^2$, an ideal class $[\aa] \in \Cl(K)$, and an ideal $\tt \mid 3$, we define the \emph{Shintani zeta function}
\[
  \xi_{K, \sigma, [\aa], \tt}(s) = \sum_{\OO} \frac{1}{\size{\Aut_K(\OO)}} N_{K/\QQ}(\Disc(\OO))^{-s}
\]
where the sum ranges over all cubic orders $\OO$ over $\OO_K$ having signature $\sigma$, Steinitz class $\aa$, and trace ideal contained in $\tt$. We also define the \emph{Shintani zeta function} with unrestricted Steinitz class
\[
  \xi_{K, \sigma, \tt}(s) = \sum_{[\aa] \in \Cl(K)} \xi_{K, \sigma, [\aa], \tt}(s).
\]
\end{defn}
\begin{rem}
We perpetuate the (confusing) tradition of denoting Shintani zeta functions by the Greek letter xi.
\end{rem}
\begin{rem}
By Minkowski's theorem on the finite count of number fields with bounded degree and discriminant, each term $n^{-s}$ has a finite coefficient, so the Shintani zeta function at least makes sense as a formal Dirichlet series. It generalizes the Shintani zeta functions for rings over $\ZZ$ mentioned in the introduction. Datskovsky and Wright \cite{DW2} study an adelic version of the Shintani zeta function; they show that $\xi_{K, \sigma, (1)}$ and $\xi_{K, \sigma, (3)}$ have meromorphic continuation with at most simple poles at $s = 1$ and $s = 5/6$, satisfying an explicit functional equation. We surmise that the same method will prove the same for $\xi_{K, \sigma, [\aa], \tt}$. However, we do not consider the analytic properties here.
\end{rem}

Then we have the following corollary, which generalizes Conjecture 1.1 of Dioses \cite{Dioses}.
\begin{cor}[\textbf{the extra functional equation for Shintani zeta functions}]\label{cor:Shintani}
Let $K$ be a number field, $\sigma \in (K\tensor_\QQ \RR)^\cross / ((K\tensor_\QQ \RR)^\cross)^2$ a signature, $[\aa] \in \Cl(K)$ an ideal class, and $\tt \mid 3$ an ideal. Then the Shintani zeta function $\xi_{K, \sigma, [\aa], \tt}(s)$ satisfies an extra functional equation
\begin{equation}\label{eq:Shintani_a_t}
  \xi_{K, \sigma, [\aa], \tt}(s)
  = 3^{3[K : \QQ]s - \#\{v|\infty : \sigma_v = 1\}}
    N_{K/\QQ}(\tt)^{1 - 6s}
    \xi_{K, -\sigma, [\aa \tt^{-3}], 3\tt^{-1}}(s),
\end{equation}
Hence, summing over all $\aa$,
\begin{equation}\label{eq:Shintani_t}
  \xi_{K, \sigma, \tt}(s)
  = 3^{3[K : \QQ]s - \#\{v|\infty : \sigma_v = 1\}}
    N_{K/\QQ}(\tt)^{1 - 6s}
    \xi_{K, -\sigma, 3\tt^{-1}}(s),
\end{equation}
\end{cor}
\begin{proof}
Fix $\aa$ and $\tt$. Sum Theorem \ref{thm:O-N_traced} over all $D \in \tt^3\aa^{-2}$ of signature $\sigma$, weighting each $D$ by 
\[
  N_{K/\QQ}\(D\aa^2\)^{-s},
\]
the norm of the discriminant of the associated cubic rings. Then the left-hand side of the summed equality matches the left-hand side of \eqref{eq:Shintani_a_t}. The right-hand side involves rings with discriminant ideal $27D\aa^2\tt^{-6}$, so a compensatory factor of
\[
  \frac{N_{K/\QQ}\(D\aa^2\)^{-s}}{N_{K/\QQ}\(27D\aa^2\tt^{-6}\)^{-s}} = \frac{3^{3[K:\QQ]s}}{N_{K/\QQ}(\tt)^{6s}}
\]
must be added to the right-hand side to pull out the desired Shintani zeta function.
\end{proof}

\subsubsection{The function field setting}
Lastly, we give an example of our results over the function field of a curve over a finite field, to illustrate that the same local reflection theorem can yield results in both settings.

\begin{thm}\label{thm:ON_fn_tame}
  Let $q \equiv 5 \mod 6$ be a prime power. Let $\A$ be a complete curve over $\FF_q$ with function field $K$, and $D$ a divisor on $\A$. Let $\B \to \A$ be a double cover (possibly singular) with branch divisor $D$, and let $L = K[\sqrt{d}]$ be its corresponding quadratic algebra over $K$. Let $L' = K[\sqrt{-3d}]$ be the $\FF_{q^2}$-twist of $L$. Then among the triple covers $\C \to \A$ with branch divisor $D$, the quadratic resolvent algebras $L$ and $L'$ occur equally often, when each cover $\C$ is weighted by the reciprocal of its number of automorphisms preserving $\A$.
\end{thm}
\begin{proof}
  Write $D = \sum_v m_v[v]$. At each $v$, choose a generator $e_v$ of the
  trace-zero summand of $\pi_*\OO_\B \tensor \OO_v$ and write
  $e_v = \alpha_v\sqrt d$. Since the discriminant of $(1,e_v)$ is
  $4d\alpha_v^2$ and $2$ is a unit, we have
  \[
  m_v = v(d) + 2v(\alpha_v).
  \]
  Thus $D$ and $\div(d)$ have the same parity, and we may define
  the adelic fractional ideal $\nn = \prod_v \nn_v$ by
  \[
  \div(\nn)
  = \frac12\bigl(D - \div(d)\bigr),
  \qquad \nn_v = \alpha_v\OO_v.
  \]
  Equivalently, $d\nn_v^2 = \pi_v^{m_v}\OO_v$; under the generic
  trivialization $\sqrt d \mapsto 1$, $\nn$ represents the trace-zero
  line bundle $\mathcal N$ of $\B$. Since $-3$ is a nonzero constant,
  $\div(-3d) = \div(d)$, so the same
  $\nn$ is attached to $D$ and $L' = K[\sqrt{-3d}]$.
  
  By the Deligne--Wood parametrization \cite[Theorem 2.1]{WQuartic}, a triple
  cover is given by a rank-$2$ bundle $E$ and a cubic map
  $E \to \Lambda^2 E$. Fixing the quadratic resolvent and choosing an
  orientation $\Lambda^2 E \isom \mathcal N$ identifies its generic cubic
  form with a point of the variety $X_d$ of cubic forms with discriminant $d$; the discriminant condition fixes
  the orientation up to sign. This variety has a composed variety structure $(\Gamma, V, X_d, X_{d,0})$ where $\Gamma = \SL_2$ acts on binary cubic forms $V = \Sym^3(2)$. Choose the adelic lattice
  \[
  \Lambda_{\nn} = \prod_v \Lambda_{v,\alpha_v,1}
  \]
  (the trace parameter $\tau = 1$ carries no information here). Then the possible oriented bundles $E$ are parametrized by
  \[
  \Gamma_{\Lambda_{\nn}}\bs\SL_2(\AA_K)/\SL_2(K),
  \]
  which is exactly the adelic class set occurring in Lemma
  \ref{lem:total_cl_num} with lattice $\Lambda_{\nn}$.
  Hence
  \[
  \frac{h(X_d)}{\size{H^0(K,M_d)}}
  \]
  is the weighted mass of oriented triple covers with branch divisor $D$ and
  quadratic resolvent algebra $L$.
  
  At every place $v$, Theorem \ref{thm:O-N_cubic_local}, with $\tau = 1$ and
  $\alpha = \alpha_v$, gives
  \[
  \widehat{\ell}_{v,\alpha_v,1,d}
  = \ell_{v,\alpha_v,3,-27d}
  = \ell_{v,\alpha_v,1,-27d},
  \]
  since $3 \in \OO_v^\cross$. Theorem \ref{thm:main_compose} therefore gives
  \[
  \frac{h(X_d)}{\size{H^0(K,M_d)}}
  = \frac{h(X_{-27d})}{\size{H^0(K,M_{-27d})}}.
  \]
  Since $K[\sqrt{-27d}] = K[\sqrt{-3d}] = L'$, the right-hand side is the
  corresponding weighted count of oriented covers with quadratic resolvent
  algebra $L'$, using the same adelic lattice $\Lambda_{\nn}$.
  
  Finally, the two admissible orientations of each triple cover differ by
  sign, and its automorphism group acts on this two-element set. Thus its
  total contribution to the oriented count is $2/\size{\Aut_\A(\C)}$.
  Dividing both sides by $2$ gives the claimed equality.
\end{proof}

\begin{rem}
  If $q \equiv 1$ mod $3$, then since $-3$ is a square, global reflection is a tautology if cast in this manner; however, function-field analogues of the discriminant-reducing formulas of Section \ref{sec:disc_red} are possible. In characteristic $2$, we instead have $L = K[\wp^{-1}(d)]$ and $L' = K[\wp^{-1}(d + 1)]$, where $\wp(x) = x^2 - x$ is the Artin-Schreier map. In characteristic $3$, things are quite different and will not be discussed in this paper.
\end{rem}
\begin{examp}
  On the curve $\A = \PP^1_{\FF_q}$, where $q \equiv 5$ mod $6$, consider the divisor
  \[
  D = 2(\infty) + 2(0).
  \]
  The following triple covers have branch divisor $D$:
  \begin{itemize}
    \item The union of three copies of $\A$, two of them glued together by simple nodes at $\infty$ and $0$. It has a twofold symmetry and quadratic resolvent $K \cross K$.
    \item The union of three copies of $\A$, one pair glued by a simple node at $\infty$ and a different pair glued by a simple node at $0$. It is asymmetric and has quadratic resolvent $K \cross K$.
    \item The union of $\A$ and the arithmetically irreducible, geometrically reducible double cover made by gluing together the two points at infinity on the conic $\mathcal{Q} = \{[X:Y:Z] \in \PP^2_{\FF_q} : X^2 + 3Y^2 = 0\}$, mapping to $\A$ via $[X:Y:Z] \mapsto [X:Z]$. This cover has a twofold symmetry $y \mapsto -y$ and has quadratic resolvent $K[\sqrt{-3}]$.
    \item $\PP^1$ itself, mapping to $\A$ by $t \mapsto t^3$. This is the only irreducible cover of the lot; it is asymmetric (there are no nontrivial cube roots of $1$ in $\FF_q$) and has quadratic resolvent $K[\sqrt{-3}]$.
  \end{itemize}
  Overall, the weighted number of covers of branch divisor $D$ is $3/2$ for each resolvent, in accord with the theorem.
\end{examp}

We now prove Theorem \ref{thm:O-N_cubic_local}, first in the \emph{tame} case where $\ch k_K \neq 3$, and then in the \emph{wild} case where $K$ is $3$-adic.

\subsection{Cubic rings}
Cubic and quartic rings have parametrizations, known as \emph{higher composition laws,} linking them to certain forms over $\OO_K$ and also to ideals in resolvent rings. The study of higher composition laws was inaugurated by Bhargava in his celebrated series of papers (\cite{B1,B2,B3,B4}), although the gist of the parametrization of cubic rings goes back to work of F.W.~Levi \cite{Levi} and was presented in greater generality by Delone and Faddeev \cite{DF}. Later work by Deligne and by Wood \cite{WQuartic, W2xnxn} has extended much of Bhargava's work from $\ZZ$ to an arbitrary base scheme. In a previous paper \cite{ORings}, the author explained how a representative sample of these higher composition laws extend to the case when the base ring $A$ is a Dedekind domain. In the present work, we will need a few more; fortunately, there are no added difficulties, and we will briefly run through the statements and the methods of proof.

\begin{thm}[\textbf{the parametrization of cubic rings}] \label{thm:hcl_cubic_ring} Let $A$ be a Dedekind domain with field of fractions $K$, $\ch A \neq 3$.
  \begin{enumerate}[$($a$)$]
    \item Cubic rings $\OO$ over $A$, up to isomorphism, are in bijection with cubic maps
    \[
      f : M \to \Lambda^2 M
    \]
    between a two-dimensional $A$-lattice $M$ and its own Steinitz class, up to isomorphism, in the obvious sense of a commutative square
    \[
    \begin{tikzcd}
      M_1 \arrow[r, "i", "\sim"'] \arrow[d, "f_1"'] & M_2 \arrow[d, "f_2"] \\
      \Lambda^2 M_1 \arrow[r, "\det i", "\sim"'] & \Lambda^2 M_2
    \end{tikzcd}
    \]
    The bijection sends a ring $\OO$ to the index form $f : \OO/A \to \Lambda^2(\OO/A)$ given by
    \[
    x \mapsto x \wedge x^2.
    \]
    \item \label{cubic:field} If $\OO$ is nondegenerate, that is, the corresponding cubic $K$-algebra $L = K \tensor_A \OO$ is \'etale, then the map $f$ is the restriction, under the Minkowski embedding, of the index form of $(K^\alg)^3$, which is
    \begin{equation}\label{eq:k3_phi}
      \begin{aligned}
        f : (K^\alg)^3 / K^\alg &\to (K^\alg)^2 / K^\alg \\
        (x;y;z) &\mapsto \bigl((x - y)(y - z)(z - x), 0\bigr).
      \end{aligned}
    \end{equation}
    \item \label{cubic:lift} Conversely, let $L$ be a cubic \'etale algebra over $K$. If $\bar{\OO} \subseteq L/K$ is a lattice such that
    $f$ sends $\bar{\OO}$ into $\Lambda^2 \bar{\OO}$, then there is a unique cubic ring $\OO \subseteq L$ such that, under the natural identifications, $\OO/A = \bar{\OO}$.
  \end{enumerate}
\end{thm}
\begin{proof}
  \begin{enumerate}[$($a$)$]
    \item 
    The proof is quite elementary, involving merely solving for the coefficients of the unknown multiplication table of $\OO$. The case where $A$ is a PID is due to Gross (\cite{cubquat}, Section 2): the cubic ring having index form
    \[
    f(x\xi + y\eta) = a x^3 + b x^2 y + c x y^2 + d y^3
    \]
    has multiplication table
    \begin{equation}\label{eq:mult_table}
      \begin{aligned}
        \xi\eta &= -a d \\
        \xi^2 &= -a c + b\xi - a\eta \\
        \eta^2 &= -b d + d\xi - c\eta.
      \end{aligned}
    \end{equation}
    For the general Dedekind case, see the author's \cite{ORings}, Theorem 7.1. It is also subsumed by Deligne's work over an arbitrary base scheme; see Wood \cite{WQuartic} and the references therein.
    \item Since the index form respects base change, we reduce to the case $K = K^\alg$. We then compute the index form of $(K^\alg)^3/K^\alg$ using the formula for the Vandermonde determinant.
    \item We have an integral cubic map $f\big|_{\bar \OO} : \bar \OO \to \Lambda^2 \bar \OO$, which is the index form $f_\OO$ of a unique cubic ring $\OO$ over $A$. But over $K$, $f_\OO$ is isomorphic to the index form of $L$. Since $L$ (as a cubic ring over $K$) is determined by its index form, we obtain an identification $\OO \tensor_A K \isom L$ for which the image of $\OO/A$ in $L/K$ coincides with $\bar{\OO}$. Under this identification, $\OO$ must be the inverse image of $\bar{\OO}$ under $L \to L/K$, which also proves uniqueness. \qedhere
  \end{enumerate}
\end{proof}
In this paper we only deal with nondegenerate rings, that is, those of nonzero discriminant, or equivalently, those that lie in an \'etale $K$-algebra. Consequently, all index forms $f$ that we will see are restrictions of \eqref{eq:k3_phi}. When cubic algebras are parametrized Kummer-theoretically, the resolvent map becomes very explicit and simple.

\begin{prop}[\textbf{self-balanced ideals in the cubic case}]
  \label{prop:hcl_cubic_sbi}
  Let $\OO_K$ be a Dedekind domain, $\ch K \neq 3$, and let $R$ be a quadratic \'etale $K$-algebra. A \emph{self-balanced triple} in $R$ is a triple $(B, I, \delta)$ consisting of a quadratic order $B \subseteq R$, a fractional ideal $I$ of $B$, and a scalar $\delta \in (KB)^\cross$ satisfying the conditions
  \begin{equation}
    \delta I^3 \subseteq B, \quad N(I) = (t) \text{ is principal}, \textand N(\delta) t^3 = 1,
  \end{equation}
  \begin{enumerate}[$($a$)$]
    \item \label{cubic:idl} 
    Fix $B$ and $\delta \in R^\cross$ with $N(\delta)$ a cube $t^{-3}$. Then the mapping
    \begin{equation}\label{eq:sbi_cubic}
      I \mapsto \OO = \OO_K + \kappa(I)
    \end{equation}
    defines a bijection between
    \begin{itemize}
      \item self-balanced triples of the form $(B, I, \delta)$, and
      \item subrings $\OO \subseteq L$ of the cubic algebra $L = K + \kappa(R)$ corresponding to the Kummer element $\delta$, such that $\OO$ is \emph{$3$-traced}, that is, $\tr(\xi) \in 3\OO_K$ for every $\xi \in \OO$.
    \end{itemize} 
    \item Under this bijection, we have the discriminant relation
    \begin{equation} \label{eq:cubic zmat disc}
      \disc \OO = -27 \disc B.
    \end{equation}
  \end{enumerate}
\end{prop}
\begin{proof}
  The mapping $\kappa$ defines a bijection between lattices $I \subseteq R$ and $\kappa(I) \subseteq L/K$. The difficult part is showing that $I$ fits into a self-balanced triple $(B, I, \delta)$ if and only if $\kappa(I)$ is the projection of a $3$-traced order $\OO$. Note that if $(B, I, \delta)$ exists, it is unique, as the requirement $[B : I] = (t)$ pins down $B$. 
  
  Rather than establish this equivalence directly, we will show that both conditions are equivalent to the symmetric trilinear form
  \begin{align*}
    \beta : I \cross I \cross I &\to \Lambda^2 R \\
    (\alpha_1, \alpha_2, \alpha_3) &\to \delta \alpha_1 \alpha_2 \alpha_3
  \end{align*}
  taking values in $t^{-1} \cdot \Lambda^2 I$.
  
  In the case of self-balanced ideals, this was done over $\ZZ$ by Bhargava \cite[Theorem 3]{B1}. Over a Dedekind domain, it follows from the parametrization of \emph{balanced} triples of ideals over $B$ \cite[Theorem 5.3]{ORings}, after specializing to the case that all three ideals are identified with one ideal $I$. It also follows from the corresponding results over an arbitrary base in Wood \cite[Theorem 1.4]{W2xnxn}.
  
  In the case of rings, we compute by Proposition \ref{prop:Kummer_resolvent_cubic} that $\beta$ is the trilinear form attached to the index form of $\kappa(I)$. By Theorem \ref{thm:hcl_cubic_ring}\ref{cubic:lift}, the diagonal restriction $\beta(\alpha, \alpha, \alpha)$ takes values in $t^{-1} \Lambda^2 (I)$ if and only if $\kappa(I)$ lifts to a ring $\OO$. We wish to prove that $\beta$ itself takes values in $t^{-1} \Lambda^2(I)$ if and only if $\OO$ is $3$-traced. Note that both conditions are local at the primes dividing $2$ and $3$, so we may assume that $\OO_K$ is a DVR. With respect to a basis $(\xi,\eta)$ of $I$ and a generator of $t^{-1} \Lambda^2(I)$, the index form of $\OO$ has the form
  \[
  f(x,y) = a x^3 + b x^2 y + c x y^2 + d y^3, \quad a,\ldots, d\in \OO_K.
  \]
  If this is the diagonal restriction of $\beta$, then $\beta$ itself can be represented as a $3$-dimensional matrix
  \[
  \bbq {a}{b/3}{b/3}{c/3}{b/3}{c/3}{c/3}{d\rlap{,}}
  \]
  which is integral exactly when $b, c \in 3\OO_K$. Since the trace ideal of $\OO$ is generated by
  \[
  \tr(1) = 3, \quad \tr(\xi) = b, \quad \tr(\eta) = -c
  \]
  (by reference to the multiplication table \eqref{eq:mult_table}), this is also the condition for $\OO$ to be $3$-traced, establishing the equivalence.
  
  The discriminant relation \eqref{eq:cubic zmat disc} follows easily from the definition of $\kappa$.
\end{proof}

\subsection{The tame case}
\label{sec:cubic_tame}

\begin{proof}[Proof of the tame case of Theorem \ref{thm:O-N_cubic_local}.]
Fix $D \in \OO_K$. Let $T = K[\sqrt{D}]$ be the corresponding quadratic algebra, and $T' = K[\sqrt{-3D}]$. For brevity we will write $H^i(T)$ for the cohomology $H^i(K, M_T)$ of the corresponding order-$3$ Galois module, and $H^i(T')$ likewise.

 Denote by $g(\sigma)$, for $\sigma \in H^1(T)$, the number of orders of discriminant $D$ in the corresponding cubic algebra $L_\sigma$; and likewise, denote by $g'(\tau)$, for $\tau \in H^1(T')$, the number of orders of discriminant $-3D$ in $L_\tau$. Our task is to prove that $g' = \widehat g$. We note that if $-3$ is a square in $\OO_K^\cross$, then $T = T'$ and $g = g'$.

Note that $g$ is even: $\sigma$ and $\sigma^{-1}$ are parametrized by the same cubic algebra with opposite orientations of its resolvent. The Fourier transform of an even, rational-valued function on a $3$-torsion group $H^1(T)$ is again even and rational-valued. So far, so good.

Our method will be first to prove the duality at $1$: that is, that
\begin{align}
  g'(1) &= \widehat g(1) \label{eq:dual_at_1 1} \\
  g(1) &= \widehat{g'}(1). \label{eq:dual_at_1 2}
\end{align}
Let us explain how this implies that $g' = \widehat g$. We compute $\size{H^1(T)}$ using the self-orthogonality of unramified cohomology:
\[
  \size{H^1(T)} = \size{H^1(T)^\ur} \cdot \size{H^1(T')^\ur}
  = \size{H^0(T)} \cdot \size{H^0(T')}.
\]
So there are basically three cases:
\begin{enumerate}[(a)]
  \item If neither $D$ nor $-3D$ is a square in $K_v$, then $H^1(T)$ and $H^1(T')$ are both trivial and \eqref{eq:dual_at_1 1} trivially implies that $g' = \widehat g$. 
  \item If one of $D$, $-3D$ is a square, then $H^1(T)$ and $H^1(T')$ are one-dimensional $\FF_3$-vector spaces. The space of even functions on each is $2$-dimensional, and
  \begin{align*}
    g &\mapsto (g(1),\widehat g(1)) \\
    g' &\mapsto (\widehat{g'}(1), g'(1))
  \end{align*}
  are corresponding systems of coordinates on them. Consequently, the two equations \eqref{eq:dual_at_1 1} and \eqref{eq:dual_at_1 2} together imply that $\widehat g = g'$.
  \item\label{it:H1=9} Finally, if $D$ and $-3D$ are both squares, then $g$ is a function on the two-dimensional $\FF_3$-space $H^1(T) \cong H^1(T')$ which we would like to prove self-dual. Note that $H^1(T)$ has four subspaces $L_1,\ldots, L_4$ of  dimension $1$. Consider the following basis for the five-dimensional space of even functions on $H^1(T)$: 
  \[
    g_1 = \1_{L_1}, \ldots, g_4 = \1_{L_4}, g_5 = \1_{\{1\}}.
  \]
  Note that $g_1,\ldots, g_4$ are self-dual (the Tate pairing is alternating, so any one-dimensional subspace is isotropic), while $g_5$ is not: indeed $\widehat g_5(1) \neq g_5(1)$. Thus if \eqref{eq:dual_at_1 1} holds, then $g$ is a linear combination of $g_1,\ldots,g_4$ only and hence $\widehat g = g$.
\end{enumerate}

We have now reduced the theorem to a pair of identities, \eqref{eq:dual_at_1 1} and \eqref{eq:dual_at_1 2}. By symmetry, it suffices to prove \eqref{eq:dual_at_1 2}, which may be written 
\begin{equation} \label{eq:dual_at_1}
  g(1) \stackrel{?}{=} \widehat{g'}(1) = \frac{1}{\size{H^0(T')}}\sum_{\tau \in H^1(T')} g'(\tau).
\end{equation}
The proof is clean and bijective.

The sum on the right-hand side of \eqref{eq:dual_at_1} counts all cubic orders of discriminant $-3D$. Any cubic order $C$ of discriminant $-3D$ can be assigned an ideal in $T$ as follows. Let $L$ be the fraction algebra of $C$. Since the Tate dual of the module associated to $T'$ is the module associated to $T$, Proposition \ref{prop:Kummer_cubic} gives
\[
  L = K + \{\xi \sqrt[3]{\delta} + \bar\xi \sqrt[3]{\bar\delta} \mid \xi \in T\}
\]
for some $\delta \in T^{N=1}/\bigl(T^{N=1}\bigr)^3$; and so, since $3$ is invertible in $\OO_K$,
\begin{equation} \label{eq:lattice C}
  C = \OO_K + \{\xi \sqrt[3]{\delta} + \bar\xi \sqrt[3]{\bar\delta} \mid \xi \in \cc\}
\end{equation}
for some lattice $\cc$ in $T$. Let $\OO_D \subseteq T$ denote the quadratic order of discriminant $D$ (equivalently $D/9$, since $3$ is a unit). By Proposition \ref{prop:hcl_cubic_sbi}, $(\OO_D,\cc,\delta)$ is a \emph{self-balanced ideal}, that is,
\begin{equation} \label{eq:balanced cubic}
  \delta \cc^3 \subseteq \OO_D, \quad N(\cc) = (t) \text{ is principal}, \textand N(\delta) t^3 = 1,
\end{equation}
Now $\cc$ need not be invertible in $\OO_D$; but let $\OO = \End \cc$. Then $\cc$ is an invertible $\OO$-ideal (here we use that $\OO$ is quadratic; we are essentially using that $\OO$ is Gorenstein over $\OO_K$). Note that $\OO = \OO_{D/\pi^{2i}}$ for some nonnegative integer $i$. Then form the \emph{shadow}
\[
  \aa = \frac{\delta \cc^3}{\pi^{i}}.
\]
Since the norm is multiplicative on invertible ideals, the properties of $\cc$ in \eqref{eq:balanced cubic} can be recast as properties of $\aa$:
\begin{equation} \label{eq:balanced aa}
  \aa \subseteq \End \aa \textand N_{\OO_D}(\aa) = 1.
\end{equation}
By sending
\[
  \aa \mapsto C_\aa = \OO_K + 0 \cross \aa \subseteq K \cross T,
\]
this $\aa$ parametrizes one of the rings counted by $g(1)$. (The first part of \eqref{eq:balanced aa} ensures that $C_\aa$ is closed under multiplication, and the second part ensures that $\disc C_\aa = D$. )

It remains to show that each shadow $\aa$ comes from exactly $\size{H^0(T')}$-many $C$, equivalently $(\cc, \delta)$, corresponding to each shadow $\aa$, up to multiplying $\delta$ by a cube and modifying $\cc$ by the appropriate equivalence relation
\begin{equation} \label{eq:equiv rel 2}
  \cc \mapsto \lambda \cc, \quad \delta \mapsto \lambda^{-3} \delta \quad (\lambda \in T^\cross).
\end{equation}
Let $\aa = \alpha\OO$ be the given shadow, where $\OO = \OO_{D/\pi^{2i}} = \End \aa$. Without loss of generality, we may scale $\alpha$ by an element of $\OO_K^\cross$ so that $N(\pi^i\alpha)$ is a cube. Clearly $\cc$ must be invertible with regard to $\OO$. The possible $(\cc,\delta)$ may be found by fixing $\cc = \OO$ and taking $\delta = \pi^i\alpha\epsilon$ where $\epsilon \in \OO^\cross$ has norm $1$. Now in the equivalence relation \eqref{eq:equiv rel 2}, the multipliers $\lambda$ preserving $\cc = \OO$ are $\lambda \in \OO^\cross$, so we must consider $\epsilon$ up to $\bigl(\OO^\cross\bigr)^3$. Since $N : \OO^\cross / \bigl(\OO^\cross\bigr)^3 \to \OO_K^\cross / \bigl(\OO_K^\cross\bigr)^3$ is surjective, the number of distinct $\epsilon$, which is the number of algebras $L_\delta$ in which we find rings, is
\[
  \frac{\Size{\OO^\cross / \bigl(\OO^\cross\bigr)^3}}{\Size{\OO_K^\cross / \bigl(\OO_K^\cross\bigr)^3}}.
\]
However, for each $\delta$, the contribution of $\aa$ to $g'(\delta)$ includes every $\cc \subseteq T$ such that $\delta\cc^3 = \delta\OO = \pi^{-i}\aa$. The number of such $\cc$ is
\[
  \frac{\Size{\OO_T^\cross[3]}}{\Size{\OO^\cross[3]}},
\]
and the total number of cubic orders we seek is the product
\begin{equation}\label{eq:ans 1}
  \frac{\Size{\OO^\cross / \bigl(\OO^\cross\bigr)^3}}{\Size{\OO_K^\cross / \bigl(\OO_K^\cross\bigr)^3}} \cdot
  \frac{\Size{\OO_T^\cross[3]}}{\Size{\OO^\cross[3]}}.
\end{equation}
To maneuver this into the required form, first note that
\[
  \frac{\Size{\OO^\cross / \bigl(\OO^\cross\bigr)^3}}{\Size{\OO^\cross[3]}} = \frac{\Size{\OO_T^\cross / \bigl(\OO_T^\cross\bigr)^3}}{\Size{\OO_T^\cross[3]}}
\]
by the Snake Lemma, since $\OO_T^\cross/\OO^\cross$ is finite; so \eqref{eq:ans 1} takes the form
\[
  \frac{\Size{\OO_T^\cross / \bigl(\OO_T^\cross\bigr)^3}}{\Size{\OO_K^\cross / \bigl(\OO_K^\cross\bigr)^3}},
\]
which we can now compute directly to equal
\[
  \Size{\OO_T^{N=1} / \bigl(\OO_T^{N=1}\bigr)^3} = \frac{\size{T^{N=1}/\bigl(T^{N=1}\bigr)^3}}{\size{H^0(T)}} = \frac{\size{H^1(T')}}{\size{H^0(T)}} = \size{H^0(T')}. \qedhere
\]
\end{proof}

\subsection{The wild case}\label{sec:wild_bij}
Recall that in Proposition \ref{prop:levels_new}, we constructed \emph{level spaces} $H^1(T) = \L_{-1} \supseteq \L_0 \supsetneq \L_1 \supsetneq \cdots \supsetneq \L_e \supseteq \L_{e+1} = \{1\}$ with the property that
\begin{equation}\label{eq:level space dual}
  \L_i^\perp = \L_{e-i}.
\end{equation}
Note that $\L_e \setminus \L_{e+1}$ is nonvanishing only when $T$ is split and then consists of the unramified field extension of $K$. Dually, $\L_{-1} \setminus \L_0$ is nonvanishing only when $T'$ is split, and then consists of the \emph{uniformizer radical extensions (UREs)} $K[\sqrt[3]{\pi}]$, the cubic extensions of maximal discriminant-valuation $3e + 2$. We have
\[
  \widehat{\1}_{\L_{e+1}} = \frac{1}{3} \1_{\L_{-1}} \textand \widehat{\1}_{\L_{-1}} = 3q^e \1_{\L_{e+1}}.
\]
We also let
\begin{align*}
  \ell_{\min} &= \begin{cases}
    -1 & T' \cong K \cross K \\
    0 & \text{otherwise},
  \end{cases} &
  \L_{\min} &= \L_{\ell_{\min}} = H^1(T), \\
  \ell_{\max} &= \begin{cases}
    e+1 & T \cong K \cross K \\
    e & \text{otherwise},
  \end{cases} &
  \L_{\max} &= \L_{\ell_{\max}} = \{1\}.
\end{align*}
We will now group all cubic rings whose resolvent torsor is $T$ into \emph{families} $\F$ with the following properties:
\begin{itemize}
  \item All rings in a family have the same discriminant and trace ideal.
  \item All rings in a family $\F$ are contained in \'etale algebras $L$ belonging to some level space $\L_i$; this $\L_i$ is called the \emph{support} $\supp(\F)$ of the family.
  \item Each $L \in \L_i$ has the same number of orders in the family; this number is called the \emph{thickness} $\th(\F)$ of the family.
\end{itemize}
The proof of Theorem \ref{thm:O-N_cubic_local} will then consist in exhibiting an involution between the families of support $\L_i$ and $\L_{e-i}$ which affects their thicknesses, discriminants, and trace ideals in such a manner that they contribute equally to both sides of the theorem. The remainder of this section will be spent in carrying this out.

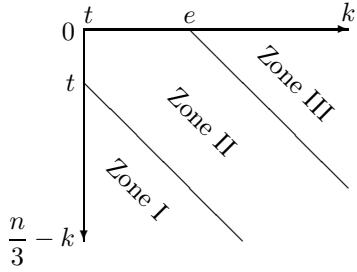
\begin{figure}
\setlength{\unitlength}{1em}
\[
\begin{picture}(13,9)(-2,-8.2)
  \put(0,0){\vector(1,0){10}}
  \put(0,0){\vector(0,-1){8}}
  \put(0,-2){\line(1,-1){6}}
  \put(4,0){\line(1,-1){6}}
  \put(-2,-9){\makebox(2,2)[r]{$\dfrac{n}{3} - k$~}}
  \put(-2,-3){\makebox(2,2)[r]{$t$~}}
  \put(-2,-1){\makebox(2,2)[r]{$0$\vphantom{$^2$}~}}
  \put( 0,0){\makebox(2,2)[bl]{$t$\vphantom{|}}}
  \put( 3,0){\makebox(2,2)[b]{$e$\vphantom{|}}}
  \put( 9,0){\makebox(2,2)[b]{$k$\vphantom{\big|}}}
  \put(2,-1){\rotatebox{-45}{\makebox(7,0){Zone II}}}
  \put(6,0){\rotatebox{-45}{\makebox(5.7,0){Zone III}}}
  \put(0,-4){\rotatebox{-45}{\makebox(5.7,0){Zone I}}}
\end{picture}
\]
\caption{Three zones used in the local counting of cubic forms}
\label{fig:zones_cubic}
\end{figure}

\begin{lem} \label{lem:cubic families}
Let $T/K$ be a quadratic \'etale algebra and let $n_T = v_K(\disc T) \in \{0,1\}$.

Then the cubic orders whose resolvent torsor is $T$ can be partitioned into families $\F_{n,k}$ indexed by the pairs of integers $(n,k)$ satisfying the conditions
\[
  0 \leq k \leq \floor{\frac{n}{3}}, \quad n \equiv n_T \mod 2,
\]
with the following properties:
\begin{enumerate}
  \item The rings in $\F_{n,k}$ have discriminant ideal $(\pi^n)$ and trace ideal $(\pi^{\min\{k,e\}})$.
  \item The support and thickness of each $\F_{n,k}$ depend on which of three zones the pair $(n,k)$ belongs to, as follows:
  \begin{equation}
    \begin{tabular}{cccc}
      Zone & $(n,k)$ & $\supp(\F_{n,k})$ & $\th(\F_{n,k})$ \\ \hline
      I & $0 \leq k < \dfrac{n}{6}$ & $\L_{\max}$ & $\size{H^0(T)}$ \\
      II & $\dfrac{n\vphantom{h}}{6} \leq k \leq \floor{\dfrac{n}{3}} - \dfrac{n}{6} + \dfrac{e}{2}$ & $\L_{e - 2k + \floor{\frac{n}{3}}}$ & $q^{\floor{\frac{n}{3}} - k}$ \\
      III & $\floor{\dfrac{n}{3}} - \dfrac{n\vphantom{h}}{6} + \dfrac{e}{2} < k \leq \floor{\dfrac{n}{3}}$ & $\L_{\min}$ & $q^{\floor{\frac{n}{3}} - k}$
    \end{tabular}
  \end{equation}
\end{enumerate}
\end{lem}
The shapes of these zones are like those in quadratic O-N, as shown in Figure \ref{fig:zones_cubic}.
\begin{rem}
Zone I, which is supported on $\L_{\max} = \{1\}$, consists precisely of those rings whose structure uses in an essential way that $L$ is \emph{split,} that is, has more than one field factor. Zone II has the approximate shape of a band of constant width,
\[
  \frac{n}{6} \leq k \lesssim \frac{n}{6} + \frac{e}{2};
\]
but the dependency on the value of $n$ modulo $3$ attests to a waviness of the boundary between zones II and III that cannot be avoided.
\end{rem}

Our first step is to compute the trace ideals $\tr(\OO_L)$ for the maximal orders of all cubic extensions $L$. If $L$ corresponds to the coclass $\sigma \in H^1(T)$, put
\[
  \ell(L) = \min\{e, \ell(\sigma)\}, \qquad k_0(L) = e - \ell(L),
\]
and define the \emph{discriminant offset} $h = h(L)$ by
\begin{equation}\label{eq:cubic_disc_offset}
  v_K(\disc L) = 3 k_0(L) + 2 - h(L).
\end{equation}
By Proposition \ref{prop:levels_new}\ref{lev:defect} and Lemma \ref{lem:disc}, one has $h \in \{1,2\}$ except for a uniformizer radical extension, for which $\ell(L) = -1$ and $h = 3$. For a split or unramified algebra, $\ell(L) = e$ and $h = 2$. For a ramified coclass with $0 \leq \ell(L) < e$, if $d = d(\sigma)$ is its defect, comparison with Lemma \ref{lem:disc} gives
\[
  d = 3\ell(L) + h(L).
\]
Consequently, \eqref{eq:dp1} gives
\[
  h(L) \equiv \theta(M_{T'}) - \ell(L) \pmod 2.
\]
Thus the discriminant offset $h(L)$ is determined by the module offset together with the level, but the two offsets are not the same invariant.

The trace ideal has a particularly simple description.
\begin{prop}\label{prop:tr_idl}
  Let $L/K$ be a cubic \'etale algebra over a $3$-adic field, with level $\ell = \ell(L)$ and discriminant offset $h = h(L)$. The trace ideal of the maximal order of $L$ is given by
  \[
    v(\tr(\OO_L)) = \min\{e, e - \ell(L)\}
  \]
  where $\ell(L) \in \{-1, 0, \ldots, e\}$ is the level. Moreover, if $L$ is totally ramified,
  \begin{itemize}
    \item If $h = 1$, $\tr(\OO_L)$ is generated by the trace of a uniformizer.
    \item If $h = 2$, $\tr(\OO_L)$ is generated by the trace of an element of valuation $2$.
  \end{itemize}
\end{prop}
\begin{proof}
  If $L$ is split or unramified, we have $\ell(L) = e$ and $\tr(\OO_L) = (1)$, so we may assume $L/K$ is a totally ramified field extension. Using \eqref{eq:cubic_disc_offset}, the result follows from the relation
  \[
    \min_{\xi \in L^\cross} v_L\(\frac{\tr \xi}{\xi} \) = v_K(\disc L) - [L : K] + 1
  \]
  due to Hyodo \cite[equation 1--4]{Hyodo}.
\end{proof}
\begin{proof}[Proof of Lemma \ref{lem:cubic families}]
We will now enumerate all the orders in every cubic $K$-algebra $L$ and arrange them into families. Let $L$ have level $e - k_0$, discriminant offset $h$, and discriminant ideal $(\pi^{n_0})$, so $\tr(\OO_L) = (3, \pi^{k_0})$ and hence, by \eqref{eq:cubic_disc_offset},
\[
n_0 = 3k_0 + 2 - h.
\]

Pick a basis $[1, \xi_0, \eta_0]$ for $\OO_L$ so that $\eta_0$ is traceless and $\tr(\xi_0)$ is a generator for the trace ideal $\tr(\OO_L)$. By Proposition \ref{prop:tr_idl}, if $L$ is totally ramified and $h \in \{1, 2\}$, we can choose the basis so that $v_L(\xi_0)=h$ and $v_L(\eta_0)=3-h$. Consequently,
\begin{equation}\label{eq:index_basis}
  v_K\bigl(f_{\OO_L}(\xi_0)\bigr) = h - 1
  \textand
  v_K\bigl(f_{\OO_L}(\eta_0)\bigr) = 2 - h,
\end{equation}
because a uniformizer generates the whole of $\OO_L$ while an element of $L$-valuation $2$ generates the unique subring of index $\mm_K$. In the cases where $L$ is not totally ramified, it is easy to arrange for \eqref{eq:index_basis} to hold as well.

Any order $C \subseteq L$ then has a unique basis of the form
\begin{equation} \label{eq:trace basis}
  [1, \quad \xi = \pi^i \xi_0 + u \eta_0, \quad \eta = \pi^j \eta_0]
\end{equation}
where $i$ and $j$ are nonnegative integers and $u$ ranges over a system of coset representatives in $\OO_K/\pi^j\OO_K$. Such a $C$ has discriminant valuation
\[
  n = v(\disc C) = n_0 + 2i + 2j
\]
and trace ideal
\[
  (3, \tr(\pi^i \xi_0)) = \pi^{\min\{e, k_0 + i\}}.
\]
Let $k = k_0 + i$. With one exception, namely when $L$ is a URE (which case we will handle later), we will place such a ring $C$ into the family $\F_{n,k}$. We must now compute the sizes of the families we have thus constructed.

Whether or not a lattice $C$ with a basis \eqref{eq:trace basis} is actually a ring is determined by the integrality of its index form. The following lemma reduces the number of coefficients to be checked from four to two.
\begin{lem}
Let $\xi, \eta \in \OO_L$ be integral elements of a nondegenerate cubic algebra $L$ over the field of fractions $K$ of a Dedekind domain $\OO_K$ such that the sublattice $C = \OO_K\<1, \xi, \eta\>$ is of full rank. If the two outer coefficients $f_C(\xi)$ and $f_C(\eta)$ of the index form of $C$ are integral, then the entire index form of $C$ is integral and $C$ is a ring.
\end{lem}
\begin{proof}
If the whole index form of $C$ is integral, then there is a ring $C' = \<1, \xi', \eta'\>$ with the same index form (with respect to its basis) as $C$. Using the identity of index forms for $L$ over $K$, we can embed $C'$ into $L$ with $\xi' = \xi + u$, $\eta' = \eta + v$ for some $u, v \in K$. But since $\xi', \xi, \eta', \eta$ are integral elements and $\OO_K$ is integrally closed, we have $u, v \in \OO_K$ so $C = C'$ is a ring.

So it suffices to prove that the index form is integral. This is a local statement, so we may assume that $\OO_K$ is a DVR. Passing to a finite extension, we may assume that $L \isom K \cross K \cross K$ is totally split. Let
\[
  \xi = (u_1; u_2; u_3) \textand \eta = (v_1; v_2; v_3).
\]
Then
\[
  \Delta = [\OO_L : C] = \det \begin{bmatrix}
    1 & 1 & 1 \\
    u_1 & u_2 & u_3 \\
    v_1 & v_2 & v_3
  \end{bmatrix}
\]
We are given that the outer coefficients
\[
  a = \frac{1}{\Delta}(u_1 - u_2)(u_2 - u_3)(u_3 - u_1) \textand
  d = \frac{1}{\Delta}(v_1 - v_2)(v_2 - v_3)(v_3 - v_1)
\]
of the index form of $C$ are integral. We wish to prove that the same applies to the two middle coefficients. By symmetry, we can consider just the $x^2y$-coefficient
\[
  b = \frac{1}{\Delta}\left[\sum_{i=1}^3(u_i - u_{i+1})(u_{i+1} - u_{i+2})(v_{i+2} - v_i)\right].
\]
(Here, and for the rest of the proof, indices are modulo $3$.) It is easy to verify that
\[
  b = -u_1 + 2u_2 - u_3 + \frac{3}{\Delta}(u_1 - u_2)(u_2 - u_3)(v_3 - v_1).
\]
So it is enough to show that
\[
  r_1 = \frac{1}{\Delta}(u_1 - u_2)(u_2 - u_3)(v_3 - v_1)
\]
is integral, or more generally any of the three
\[
  r_i = \frac{1}{\Delta}(u_i - u_{i+1})(u_{i+1} - u_{i+2})(v_{i+2} - v_i).
\]
A direct calculation also gives
\[
  r_1-r_2=u_3-u_2, \qquad r_2-r_3=u_1-u_3, \qquad r_3-r_1=u_2-u_1,
\]
so the differences $r_i-r_j$ are integral. Hence, if one $r_i$ had negative valuation, all three would have the same negative valuation. But
\[
  a^2 d = r_1 r_2 r_3,
\]
and the left-hand side is integral. Thus every $r_i$ is integral, completing the proof.
\end{proof}
\begin{rem}
The hypothesis that $L$ be nondegenerate is likely inessential to the truth of this lemma.
\end{rem}

We can now resume the proof of Lemma \ref{lem:cubic families}. Assume that $L$ is \emph{not} a URE, so $h \in \{1,2\}$. Let the index form of $\OO_L$ be
\[
  f_{\OO_L}(x\xi_0 + y\eta_0) = a x^3 + b x^2 y + c x y^2 + d y^3.
\]
Of the coefficients of the index form of $C$, we focus on the outer coefficients,
\[
  f_C (\xi)
  = \frac{a \pi^{3i} + b \pi^{2i} u + c \pi^{i} u^2 + d u^3}{\pi^{i+j}}
  \textand f_C (\eta) = d \pi^{2 j - i}.
\]
The latter coefficient is the simpler one, depending only on $i$ and $j$. Due to the tracelessness of $\eta_0$, we have $3|c$. By \eqref{eq:index_basis}, we have $v_K(d) = 2 - h$.
Thus the condition $f_C(\eta) \in \OO_K$ comes out to $2j - i + 2 - h \geq 0$, which simplifies to $n \geq 3k$.

(Incidentally, when $k \leq e$, the relation $n \geq 3k$ expresses an important relation between the trace ideal of a ring and its discriminant, generalizing the observation that an integer-matrix cubic form has discriminant divisible by $27$.)

There thus remains the $\xi$-condition $f_C (\xi) \in \OO_K$, which informally states that $\xi$ is a root of $f_{\OO_L}$ modulo $\pi^{i + j}$. Now $f_{\OO_L}$ is a homogeneous binary form, and it is natural to consider its roots on the projective lines $\PP^1(\OO_K/\pi^m)$; the factorization over the field $\OO_K/\mm$, for instance, gives the splitting type $\sigma(\OO_L)$. However, in our situation there is a distinguished point on this projective line, at least for $m < e - k_0$, namely the traceless point $\eta_0$: and the line is thereby subdivided into an affine line and a portion at infinity. The point at infinity mod $\mm$ is a root of $f_{\OO_L}$ if and only if $h < 2$. This is the motivation for the calculations to be undertaken now.

Suppose first that we are in zone I, that is, $n > 6k$, which translates into $j > 2i + 2 - h$. Then the $f_C(\xi)$ condition
\[
  \pi^{j - 2i} \mid a + b u \pi^{-i} + c u^2 \pi^{-2i} + d u^3 \pi^{-3i}
\]
is clearly dominated by a non-integral last term unless $u$ is of the form $\pi^i u'$, in which case it simplifies to
\[
  \pi^{j - 2i} \mid a + b u' + c u'^2 + d u'^3.
\]
In other words, $\xi_0 + u' \eta_0$ must be a root of $f_{\OO_L}$ modulo $\mm^{j - 2i}$. The condition $j - 2i > 2 - h$ rules out any contribution from a multiple root modulo $\mm$, which never lifts to mod $\mm^2$ (or else $\OO_L$ would be nonmaximal). So the only roots that contribute are the \emph{simple} roots occurring if $L$ has splitting type $111$, $12$, or $1^21$. There are $\size{H^0(T)}$ simple roots, and none of them are traceless (to be explicit, they are at $(1;0)$ for each decomposition $L \cong K \cross T$ into a linear and a quadratic algebra). By Hensel's lemma, each simple root has a unique lift to any modulus. So the solutions $u'$ form a union of $\size{H^0(T)}$ congruence classes modulo $\mm^{j - 2i}$. Since $u' = u/\pi^i$ is defined modulo $\mm^{j - i}$, there are
\[
  \size{H^0(T)} \cdot q^i = \size{H^0(T)} \cdot q^k
\]
rings for each pair $(i,j)$. This completes the construction of the families $\F_{n,k}$ in zone I.

Now suppose that $(n,k)$ is in zone II or III, still assuming that $L$ is \emph{not} a URE: we have
\[
  3k \leq n \leq 6k
\]
or, in the $(i,j)$ coordinates,
\begin{equation} \label{eq:zones II-III}
  i \leq 2j + 2 - h \textand j \leq 2i + \frac{n_0}{2} - 2 + h.
\end{equation}
We claim that there are rings for this pair $(i,j)$ if and only if $L \in \L_{e - 2k + \floor{\frac{n}{3}}}$, which may be also written as $k_0 \leq 2k - \frac{n - 2}{3}$ or in $(i,j)$ coordinates as
\begin{equation} \label{eq:zone supp}
  j \leq 2i - 1 + h.
\end{equation}
Assume first that $k_0 > 0$, that is, $L$ has splitting type $1^3$. Then the fact that $\xi_0$ is our chosen trace-ideal generator gives $v_L(\xi_0) = h$, and hence $v_K(a) = h - 1$. Now it is easy to show that
\[
  v(f_{\OO_L}(\xi)) = v(a\pi^{3i} + bu \pi^{2i} + cu^2 \pi^i + du^3) = \min\{3i + h - 1, 3 v(u) + 2 - h \}
\]
because the sum is dominated by its first term if $\pi^{i+h-1} | u$ and by its last term otherwise. So a necessary condition for there to be rings is that $i + j \leq 3i + h-1$, which is equivalent to \eqref{eq:zone supp}. If this condition holds, then the $\xi$-condition simply becomes
\begin{equation} \label{eq:vu}
  v(u) \geq \frac{i + j - 2 + h}{3}
\end{equation}
and we get
\[
  q^{j - \ceil{\frac{i+j-2+h}{3}}} = q^{\floor{\frac{2j + 2 - h - i}{3}}} = q^{\floor{\frac{n}{3} - k}}
\]
solutions.

If $L$ has one of the other splitting types, then our task is simplified by the facts that $k_0 = 0$ and $n_0 = 2 - h \in \{0,1\}$. Note that the second inequality of \eqref{eq:zones II-III} implies \eqref{eq:zone supp}, so we are only trying to prove that there \emph{are} solutions in this case. If \eqref{eq:vu} does \emph{not} hold, then the $d u^3$ term dominates in $f_{\OO_L}(\xi)$ and we do not get a solution. If \eqref{eq:vu} holds, we leave it to the reader to check the inequalities that imply (even without knowing anything about $a$, $b$, and $c$) that $\pi^{i+j}$ divides each term of $f_{\OO_L}(\xi)$. So we get the same number of solutions as in the preceding case.

Lastly, we must address the exceptional case that $h = 3$, that is, $L = K[\sqrt[3]{\pi}]$ is a URE. We take $\xi = \sqrt[3]{\pi}$ and $\eta = (\sqrt[3]{\pi})^2$, which are both traceless; and we have the explicit index form
\[
  f_{\OO_L}(x\xi + y\eta) = x^3 - \pi y^3.
\]
This resembles the index form for a ramified $L$ with $h = 1$, and analogously to that case, we compute that rings appear only for the pairs $(i,j)$ with
\begin{equation} \label{eq:URE zone}
  i \leq 2j + 1 \textand j \leq 2i,
\end{equation}
each such $(i,j)$ yielding $q^{\floor{\frac{2j + 1 - i}{3}}}$ solutions. We must place these rings into the families $\F_{n,k}$ of zone III so that the discriminant valuations and thicknesses match up. There is a unique choice:
\begin{align*}
  n &= n_0 + 2i + 2j = 3e + 2i + 2j + 2 \\
  k &= \floor{\frac{n}{3}} - \floor{\frac{2j + 1 - i}{3}}.
\end{align*}
We leave it to the reader that this establishes a bijection between the pairs $(i,j)$ in the region \eqref{eq:URE zone} with the pairs $(n,k)$ in zone III. In this way, each family in zone III is supported on the whole of $\L_{-1}$. The discriminant and thickness are correct by construction, completing the proof.
\end{proof}

A by-product of the foregoing is a formula for the number of traced orders in a given cubic algebra:

\begin{thm}[\textbf{the traced subring zeta function}]\label{thm:traced_count}
  Let $L$ be a cubic algebra over a $3$-adic field $K$ with discriminant $\disc L = \mm_K^{d_0}$. Let $d$ and $t$ be integers satisfying the necessary restrictions
  \[
  d \geq d_0, \quad d \equiv d_0 \mod 2, \quad 0 \leq t \leq e_K.
  \]
  Then the number $g(L, d, t)$ of $\mm_K^t$-traced orders of discriminant $\mm_K^d$ in $L$ is a linear combination of the three functions
  \begin{align*}
    g^{1^3}(d_0, d, t) &= \frac{q^r - 1}{q - 1}, \quad r = \begin{cases}
      0, & d < 3t \\
      \floor{\frac{d}{3}} - t + 1, & 3t \leq d \leq 6t - d_0 \\
      \floor{\frac{d - d_0}{6}} + 1, & d \geq 6t - d_0
    \end{cases} \\
    g^{3}(d_0 = 0, d, t) &= \frac{q^r - 1}{q - 1}, \quad r = \begin{cases}
      0, & d < 3t \\
      \floor{\frac{d}{3}} - t + 1, & 3t \leq d \leq 6t \\
      \frac{d}{2} - 2 \ceil{\frac{d}{6}} + 1, & d \geq 6t
    \end{cases} \\
    g^{\SR}(d_0 = 0 \text{ or } 1, d, t) &= \frac{q^r - q^t}{q - 1}, \quad r = \begin{cases}
      t, & d < 6t \\
      \floor{\frac{d}{6}} - t + 1, & d \geq 6t
    \end{cases}
  \end{align*}
  in a manner dependent on the splitting type of $L$:
  \begin{itemize}
    \item $\sigma(L) = 1^3$: $g = g^{1^3}$
    \item $\sigma(L) = 1^21$: $g = g^{1^3} + g^{\SR}$
    \item $\sigma(L) = 3$: $g = g^3$
    \item $\sigma(L) = 12$: $g = g^3 + g^{\SR}$
    \item $\sigma(L) = 111$: $g = g^3 + 3g^{\SR}.$
  \end{itemize}
\end{thm}
\begin{rem}
  When $t = 0$ (and here we need no longer assume that $\ch k_K = 3$), we recover the formulas for orders in a cubic field computed by Datskovsky and Wright \cite{DW2} and written more explicitly by Nakagawa \cite{Nakagawa} and the author \cite{OOnARemarkable}.
\end{rem}

\begin{proof}[Proof of Theorem \ref{thm:O-N_cubic_local}]
  Once Lemma \ref{lem:cubic families} is proved, we can prove Theorem \ref{thm:O-N_cubic_local} quite simply by sending the family $\F = \F_{n,k}$ to $\F' = \F_{n',k'}$, where
  \begin{align*}
    n' = n + 3e - 6t \\
    k' = \floor{\frac{n}{3}} - k + e - t.
  \end{align*}
  If the original $\F$ satisfied the bounds $t \leq k \leq \floor{\frac{n}{3}}$, then it is easy to see that $t' \leq k' \leq \floor{\frac{n}{3}}$ where $t' = e - t$, and likewise $n \equiv n_T$ mod $2$ implies $n' \equiv n_{T'}$. Thus $\F'$ is a family of rings of resolvent torsor $T'$ whose trace ideal is contained in $(\pi^{e-t})$. It is not hard to see that $\F'$ lies in zone III, II, or I according as $\F$ lies in zone I, II, or III. We leave it to the reader to check the needed identities
  \begin{align*}
    \supp(\F') &= \supp(\F)^\perp \\
    \th(\F') &= \frac{\size{\supp(\F)}}{\size{H^0(T)}} \cdot \th(\F).
  \end{align*}
\end{proof}

\begin{rem}
The method of the above proof can also be adapted to the tame case.
\end{rem}

\section{More general weightings; discriminant reduction}
\label{sec:non-natural}

While the natural weighting $w_v = 1$ yields elegant reflection theorems, the power of this technique can be amplified by selecting a different weighting at one prime $p$, or a small selection of primes, keeping the known local reflection at all others. This has applications to counting cubic rings satisfying local conditions, and it sometimes has the benefit of \emph{discriminant reduction} whereby the cubic forms on one side have smaller discriminant and are therefore easier to count.

For simplicity we work over $\QQ$, though the techniques extend. We take $(\Gamma, V, X_D, X_{D,0})$ to be the composed variety of binary cubic forms of discriminant $D$ as before. Its stabilizer $M_D$ is the group $C_3$ with Galois action given by the quadratic character corresponding to $\QQ(\sqrt{D})$.

As in Section \ref{sec:composed}, if $w : V(\OO_{\AA_\QQ}) \to \CC$ is a locally constant weighting invariant under $\GL_2(\OO_{\AA_\QQ})$, we denote by $h(D, w)$ the number of $\GL_2 \ZZ$-classes of binary cubic forms over $\ZZ$ of discriminant $D$, each form $f$ weighted by
\[
\frac{w(f)}{\lvert \Stab_{\GL_2\ZZ} f \rvert}.
\]
If $w = \prod_v w_v$ is a product of local weightings, then our local-to-global reflection engine (Theorem \ref{thm:main_compose}) produces identities relating different $h(D, w)$, if we can find a dual for each $w_v$.

\subsection{Local weightings given by splitting types}

Let $\sigma \in \{111,12,3,1^21, 1^3, 0\}$ be one of the six \emph{splitting types} a binary cubic form can have at a prime. Let
\[
  T(\sigma) = T_p(\sigma) : V(\ZZ_p) \to \{0, 1\}
\]
be the selector that takes the value $1$ on binary cubic forms of splitting type $\sigma$. Then the associated weighted local lattice-count function
\[
  \ell_{p,D}^{T_p(\sigma)} : H^1\(\QQ_p, M_D\) \to \NN
\]
attaches to each cubic algebra $L/\QQ_p$ whose quadratic resolvent is $\QQ_p(\sqrt{D})$ its number of orders of discriminant $D$ and splitting type $\sigma$.

There is another construction of interest to us. Let $\Lambda_\ZZ = V(\ZZ)$ be the standard lattice. For $a \in \QQ^\cross$, the integral slices
\begin{equation}\label{eq:shifted_slices}
  X_{D_0}(\QQ) \intsec \Lambda_\ZZ
  \textand
  X_{a^2D_0}(\QQ) \intsec \Lambda_\ZZ
\end{equation}
need not look alike, although the rational varieties $X_{D_0}$ and $X_{a^2D_0}$ are isomorphic. Indeed, under the twisted action \eqref{eq:twisted_action}, any $g \in \GL_2(\QQ)$ with $\det g = a$ carries $X_{D_0}$ to $X_{a^2D_0}$. Pulling the standard lattice back through such an isomorphism therefore lets us regard the two sets in \eqref{eq:shifted_slices} as arising from two lattice choices on the same rational composed variety. Theorem \ref{thm:main_compose} allows this change to be made at one place while the other local lattices are left fixed.

For a prime $p$, put
\[
  z_p = \begin{bmatrix}
    1/p & \\
        & 1
  \end{bmatrix},
\]
so that the twisted action of $z_p$ carries $X_D$ to $X_{Dp^{-2}}$. Let $\Lambda_p = V(\ZZ_p)$ be the standard local lattice. If $T$ is a local weighting on $\Lambda_p$ and $k \in \ZZ$, define the \emph{$k$-shift} of the weighted lattice datum $(\Lambda_p,T)$ by
\[
  \Lambda_p^{[k]} = z_p^{-k}\Lambda_p
  \textand
  T^{[k]}(x) = T(z_p^k x), \qquad x \in X_D(\QQ_p) \intsec \Lambda_p^{[k]}.
\]
Thus $T$ is always evaluated on the standard lattice, even when $T^{[k]}$ is regarded as a weighting on the shifted lattice. Transport by $z_p^k$ gives
\begin{equation}\label{eq:shifted_weighting}
  \ell_{p,D}^{T^{[k]}} = \ell_{p,Dp^{-2k}}^T.
\end{equation}
We use the same notation globally when the local datum at $p$ is shifted and all other local data are left fixed; in particular,
\[
  h(D,T^{[k]}) = h(Dp^{-2k},T).
\]
We also allow finite formal linear combinations of weighted lattice data and extend $\ell$ and $h$ linearly.
Products of local data at distinct primes have their evident meaning. Thus, for example,
\[
  R_p^{[1]} + (1-R_p)^{[2]}
\]
is understood as a formal sum of two weighted lattice data on $X_D$. Lastly, if the specification $T$ does not involve the prime $3$, let
\[
h_3(D, T) = h\(D, T \cdot \(T_3(0) + T_3(1^3)\)\)
\]
count only the $3$-traced forms.

For a general nonarchimedean local field $K$ with uniformizer $\pi$, we define $T^{[k]}$ in the same way using $z_\pi = \operatorname{diag}(\pi^{-1},1)$.

\begin{lem} \label{lem:disc red} Let $K$ be a local field, $\ch k_K \neq 3$, and let $D \in \OO_K^\cross$. On the standard lattices, the weighting $T(1^3)$ on $X_{-3\pi^2D}$ and the weighting $2 \cdot T(111) - T(3)$ on $X_D$ are dual with duality constant $1$; that is,
\begin{equation} \label{eq:disc red}
  \widehat{\ell}_{v,-3 \pi^2 D}^{T(1^3)} = 2 \ell_{v,D}^{T(111)} - \ell_{v,D}^{T(3)}.
\end{equation}
\end{lem}
\begin{proof}
The right-hand side of \eqref{eq:disc red} can be written as
\[
  \size{H^0(M_D)} \cdot \1_{\{1\}} - \1_{H^1_\ur},
\]
so it suffices to show that the left-hand side is the Fourier transform of this, namely
\[
  1 - \1_{H^1_\ur} = \1_{H^1_{\ram}}.
\]
Look at binary cubic forms $f(x,y)$ of splitting type $1^3$ and discriminant $-3 \pi^2 D$. Changing coordinates, we can assume
\[
  f(x,y) = a x^3 + b x^2 y + c x y^2 + d y^3 \equiv x^3 \mod \pi.
\]
Then note that $\disc f \equiv - 27a^2d^2 \equiv -27d^2$ mod $\pi^3$, so the only way that $f$ can have discriminant $-3 \pi^2 D$ is if $\pi^2 \nmid d$. Then $f$ is an Eisenstein polynomial, the index form of a maximal order in a totally ramified extension $L$. Hence the weighting counting such $f$ is $\1_{H^1_{\ram}}$, as desired.
\end{proof}

Plugging this, together with the natural duality of Theorem \ref{thm:O-N_cubic_local} at the other primes, into Theorem \ref{thm:main_compose} yields results such as the following:
\begin{thm}\label{thm:1^3}
    Let $p \in \ZZ$ be a prime, $p \neq 3$. For all integers $D$ such that $p\nmid D$,
    \begin{align}
    \frac{1}{c_\infty} h_3(-27 p^2 D, T_p(1^3)) &= 2 h(D, T_p(111)) - h(D, T_p(3)) \label{eq:1^3 1} \\
    c_\infty h(p^2 D, T_p(1^3)) &= 2 h_3(-27 D, T_p(111)) - h_3(-27 D, T_p(3)) \label{eq:1^3 2}
    \end{align}
    where $c_\infty = 3$ for $D > 0$, $c_\infty = 1$ for $D < 0$.
\end{thm}

\subsection{Discriminant reduction}
\label{sec:disc_red}
This can be used to improve a step that often occurs in arithmetic statistics, namely the production of \emph{discriminant-reducing} identities that express the number of forms with certain non-squarefree discriminant recursively in terms of lower discriminants.

For $N$ a positive integer, let $ R_N $ be the local weighting that assigns to a binary cubic form its number of roots in $\PP^1(\ZZ/N\ZZ)$. This is a local weighting at the primes dividing $N$; for instance, if $N = p$ is prime, then
\[
  R_p = 3 T_p(111) + T_p(12) + 2 T_p(1^21) + T_p(1^3) + (p + 1) T_p(0).
\]
The same weighting may be constructed by taking the lattice $\Lambda_N$ of binary cubic $111N$ forms, whose stabilizer $\Gamma_{\Lambda_N}$ is a certain congruence subgroup. This enables us to state succinctly the following theorem.
\begin{thm}[\textbf{discriminant reduction}]\label{thm:disc_red}
Let $p \nmid 6$ be a prime, and $D$ an integer divisible by $p^2$. Then
\[
  h(D) = h\(\frac{D}{p^2}, R_p\) + h\(\frac{D}{p^4}\) - h\(\frac{D}{p^4}, R_p\)
    + \frac{1}{c_\infty}\(2 h_3\(\frac{-27 D}{p^2}, T_p(111)\) - h_3\(\frac{-27 D}{p^2}, T_p(3)\)\),
\]
where $c_\infty = 3$ for $D > 0$, $c_\infty = 1$ for $D < 0$.
\end{thm}
\begin{proof}
The cubic rings $C$ counted by the left-hand side can be divided into maximal and nonmaximal at $p$. If $C$ is maximal at $p$, then $C \tensor_\ZZ \ZZ_p$ is the ring of integers of a totally tamely ramified cubic extension of $\QQ_p$ and $p^2 \parallel D$. By Theorem \ref{thm:1^3}, such rings are counted by the last term.

If $C$ is nonmaximal at $p$, then $C$ sits with $p$-power index inside an overring $C'$. By considering $C' + p C$, we can take an inclusion $C \subset C'$ of one of the following forms:
\begin{itemize}
  \item $C$ has index $p$ in a $C' = C_1$ of discriminant $D/p^2$. Here the index form of $C_1$ must have a marked root modulo $p$ so that the transformation
  \[
    f_{C_1}(x,y) = a x^3 + b x^2 y + c x y^2 + d y^3
    \longmapsto f_C(x,y) = p^2 a x^3 + p b x^2 y + c x y^2 + \frac{d}{p} y^3
  \]
  keeps the form integral. This accounts for the term $h(D/p^2, R_p)$.
  \item $C$ has index $p^2$ in a $C' = C_2$ of discriminant $D/p^4$ with $C_2/C \isom (\ZZ/p\ZZ)^2$. This requires that the index form of $C$ have \emph{content} divisible by $p$; we have
  \[
    f_C = p f_{C_2}.
  \]
  This accounts for the term $h(D/p^4)$.
\end{itemize}
Observe that $C_2$ is unique if it exists. A choice of $C_1$ corresponds to a choice of multiple root of $f_C$, which is unique if $f_C$ is nonzero modulo $p$. Thus, the only chance of overcounting occurs when a $C$ admits both a $C_2$ and one or more $C_1$'s. The $C_1$'s are all the subrings of index $p$ in $C_2$ and thus correspond to the roots of $f_{C_2}$ modulo $p$. So we subtract $1$ (more precisely, $1/\size{\Aut C_2}$) for each root of a form $f_{C_2}$ counted in the term $h(D/p^4)$. That is, we subtract $h(D/p^4, R_p)$, yielding the claimed total.
\end{proof}

More generally, we can reduce at multiple primes at once. Let $T_p^{\max} : V(\OO_{\AA_\QQ}) \to \ZZ$ be the selector for rings maximal at $p$.

\begin{thm}[\textbf{discriminant reduction}]\label{thm:disc_red_gen}
Let $q = p_1\cdots p_r$ be a squarefree integer, $\gcd(q, 6) = 1$. If $D$ is a nonzero integer divisible by $q^2$, then for any $t < q$,
\begin{align*}
  h(D) &= \sum_{\substack{q = q_1q_2q_3 \\ q_1 \leq t}} h\(\frac{D}{q_2^2 q_3^4}, \prod_{p|q_1} T_p^{\max} \prod_{p|q_2} R_p \prod_{p|q_3} (1 - R_p)\) + {} \\
  & \quad + \frac{1}{c_\infty} \sum_{\substack{q = q_1q_2q_3 \\ q_1 > t}} h_3 \(\frac{-27 D}{q_1^2 q_3^2}, \prod_{p|q_1} \1_{p^2\parallel D}(R_p - 1) \prod_{p|q_3} \1_{p^2 \parallel D}(1 - R_p) \) \\
  \intertext{and}
  h_3(-27 D) &= \sum_{\substack{q = q_1q_2q_3 \\ q_1 \leq t}} h_3\(\frac{-27 D}{q_2^2 q_3^4}, \prod_{p|q_1} T_p^{\max} \prod_{p|q_2} R_p \prod_{p|q_3} (1 - R_p)\) + {} \\
  & \quad + {c_\infty} \sum_{\substack{q = q_1q_2q_3 \\ q_1 > t}} h_3 \(\frac{D}{q_1^2 q_3^2}, \prod_{p|q_1} \1_{p^2\parallel D}(R_p - 1) \prod_{p|q_3} \1_{p^2 \parallel D}(1 - R_p) \).
\end{align*}
\end{thm}
\begin{rem}
If we take $t = \sqrt{q}$, we find that all discriminants appearing have absolute value at most $27|D|/q$.
\end{rem}
\begin{proof}
Since a ring of discriminant $D$ is maximal or nonmaximal at each of the primes dividing $q$, we have
\[
  h(D) = \sum_{q_1 q_2' = q} h\(D, \prod_{p|q_1} T_p^{\max}(1^3) \prod_{p|q_2'} T_p^{\nonmax}\),
\]
where $T_p^{\max}(1^3) = T_p^{\max} \cdot T_p(1^3)$ and $T_p^{\nonmax} = 1 - T_p^{\max}$ (as is natural).
We transform each term in one of two ways, depending on whether $q_1 \leq t$.

If $q_1 \leq t$, we replace each $T_p^{\nonmax}$ by the formal sum
\[
  R_p^{[1]} + (1 - R_p)^{[2]}
\]
by the method of the preceding theorem. Since the shift is local at $p$, this respects local conditions at the other primes. We get
\begin{align*}
  & h\(D, \prod_{p|q_1} T_p^{\max}(1^3) \prod_{p|q_2'} T_p^{\nonmax}\) \\
  &= h\(D, \prod_{p|q_1} T_p^{\max}(1^3) \prod_{p|q_2'} \bigl(R_p^{[1]} + (1 - R_p)^{[2]}\bigr)\) \\
  &= \sum_{q_2' = q_2q_3}  h\(\frac{D}{q_2^2 q_3^4}, \prod_{p|q_1} T_p^{\max} \prod_{p|q_2} R_p \prod_{p|q_3} (1 - R_p)\),
\end{align*}
where the last equality is \eqref{eq:shifted_weighting} applied at each shifted prime.

If $q_1 > t$, we reflect. By Lemma \ref{lem:disc red}, a dual of $T_p^{\max}$, when restricted to discriminants $D$ divisible by $p^2$, is
\[
  \1_{p^2 \parallel D}(R_p - 1)^{[1]}.
\]
Hence a dual of $T_p^{\nonmax} = 1 - T_p^{\max}$ on the same discriminants is the formal sum
\[
  1 + \1_{p^2 \parallel D}(1 - R_p)^{[1]}.
\]
Applying the reflection theorem termwise and then using \eqref{eq:shifted_weighting},
\begin{align*}
  & h\(D, \prod_{p|q_1} T_p^{\max}(1^3) \prod_{p|q_2'} T_p^{\nonmax}\) \\
  &= \frac{1}{c_\infty} h_3\(-27D, \prod_{p|q_1} \1_{p^2\parallel D}(R_p - 1)^{[1]} \cdot \prod_{p|q_2'} \bigl(1 + \1_{p^2 \parallel D}(1 - R_p)^{[1]}\bigr)\) \\
  &= \frac{1}{c_\infty} \sum_{\substack{q_2' = q_2q_3}} h_3 \(\frac{-27 D}{q_1^2 q_3^2}, \prod_{p|q_1} \1_{p^2\parallel D}(R_p - 1) \prod_{p|q_3} \1_{p^2 \parallel D}(1 - R_p) \).
\end{align*}
Summing over $q_1$ yields the first identity. The second is proved in the same way.
\end{proof}

\subsection{Binary cubic forms over \texorpdfstring{$\ZZ[1/N]$}{Z[1/N]}}
\label{sec:Z[1/N]}

The local-to-global reflection engine has a useful variant in which integrality is dropped at one set of places and, at a disjoint set, the local orbit is required to be the integral base orbit. We record it here in the form needed below.

Let $(\Gamma, V, X, X_0)$ be a composed variety over a global field $K$ with adelic lattice $\Lambda = \prod_v \Lambda_v$. If $\Sigma$ is a finite set of nonarchimedean places, set
\[
  \Lambda^\Sigma
  = \prod_{v \in \Sigma} V(K_v) \times \prod_{v \notin \Sigma} \Lambda_v
\]
and
\[
  \Gamma^\Sigma
  = \Gamma(K) \intsec
    \left(\prod_{v \in \Sigma} \Gamma(K_v)
      \times \prod_{v \notin \Sigma} \Gamma_{\Lambda_v}\right).
\]
Thus $\Gamma^\Sigma$ is the group of rational elements preserving the $\Sigma$-integral subgroup $\Lambda^\Sigma$. If $\Theta$ is a second finite set of nonarchimedean places disjoint from $\Sigma$, choose, for each $v \in \Theta$, an integral point $x_{0,v} \in X_0(K_v) \intsec \Lambda_v$, and let
\[
  X_\Theta^\Sigma(K)
  = \left\{x \in X(K) :
      x_v \in \Lambda_v \text{ for } v \notin \Sigma,
      \quad x_v \in \Gamma_{\Lambda_v}x_{0,v} \text{ for } v \in \Theta
    \right\}.
\]
If $w = \prod_{v \notin \Sigma \cup \Theta} w_v$ is a product of local weightings away from the exceptional places, extended by the trivial weighting there, define
\begin{equation}\label{eq:S_class_number}
  h_\Theta^\Sigma(X,w)
  = \sum_{x \in \Gamma^\Sigma \bs X_\Theta^\Sigma(K)}
      \frac{w(x)}{\size{\Stab_{\Gamma^\Sigma}x}}.
\end{equation}

\begin{thm}[\textbf{$S$-integral reflection engine}]\label{thm:S_reflection}
  For $i = 1,2$, let
  $\bigl(\Gamma^{(i)}, V^{(i)}, X^{(i)}, X_0^{(i)}\bigr)$ be Hasse composed varieties over $K$ whose stabilizers $M^{(1)}$ and $M^{(2)}$ are Tate duals, and put
  \[
    m = \size{M^{(1)}} = \size{M^{(2)}}.
  \]
  Let $\Sigma^{(1)}$ and $\Sigma^{(2)}$ be disjoint finite sets of nonarchimedean places. Equip the two varieties with adelic lattices $\Lambda^{(i)}$ and suppose that both varieties are full over $K_v$ for every $v \in \Sigma^{(1)} \cup \Sigma^{(2)}$. Choose integral basepoints $x_{0,v}^{(i)} \in X_0^{(i)}(K_v) \intsec \Lambda_v^{(i)}$ at the exceptional places.

  Assume that, for each $i$, the $\Sigma^{(i)}$-adelic class set
  \[
    \left(\prod_{v \in \Sigma^{(i)}} \Gamma^{(i)}(K_v)
      \times \prod_{v \notin \Sigma^{(i)}} \Gamma_{\Lambda_v^{(i)}}^{(i)}\right)
    \bs \Gamma^{(i)}(\AA_K) / \Gamma^{(i)}(K)
  \]
  has one element. Away from $\Sigma^{(1)} \cup \Sigma^{(2)}$, let the local weightings $w_v^{(1)}$, $w_v^{(2)}$ be dual with constant $c_v$ as in Definition \ref{defn:dual}. Then
  \begin{align}\label{eq:S_reflection}
    h_{\Sigma^{(1)}}^{\Sigma^{(2)}}(X^{(2)},w^{(2)})
    &= \left(\prod_{v \notin \Sigma^{(1)} \cup \Sigma^{(2)}} c_v^{-1}\right)
       \left(\prod_{v \in \Sigma^{(1)}}
          \frac{1}{
            \size{\Stab_{\Gamma_{\Lambda_v^{(2)}}^{(2)}}x_{0,v}^{(2)}}
            \size{\OO_v/m\OO_v}}
       \right) \notag\\
    &\qquad \cdot
       \left(\prod_{v \in \Sigma^{(2)}}
          \size{\Stab_{\Gamma_{\Lambda_v^{(1)}}^{(1)}}x_{0,v}^{(1)}}
       \right)
       h_{\Sigma^{(2)}}^{\Sigma^{(1)}}(X^{(1)},w^{(1)}).
  \end{align}
\end{thm}
\begin{proof}
  The proof is the same Poisson-summation argument as Theorem \ref{thm:main_compose}; we only record the local factors at the exceptional places. The class-number-one hypotheses and the same regrouping as in Lemma \ref{lem:total_cl_num} identify the two sides of \eqref{eq:S_class_number} with the corresponding normalized sums over global coclasses.

  At a place where denominators are allowed, fullness makes the local cohomological counting function identically $1$. At a place where we restrict to the integral base orbit, the local counting function is
  \begin{equation}\label{eq:base_orbit_local_count}
    \frac{\size{H^0(K_v,M)}}
         {\size{\Stab_{\Gamma_{\Lambda_v}}x_{0,v}}}
    \1_{\{1\}}.
  \end{equation}
  Indeed, the elements of $\Gamma(K_v)$ carrying $x_{0,v}$ into its $\Gamma_{\Lambda_v}$-orbit form $\Gamma_{\Lambda_v}\Stab_{\Gamma(K_v)}x_{0,v}$, and Proposition \ref{prop:composed}\ref{comp:stab} identifies the latter stabilizer with $H^0(K_v,M)$.

  For $v \in \Sigma^{(1)}$, the local function for the first variety is therefore $1$, while the desired function for the second variety is \eqref{eq:base_orbit_local_count}. By the local Euler characteristic formula and Tate duality,
  \begin{equation}\label{eq:local_Euler_size}
    \size{H^1(K_v,M^{(1)})}
    = \size{H^0(K_v,M^{(1)})}
      \size{H^0(K_v,M^{(2)})}
      \size{\OO_v/m\OO_v}.
  \end{equation}
  Hence
  \[
    \widehat{1}
    = \size{H^0(K_v,M^{(2)})}
      \size{\OO_v/m\OO_v}\,\1_{\{1\}},
  \]
  and the exceptional local duality constant is
  \[
    \size{\Stab_{\Gamma_{\Lambda_v^{(2)}}^{(2)}}x_{0,v}^{(2)}}
    \size{\OO_v/m\OO_v}.
  \]
  For $v \in \Sigma^{(2)}$, the first local function is \eqref{eq:base_orbit_local_count} and its Fourier transform is
  \[
    \frac{1}{\size{\Stab_{\Gamma_{\Lambda_v^{(1)}}^{(1)}}x_{0,v}^{(1)}}},
  \]
  so the exceptional local duality constant is
  $\size{\Stab_{\Gamma_{\Lambda_v^{(1)}}^{(1)}}x_{0,v}^{(1)}}^{-1}$.
  Taking the inverses of these constants, together with the inverses of the ordinary
  $c_v$, in Theorem \ref{thm:main_compose} proves \eqref{eq:S_reflection}.
\end{proof}

For $N$ a squarefree integer, it is natural to ask what happens if we invert finitely many primes and count binary cubic forms of discriminant $D \neq 0$ over $\ZZ[1/N]$, up to the action of the relevant group $\SL_2(\ZZ[1/N])$. There are still only finitely many for each degree, owing to Hermite's theorem on the finiteness of the number of number fields with prescribed degree and set of ramified primes. However, Cremona's reduction theory \cite{Crem_Redn,Crem_Redn_2} for binary forms over $\ZZ$ does not carry over so easily to this setting. In this section, we prove a reflection theorem that reduces counting forms over $\ZZ[1/N]$ to counting forms over $\ZZ$ with certain conditions at the primes dividing $3N$.

Note that $D$ matters only up to multiplication by the squares in $\ZZ[1/N]^\cross$; hence we can restrict our attention to $D\in \ZZ$ that are \emph{fundamental} at each prime $p \mid N$. (If $p \neq 2$, this means that $p^2 \nmid D$. If $p = 2$, this means that $D \equiv 1 \mod 4$ or $D \equiv 8, 12 \mod 16$. However, we allow $D$ to be non-fundamental at primes not dividing $N$.)

\begin{defn}
Let $N > 1$ be a squarefree integer. Recall that $R_N$ weights a binary cubic form by its number of roots in $\PP^1(\ZZ/N\ZZ)$. Denote by $R_N^\cross$ the weighting by the number of roots that are simple modulo every prime dividing $N$; since $N$ is squarefree,
\[
  R_N^\cross = \prod_{p \mid N} \bigl(3 T_p(111) + T_p(12) + T_p(1^21)\bigr).
\]
Thus $h(D, R_N^\cross)$ is half the number of integer binary cubic forms
\[
  f(x,y) = a x^3 + b x^2 y + c x y^2 + d y^3
\]
such that $N \mid d$ and $\gcd(c, N) = 1$, up to the action of the familiar congruence subgroup
\[
\Gamma^0(N) = \left\{\begin{bmatrix}
  a & Nb \\
  c & d
\end{bmatrix} \in \SL_2 \right\}.
\]
\end{defn}
\begin{rem}
The group $\Gamma^0(N)$ acts freely on nondegenerate forms for $N \geq 2$, so we have elided the mention of stabilizers. The factor of $2$ appears because of transitioning from $\GL_2$ to $\SL_2$.
\end{rem}

\begin{thm}
Let $N$ be a squarefree integer.

For $0 \neq D \in \ZZ[1/N]$, let $h_{\ZZ[1/N]}(D)$ be the number of $\SL_2(\ZZ[1/N])$-orbits of integral binary cubic forms over $\ZZ[1/N]$, each weighted by the reciprocal of its stabilizer in $\SL_2(\ZZ[1/N])$. If $3 \nmid N$, define $h_{3, \ZZ[1/N]}(D)$ to be the same count, counting only $1331$-forms (that is, forms whose middle two coefficients belong to the ideal $3\ZZ[1/N] \subsetneq \ZZ[1/N])$.

Now let $D \in \ZZ$ be a discriminant. As usual, let
\[
  c_{D, \infty} = \begin{cases}
    3 & D > 0 \\
    1 & D < 0.
  \end{cases}
\]
Then:
\begin{enumerate}[$($a$)$]
  \item\label{it:nmb} If $3 \nmid N$, then
  \[
    h_3(-27D, R_N^\cross)
      = \frac{3}{2c_{-D,\infty}} \cdot h_{\ZZ[1/N]}(D).
  \]
  \item\label{it:nma} If $3 \nmid N$, then
  \[
    h_{3, \ZZ[1/N]}(-27D)
      = 2c_{D, \infty} \cdot h(D, R_N^\cross).
  \]
  \item\label{it:3|N} If $3 \mid N$, then
  \[
    h_{\ZZ[1/N]}(-3 D) = 2c_{D,\infty} \cdot h(D, R_N^\cross).
  \]
\end{enumerate}
\end{thm}

\begin{proof}
  Let $\Sigma_N$ be the set of primes dividing $N$. At a prime $p \mid N$, the weighting $R_N^\cross$ can be realized by using the lattice of $111p$-forms and marking the simple root modulo $p$. Because $D$ is fundamental at $p$, the marked forms in the trivial coclass form the integral base orbit, and its stabilizer in the corresponding congruence subgroup $\Gamma^0(p)$ is trivial. Thus the condition at every $p \mid N$ is exactly the base-orbit condition in \eqref{eq:S_class_number}. The resulting $\SL_2$-orbit mass is twice $h(D,R_N^\cross)$, or twice $h_3(D,R_N^\cross)$ in the $3$-traced case, because $h(D,\cdot)$ and $h_3(D,\cdot)$ are defined using $\GL_2(\ZZ)$-orbits.

  On the reflected side we allow denominators at the primes in $\Sigma_N$. Strong approximation for $\SL_2$ gives the class-number-one hypotheses of Theorem \ref{thm:S_reflection} for both the congruence lattice on the first side and the $\Sigma_N$-integral subgroup on the second. We apply Theorem \ref{thm:S_reflection} with $\Sigma^{(1)} = \varnothing$ and $\Sigma^{(2)} = \Sigma_N$, taking the first variety to be the side with the marked-root condition. Since the integral basepoint stabilizer in $\Gamma^0(p)$ is trivial, the exceptional product in \eqref{eq:S_reflection} is $1$.

  If $3 \nmid N$, the remaining local reflection constants are exactly those of Theorem \ref{thm:O-N_traced_intro}. Applying the identity in the reverse direction to the pair of discriminants $D$ and $-27D$ gives
  \[
    h_3(-27D,R_N^\cross)
      = \frac{3}{2c_{-D,\infty}} \cdot h_{\ZZ[1/N]}(D),
  \]
  proving \ref{it:nmb}. Applying it from ordinary forms to $3$-traced forms gives
  \[
    h_{3,\ZZ[1/N]}(-27D)
      = c_{D,\infty} \cdot 2h(D,R_N^\cross),
  \]
  proving \ref{it:nma}.

  Finally suppose that $3 \mid N$. Since $3$ is a unit in $\ZZ[1/N]$, the $1111$- and $1331$-lattices become the same after base change to $\ZZ[1/N]$, and discriminants $-27D$ and $-3D$ differ by the square of a unit. Thus the same application of Theorem \ref{thm:S_reflection} gives
  \[
    h_{\ZZ[1/N]}(-3D)
    = c_{D,\infty} \cdot 2h(D,R_N^\cross),
  \]
  proving \ref{it:3|N}.
\end{proof}

\begin{examp}
Take $N = 15$ and $D = 1$. The unique form $f$ of discriminant $1$ has $3 \cdot 3 = 9$ roots modulo $N$, all simple, but because of the sixfold symmetry of $f$, we get $h(D, R_N^\cross) = 3/2$. Part \ref{it:3|N} of the theorem then assures us that the number of $\SL_2(\ZZ[1/15])$-classes of forms over $\ZZ[1/15]$ of discriminant $-3$ is $2 \cdot 3 \cdot 3/2 = 9$, when forms are weighted by stabilizer. In fact there are nine, none having any stabilizer:
\[
  3x^3 - \frac{1}{3^2}y^3, \quad x^3 \pm \frac{1}{3}y^3, \quad 5 x^3 \pm \frac{1}{3 \cdot 5} y^3, \quad \frac{5}{3} x^3 \pm \frac{1}{5} y^3, \quad 5 \cdot 3 x^3 \pm \frac{1}{5 \cdot 3^2} y^3.
\]
\end{examp}

\bibliography{Master,me}
\bibliographystyle{alpha}

\end{document}